\documentclass[11pt,oneside,reqno]{amsart}
\usepackage{amssymb}
\usepackage{color}
\usepackage{tikz-cd}
\usetikzlibrary{decorations.markings}

\usepackage{xcolor}
\definecolor{deepgreen}{cmyk}{0.99998,0,1,0}

\usepackage{hyperref}
\hypersetup{hypertex=true,
	colorlinks=true,
	linkcolor=blue,
	anchorcolor=blue,
	citecolor=blue}
\usepackage{comment}
\usepackage{mathrsfs}

\allowdisplaybreaks[4]

\usepackage{etoolbox}

\usepackage{bm}
\usepackage{bbm}

\usepackage[new]{old-arrows}

\theoremstyle{definition}
\newtheorem{defi}{$\mathbf{Definition}$}[section]
\newtheorem*{pro}{$\mathbf{Proof}$}

\theoremstyle{plain}
\newtheorem{theo}[defi]{$\mathbf{Theorem}$}
\newtheorem{lemma}[defi]{$\mathbf{Lemma}$}
\newtheorem{coro}[defi]{$\mathbf{Corollary}$}
\newtheorem{remark}[defi]{$\mathbf{Remark}$}
\newtheorem{prop}[defi]{$\mathbf{Proposition}$}

\numberwithin{equation}{section}

\newcommand{\de}{{d}}

\newcommand{\pa}{\partial}

\newcommand{\lk}{\left(}
\newcommand{\rk}{\right)}

\newcommand{\lv}{\left\vert}
\newcommand{\rv}{\right\vert}
\newcommand{\lV}{\left\Vert}
\newcommand{\rV}{\right\Vert}

\newcommand{\blk}{\big(}
\newcommand{\brk}{\big)}
\newcommand{\bli}{\big\langle}
\newcommand{\bri}{\big\rangle}

\newcommand{\bv}{\big\vert}
\newcommand{\Bv}{\Big\vert}
\newcommand{\bbv}{\bigg\vert}

\newcommand{\bV}{\big\Vert}
\newcommand{\BV}{\Big\Vert}

\newcommand{\cref}{\S\ \ref}
\newcommand{\Cref}{\S\ \ref}

\numberwithin{equation}{section}

\newcommand{\supp}{\mathrm{supp}}

\newcommand{\interior}[1]{%
	{\kern0pt#1}^{\mathrm{\,o}}%
}
\definecolor{pink}{RGB}{249,164,186}
\definecolor{grassgreen}{RGB}{128,255,0}

\numberwithin{equation}{section}
\title[]{Small Scale Index Theory, Scalar Curvature, and Gromov's Simplicial norm}

\author{Qiaochu Ma$^{1}$}
\author{Guoliang Yu$^{1}$}
\address{$^1$Department of Mathematics, Texas A\&M University}
\date{}

\makeatletter
\newcommand*{\rom}[1]{{\mathscr{B}_{\varepsilon}}andafter\@slowromancap\romannumeral #1@}
\makeatother

\begin{document} 
	\clearpage
	
	\begin{abstract}
		\fontsize{10pt}{11pt}\selectfont 
		
		In this article, we study the topological complexity of manifolds with a lower scalar curvature bound. We introduce a small scale index theorem to establish an upper bound for Gromov's simplicial norm of the Poincaré dual of the A-hat class for manifolds with spin universal covering, in terms of a scalar curvature lower bound, volume upper bound, and injectivity radius lower bound of the universal covering. This result can be viewed both as a generalization of Lichnerowicz vanishing theorem and as a scalar curvature analogue to Cheeger finiteness theorem.
		
		\normalsize	
	\end{abstract}
	\maketitle
	\section{Introduction}

	\subsection{Backgrounds}\label{A1''}

	The interaction between topology and curvature, including sectional, Ricci, and scalar curvatures, is a major and foundational theme in differential geometry.

	Sectional and Ricci curvatures often exhibit rigidity. For example, the Cheeger finiteness theorem \cite{MR263092} asserts that in a fixed dimension, if the sectional curvature is bounded on both sides, the diameter is bounded above, and the volume is bounded below, then there exist only finitely many diffeomorphism types. This statement can also be reformulated using a lower bound on the injectivity radius and an upper bound on the volume, in place of the diameter upper bound and volume lower bound.

	In comparison, scalar curvature is to some extent the weakest curvature invariant and behaves with flexibility. For instance, Kazdan-Warner theorem \cite{MR375153} states that any smooth function that is somewhere negative can be realized as the scalar curvature of a Riemannian metric. However, the flexibility of scalar curvature is not without restrictions. For example, Lichnerowicz vanishing theorem \cite{MR156292} shows that, if a spin manifold admits a positive scalar curvature metric, then its A-hat genus must vanish. This subtle balance between flexibility and rigidity has been explored in the pioneering works of Schoen-Yau \cite{MR541332,MR526976}, Gromov-Lawson \cite{MR577131,MR569070}, and Witten \cite{MR626707}.

	Both Cheeger finiteness and Lichnerowicz vanishing are about the topological complexity of manifolds. A useful tool for measuring such complexity is Gromov’s simplicial norm \cite{MR686042}, defined as the infimum of the sum of the absolute values of the coefficients in a singular chain representing a given homology class. Intuitively, it quantifies how efficiently a homology class can be represented by a combination of singular simplices. In \cite[\S\,3.13]{MR859578} and \cite[\S\,3.A]{MR4577903}, Gromov conjectured that in a fixed dimension, the simplicial norm of the fundamental class, also known as the simplicial volume, and the Pontryagin numbers of aspherical manifolds can be controlled in terms of a scalar curvature lower bound and a volume upper bound. These conjectures suggest that manifolds with a scalar curvature lower bound exhibit only finite topological complexity.

	In this article, we show that in a fixed dimension, the simplicial norm of the Poincaré dual of the $\widehat{A}$-class of a manifold with spin universal covering can be controlled in terms of a scalar curvature lower bound, injectivity radius of the universal covering, and volume. This theorem can be regarded both as a generalization of Lichnerowicz vanishing theorem and as a scalar curvature analogue of Cheeger finiteness theorem.
	
	
	We now explain in detail.

	\subsection{Simplicial norm estimate}

	Let $(M, g^{TM})$ be an $m$-dimensional closed oriented Riemannian manifold and $R^{TM}$ the associated Riemannian curvature. The $\widehat{A}$-class $\widehat{A}(TM)\in H^*(M)$ of $M$ is a polynomial of the Pontryagin class $p(TM)$ of $M$ and can be represented by
	\begin{equation}
		\begin{split}
			\widehat{A}(TM)
			&=\det{}^{1/2}\Bigg(\frac{\tfrac{1}{2\pi i}R^{TM}/2}{\sinh\Big(\tfrac{1}{2\pi i}R^{TM}/2\Big)}\Bigg)\\
			&=\widehat{A}^{[0]}(TM)+\widehat{A}^{[4]}(TM)+\widehat{A}^{[8]}(TM)+\cdots\\
			&=1-\frac{1}{24}p_1(TM)+\Big(-\frac{1}{1440}p_2(TM)+\frac{7}{5760}p_1(TM)^2\Big)+\cdots
		\end{split}
	\end{equation}
	where $\widehat{A}^{[i]}(TM)\in H^{i}(M)$ is the $i$-th degree component of $\widehat{A}(TM)$ and $p_i(TM)$ is the $i$-th Pontryagin class of $M$. The $\widehat{A}$-genus $\widehat{A}(M)$ of $M$ is the integral of $\widehat{A}$-class,
	\begin{equation}
		\widehat{A}(M)=\int_M\widehat{A}(TM)=\int_M\widehat{A}^{[m]}(TM).
	\end{equation}

	Let $[M]\in H_m(M)$ be the fundamental homology class of $M$. Then for any $i\in\mathbb{N}$, the cap product $\widehat{A}^{[i]}(TM)\cap [M]\in H_{m-i}(M)$ defines the Poincaré dual of $\widehat{A}^{[i]}(TM)$. Note that in particular, we have $\widehat{A}^{[m]}(TM)\cap [M]=\widehat{A}(M)$ and $\widehat{A}^{[0]}(TM)\cap [M]=[M]$. Now we state our main result.

	\begin{theo}\label{t1.1111}
		For any integers $m\in \mathbb{N}$, there exists a constant $C_m>0$ such that the following holds. Let $(M,g^{TM})$ be an $m$-dimensional closed oriented Riemannian manifold with its universal covering $\widetilde{M}$ spin. If
		\begin{equation}
			\inf_{x\in M} \mathrm{Sc}_{g^{TM}}(x)\geqslant k_0,\ \ \mathrm{Inj}(\widetilde{M})\geqslant i_0,\ \ \sup_{\widetilde{x}\in \widetilde{M}}\mathrm{Vol}\big(B^{\widetilde{M}}(\widetilde{x},i_0)\big)\leqslant V_0, \ \ \mathrm{Vol}(M)\leqslant V_1,
		\end{equation}
		then the simplicial norm $\lV\cdot\rV_{\ell^1}$ of $\widehat{A}^{[i]}(TM)\cap [M]$ satisfies
		\begin{equation}\label{1.11111}
			\lV\widehat{A}^{[i]}(TM)\cap[M]\rV_{\ell^1}\leqslant C_m\big(\lv k_0\rv^{\tfrac{m}{2}}+i_0^{-m}\big)^{m-i+1}V_0^{m-i}V_1.
		\end{equation}
		Here $k_0 \in \mathbb{R}$, $i_0,V_0,V_1>0$, $\mathrm{Sc}_{g^{TM}}(x)$ denotes scalar curvature, $\mathrm{Vol}(\cdot)$ the Riemannian volume, $\mathrm{Inj}(\widetilde{M})$ the injectivity radius of $\widetilde{M}$, and $B^{\widetilde{M}}(\widetilde{x},i_0)$ the geodesic ball in $\widetilde{M}$ of radius $i_0$ centered at $\widetilde{x}$, 
	\end{theo}

	\subsection{Small scale index theorem}

	Our simplicial norm estimate is derived from a small scale index theorem for Dirac operators, itself inspired by three primary sources. First, Connes’ noncommutative geometry \cite{MR823176} provides the conceptual framework for analyzing the involved geometric structures through spectral triples. Second, the quantitative operator K-theory developed in \cite{MR1626745,MR3590530} builds a bridge between large scale topology and microlocal geometry. Third, Bismut’s local index theorem \cite{MR813584} introduces powerful asymptotic tools that support our analytic arguments. The interplay of these ideas gives our result special interest and suggests a variety of potential applications.

	The small scale index theorem can be regarded as a quantitative refinement of Connes-Moscovici higher index theorem \cite{MR1066176}, and the proof depends heavily on the analytic localization techniques of Bismut-Lebeau \cite{MR1188532}. The small scale index of the Dirac operator captures information about the Poincaré dual $\widehat{A}(TM)\cap[M]$, of the $\widehat{A}$-class, and lives in a quantitative homology group 
	$H_*^\varepsilon(M)$ called epsilon homology, in which simplices intuitively have diameters less than or equal to epsilon. The small scale index theorem allows us to find a representative $\mathrm{Ch}(\mathrm{Ind}_\varepsilon(D^{S_{\widetilde{M}}}))\in H_*^\varepsilon(M)$ of $\widehat{A}(TM)\cap[M]$ in terms of the kernels of certain operators with small propagation, constructed via functional calculus of the Dirac operator $D^{S_{\widetilde{M}}}$. Using Moser iteration \cite{MR159138} and Croke's isoperimetric constant lower bound	\cite{MR608287} of a region inside a ball whose radius is smaller than the injectivity radius, these kernels can be estimated using a scalar curvature lower bound and injectivity radius. By combining these ingredients, we can obtain an upper bound for the simplicial norm of $\widehat{A}(TM)\cap[M]$.

	Connes' noncommutative geometry and the analytic techniques of Bismut-Lebeau have also inspired Zhang's work \cite{MR3664818} on the generalization of the vanishing theorems of Lichnerowicz \cite{MR156292} and Hitchin \cite{MR358873} to foliations, by combining Connes’ fibration \cite{MR866491} with the sub-Dirac operator of Liu-Zhang \cite{MR1850748}.

	Note that in this work, we focus on the even dimensional case. The odd dimensional case can be treated using standard techniques involving a suspension argument.

	\subsection{A $C^0$-stable bound}

	In our small-scale index theorem, different choices of the scale parameter $\varepsilon$ yield different estimates for the simplicial norm.
	
	\begin{theo}\label{tt1.1}
		For any integers $m\in \mathbb{N}$, there exists a constant $C_m>0$ such that the following holds. Let $(M,g^{TM})$ be an $m$-dimensional closed oriented aspherical Riemannian manifold, that is, its universal covering $\widetilde{M}$ is contractible. If
		\begin{equation}
			\sup_{x\in M} \mathrm{Sc}_{g^{TM}}(x)\geqslant k_0,\ \ \inf_{x\in\widetilde{M}}\mathrm{ID}\big(B^{\widetilde{M}}(x,1)\big)\geqslant D, \ \  \sup_{\widetilde{x}\in \widetilde{M}}\mathrm{Vol}\big(B^{\widetilde{M}}(\widetilde{x},1)\big)\leqslant V_0,\ \ \mathrm{Vol}(M)\leqslant V_1, 
		\end{equation}
		then we have
		\begin{equation}\label{1.4.'}
			\lV\widehat{A}^{[i]}(TM)\cap[M]\rV_{\ell^1}\leqslant C_m\bigg(\frac{\lv k_0\rv^{\tfrac{m}{2}}+1}{D^m}\bigg)^{m-i+1}V_0^{m-i}V_1.
		\end{equation}
		Here $k_0 \in \mathbb{R}$ and $D,V_0,V_1>0$, $B^{\widetilde{M}}(x,1)$ is the unit ball with center $x\in\widetilde{M}$ and $\mathrm{ID}\big(B^{\widetilde{M}}(x,1)\big)$ is the Dirichlet isoperimetric constant given by
		\begin{equation}\label{2.5'}
			\mathrm{ID}\big(B^{\widetilde{M}}(x,1)\big)=\inf_{\begin{subarray}{c}
					\Sigma\subset B^{\widetilde{M}}(x,1),\\ \partial \Sigma\cap \partial B^{\widetilde{M}}(x,1)=\emptyset
			\end{subarray}}\frac{\mathrm{Area}(\partial\Sigma)}{\mathrm{Vol}(\Sigma)^{(m-1)/m}},
		\end{equation}
		where the infimum is taken over all subdomains $\Sigma\subset B^{\widetilde{M}}(x,1)$ with the property that $\partial\Sigma$ is a hypersurface not intersecting $\partial B^{\widetilde{M}}(x,1)$.
		
	\end{theo}

	Compared with \eqref{1.11111}, this bound depends also on the volume $\mathrm{Vol}\big(B^{\widetilde{M}}(x,1)\big)$ of radius one balls $B^{\widetilde{M}}(x,1)$ in $\widetilde{M}$ and their Dirichlet isoperimetric constant $\mathrm{ID}\big(B^{\widetilde{M}}(x,1)\big)$. We emphasize that this estimate is essentially $C^0$ in nature. To see this, recall Gromov’s result \cite[\S\,1.8]{MR3201312} that a uniform lower bound on scalar curvature is a $C^0$-property.
	
	In earlier work with Wang and Zhu \cite{ma2024boundingahatgenususing}, we proved that the $\widehat{A}$-genus $\widehat{A}(M)$ of a closed spin manifold can be controlled in terms of the scalar curvature lower bound $\sup_{x\in M}\mathrm{Sc}^-_{g^{TM}}(x)$, Neumann isoperimetric constant $\mathrm{IN}(M)$, and $\mathrm{Vol}(M)$. Examples constructed there showed that the isoperimetric constant is necessary in such estimates.

	
	The volume $\mathrm{Vol}(B^{\widetilde{M}}(x,1))$ of a radius-one ball in $\widetilde{M}$ is called the macroscopic scalar curvature. For further background, we refer to Guth’s survey \cite{MR2827817}. Inspiring results in this direction were made by Guth \cite{MR2753599} in the hyperbolic case and Braun-Sauer \cite{MR4386411} in the general case, where they proved that the simplicial volume $\lV[M]\rV_{\ell^1}$ can be controlled in terms of $\sup_{x\in\widetilde{M}}\mathrm{Vol}\big(B^{\widetilde{M}}(x,1)\big)$ and $\mathrm{Vol}(M)$.

	\subsection{Organization of the article}

	The structure of this paper is as follows. In \cref{SQH}, we introduce a quantitative homology theory called epsilon homology. In \cref{SSS}, we present the small scale index theorem. In \cref{SNE}, we apply the small scale index theorem to control simplicial norms. In \cref{SCCC}, we construct the Connes-Chern character of the Dirac operator within the framework of epsilon homology. In \cref{SLIC}, we provide a proof of the small-scale index theorem, which computes the Connes-Chern character of the Dirac operator in epsilon homology.
	In \cref{FAM}, we validate some asymptotic expansions used in \cref{SLIC}. In \cref{SEE}, we estimate the kernels of small propagation operators constructed from the Dirac operators using functional calculus in terms of a scalar curvature lower bound and the Dirichlet isoperimetric constant.

	\subsection{Acknowledgment}

	We are grateful to Zhizhang Xie, Ruobin Zhang, and Bo Zhu for helpful discussions. We would like to thank Xianzhe Dai, Jinmin Wang, Guofang Wei for	insightful suggestions, particularly for pointing us to Grove-Petersen-Wu \cite{MR1029396} and Croke \cite{MR608287}. Qiaochu Ma is partially supported by NSF grants DMS-1546917 and DMS-2247313, and Guoliang Yu is partially supported by NSF grant DMS-2247322.

	\section{Quantitative homology}\label{SQH}

	In this section, we discuss Alexander-Spanier homology and its quantitative version, called epsilon homology. In \cref{QH1}, we review the Alexander-Spanier homology following Connes-Moscovici \cite[\S\,1]{MR1066176} and Moscovici-Wu \cite[\S\,1]{MR1292020}. In \cref{QH2}, we introduce the epsilon homology following the algebraic formulation given in \cite[\S\,4]{MR3590530}. In \cref{QH3}, we define a quantity called the homological radius, which plays a key role in the convergence of epsilon homology.

	We shall use the notation that $M$ is a compact manifold without boundary, $\Gamma=\pi_1(M)$ is its fundamental group and $\widetilde{M}$ is its universal covering.
	\subsection{Alexander-Spanier (co)homology}\label{QH1}

	Alexander-Spanier homology can be thought of as describing the homology of a manifold using extremely small simplices. When the simplices are made sufficiently small, the manifold is decomposed into a very fine net, from which we can recover the full homological information.

	For any finite open covering $\mathscr{U}$ of $M$, let $C^k_\mathscr{U}(M)$ be the space of continuous functions on the subspace $\cup_{U\in\mathscr{U}}U^{k+1}\subseteq M^{k+1}$. The coboundary map $\pa\colon C^k_\mathscr{U}(M)\to C^{k+1}_\mathscr{U}(M)$ is defined by
	\begin{equation}\label{2.1}
		(\pa f)(x_0,\cdots,x_k)=\sum_{i=0}^{k}(-1)^if(x_0,\cdots,\widehat{x}_i,\cdots,x_k),
	\end{equation}
	where the ‘hat’ symbol over $x_i$ indicates that this variable is deleted.
	
	Let $H^*_\mathscr{U}(M)$ denote the cohomology of the cochain complex $(\pa,C^*_\mathscr{U}(M))$. Using the natural cochain homomorphism $H^*_\mathscr{U}(M)\to H^*_\mathscr{V}(M)$ of restriction given by the inclusion of a refinement $\mathscr{V}\to \mathscr{U}$, we see that $\{H^*_\mathscr{U}(M)\}_{\mathscr{U}}$ forms a direct system and the Alexander-Spanier cohomology $H_{\mathrm{AS}}^{*}(M)$ of $M$ is defined by a direct limit
	\begin{equation}\label{2.2}
		H_{\mathrm{AS}}^{*}(M)=\varinjlim \{H^*_\mathscr{U}(M)\}_{\mathscr{U}}.
	\end{equation}

	Let $C_k^\mathscr{U}(M)$ be the space of Borel measures on $M^{k+1}$ with their supports contained in the subset $\cup_{U\in\mathscr{U}}U^{k+1}$. The boundary map $\pa\colon C_{k}^\mathscr{U}(M)\to C_{k-1}^\mathscr{U}(M)$ is dually by
	\begin{equation}\label{2.3'}
		\langle\pa\mu,f\rangle=\langle\mu,\pa f\rangle
	\end{equation}
	for any $\mu\in C_k^\mathscr{U}(M)$ and $f\in C^k_\mathscr{U}(M)$, where $\langle\cdot,\cdot\rangle$ denotes the natural pairing between measures and continuous functions. We can also write the boundary map explicitly by
	\begin{equation}\label{2.4}
		(\pa\mu)(x_0,\cdots,x_k)=\sum_{i=0}^{k}(-1)^i\int_M\mu(x_0,\cdots,x_{i-1},x,x_{i+1}\cdots,x_k),
	\end{equation}
	a formal alternating sum of the integration over the $i$-th variable.

	Let $H_*^\mathscr{U}(M)$ denote the homology of the chain complex $(\pa,C_*^\mathscr{U}(M))$. Using the natural chain homomorphism $H_*^\mathscr{V}(M)\to H_*^\mathscr{U}(M)$ of extension given by the inclusion of a refinement $\mathscr{V}\to \mathscr{U}$, we see that $\{H_*^\mathscr{U}(M)\}_{\mathscr{U}}$ forms a inverse system and the Alexander-Spanier homology $H^{\mathrm{AS}}_{*}(M)$ of $M$ is defined by a inverse limit
	\begin{equation}\label{2.5}
		H^{\mathrm{AS}}_{*}(M)=\varprojlim \{H_*^\mathscr{U}(M)\}_{\mathscr{U}}.
	\end{equation}

	\subsection{Epsilon (co)homology}\label{QH2}

	Alexander-Spanier homology involves a limit process and relies only on the vertices of a simplex rather than the entire simplex. This motivates two adjustments, introducing a parameter to analyze convergence quantitatively, and lifting to the universal cover so that simplices with the same vertices but not homologous can be distinguished.

	For any $(x_0,x_1,\cdots,x_k)\in \widetilde{M}^{k+1}$, we set
	\begin{equation}
		\mathrm{diam}(x_0,x_1,\cdots,x_k)=\max_{0\leqslant i,j\leqslant k}d_{\widetilde{M}}(x_i,x_j).
	\end{equation}

	For any $\varepsilon>0$, let $C^{k}_\varepsilon(\widetilde{M})$ be the space of continuous functions defined on the $\varepsilon$-diagonal $\mathrm{diag}_\varepsilon(\widetilde{M}^{k+1})\subseteq \widetilde{M}^{k+1}$ given by
	\begin{equation}
		\mathrm{diag}_\varepsilon(\widetilde{M}^{k+1})=\{(x_0,\cdots,x_k)\in \widetilde{M}^{k+1}\mid \mathrm{diam}(x_0,x_1,\cdots,x_k)\leqslant\varepsilon\}.
	\end{equation}
	Let $C^{k}_\varepsilon(\widetilde{M},\Gamma)\subset C^{k}_\varepsilon(\widetilde{M})$ be the subspace of $\Gamma$-invariant functions, where $\Gamma$ acts diagonally on $\widetilde{M}^{k+1}$. The coboundary map \eqref{2.1} restricts to a coboundary map on $C^{k}_\varepsilon(\widetilde{M},\Gamma)$, hence we obtain a cochain complex $(\pa,C^{k}_\varepsilon(\widetilde{M},\Gamma))$.
	
	\begin{defi}
		We define the epsilon cohomology  $H^k_\varepsilon(\widetilde{M},\Gamma)$ to be the cohomology of the cochain complex $(\pa,C^{k}_\varepsilon(\widetilde{M},\Gamma))$.
		
	\end{defi}
	
	Similar to \eqref{2.2}, we define
	\begin{equation}\label{2.8..}
		H_{\mathrm{AS}}^{*}(\widetilde{M},\Gamma)=\varinjlim \{H^*_\varepsilon(\widetilde{M},\Gamma)\}_{\varepsilon>0}.
	\end{equation}

	Likewise, let $C^{*}_\varepsilon(\widetilde{M})$ be the space of measures supported on the $\varepsilon$-diagonal $\mathrm{diag}_\varepsilon(\widetilde{M}^{k+1})\subseteq \widetilde{M}^{k+1}$, then from the subspace $C_{*}^\varepsilon(\widetilde{M},\Gamma)\subset C^{*}_\varepsilon(\widetilde{M})$ of $\Gamma$-invariant measures and the boundary map \eqref{2.3'}, we obtain a chain complex $(\pa,C_{*}^\varepsilon(\widetilde{M},\Gamma))$.
	
	\begin{defi} We define the epsilon homology
		$H_*^\varepsilon(\widetilde{M},\Gamma)$ to be the homology of the chain complex $(\pa,C_{*}^\varepsilon(\widetilde{M},\Gamma))$. 
	\end{defi}
	
	Similar to \eqref{2.5}, define
	\begin{equation}\label{2.9..}
		H^{\mathrm{AS}}_{*}(\widetilde{M},\Gamma)=\varprojlim \{H_*^\varepsilon(\widetilde{M},\Gamma)\}_{\varepsilon>0}.
	\end{equation}

	By viewing $f\in H^*_\varepsilon(\widetilde{M},\Gamma)$ a continuous function on $\Gamma\backslash \widetilde{M}^{*+1}$ and $\mu\in H^*_\varepsilon(\widetilde{M},\Gamma)$ a measure on $\Gamma\backslash \widetilde{M}^{*+1}$, we can define a natural pairing $\langle\cdot,\cdot\rangle_\Gamma$ between $H^*_\varepsilon(\widetilde{M},\Gamma)$ and $H_*^\varepsilon(\widetilde{M},\Gamma)$ by the integral of $f$ with respect to $\mu$ over $\Gamma\backslash \widetilde{M}^{*+1}$. In other words, let $D\subset \widetilde{M}$ be a fundamental domain of $\Gamma$-action, then $D\times \widetilde{M}^k$ is a fundamental domain of the diagonal $\Gamma$-action on $\widetilde{M}^{k+1}$, and
	\begin{equation}\label{2.10''}
		\langle f,\mu\rangle_\Gamma=\int_{D\times \widetilde{M}^k}f\cdot d\mu.
	\end{equation}
	Note that $\mathrm{diag}_\varepsilon(\widetilde{M}^{k+1})\cap \big(D\times \widetilde{M}^k\big)$ is compact, so the integral \eqref{2.10''} is finite. Since the $\Gamma$-action on $C_{*}^\varepsilon(\widetilde{M})$ and $C^{*}_\varepsilon(\widetilde{M})$ commutes with the (co)boundary operator $\pa$, we see that
	\begin{equation}
		\langle \pa f,\mu\rangle_\Gamma=\langle f,\pa\mu\rangle_\Gamma.
	\end{equation}
	By \eqref{2.8..} and \eqref{2.9..}, this pairing also induces a pairing between $H^*_{\mathrm{AS}}(\widetilde{M},\Gamma)$ and $H_*^{\mathrm{AS}}(\widetilde{M},\Gamma)$, and we still denote it by $\langle\cdot,\cdot\rangle_\Gamma$.

	Now we compare $H^*_{\mathrm{AS}}(\widetilde{M},\Gamma)$ and $H_*^{\mathrm{AS}}(\widetilde{M},\Gamma)$ with $H^*_{\mathrm{AS}}(M)$ and $H_*^{\mathrm{AS}}(M)$. Using pullback functions through $\widetilde{M}\to M$, we get a morphism
	\begin{equation}\label{2.12}
		C^{*}_{\varepsilon}(M)\to C^{*}_{\varepsilon}(\widetilde{M},\Gamma).
	\end{equation}
	Since $\Gamma$ acts discretely on $\widetilde{M}$, when $\varepsilon$ is sufficiently small, the projection map $\Gamma\backslash \mathrm{diag}_\varepsilon(\widetilde{M}^{k+1})\to \mathrm{diag}_\varepsilon(M^{k+1})$ is bijective. In this case, \eqref{2.12} becomes an isomorphism, and therefore
	\begin{equation}\label{2.13}
		H^*_{\mathrm{AS}}(M)\cong H^*_{\mathrm{AS}}(\widetilde{M},\Gamma).
	\end{equation}
	Similarly, by pushforward measures, we have a morphism
	\begin{equation}\label{2.14}
		C_{*}^{\varepsilon}(\widetilde{M},\Gamma)\to C_{*}^{\varepsilon}(M),
	\end{equation}
	which becomes isomorphism when $\varepsilon$ small, and 
	\begin{equation}\label{2.15}
		H_*^{\mathrm{AS}}(M)\cong H_*^{\mathrm{AS}}(\widetilde{M},\Gamma).
	\end{equation}

	\subsection{Homological radius}\label{QH3}

	From \eqref{2.13} and \eqref{2.15}, we see that, as $\varepsilon\to0$, epsilon (co)homology coincides with Alexander-Spanier (co)homology. We now analyze this convergence process in detail.

	Let $C_{*}(M)$ denote the singular complex of $M$ and $C_*(\widetilde{M},\Gamma)$ the $\Gamma$-invariant singular complex of $\widetilde{M}$, then from the path lifting property \cite[Proposition 1.30]{MR1867354} we get
	\begin{equation}\label{2.16}
		C_{*}(M)\cong C_*(\widetilde{M},\Gamma).
	\end{equation}
	
	For any $\varepsilon>0$, we can construct a chain map
	\begin{equation}\label{2.17}
		\begin{split}
			\chi_\varepsilon\colon C_*(\widetilde{M},\Gamma)&\to C_*^\varepsilon (\widetilde{M},\Gamma)\\
			\Delta_{(x_0,\cdots,x_k)}&\mapsto \delta_{(x_0,\cdots,x_k)},
		\end{split}
	\end{equation}
	where $\Delta_{(x_0,\cdots,x_k)}$ is a $k$-simplex having vertices $(x_0,\cdots,x_k)\in\Gamma\backslash \widetilde{M}^{k+1}$ with $\mathrm{diam}(x_0,\cdots,x_k)\leqslant\varepsilon$, and $\delta_{(x_0,\cdots,x_k)}$ is the delta measure on $(x_0,\cdots,x_k)$, or after sufficiently many barycentric subdivisions, $\Delta_{(x_0,\dots,x_k)}$ can be written as a sum of simplices whose vertices have diameter at most $\varepsilon$, then we send each such simplex to the corresponding delta measure.

	Now we try to construct the inverse of \eqref{2.17}
	\begin{equation}\label{2.18}
		\chi^{-1}_\varepsilon\colon H_*^\varepsilon (\widetilde{M},\Gamma)\rightarrow H_*(\widetilde{M},\Gamma)
	\end{equation}
	at the level of homology. Similar to the discussion in Gromov \cite[\S\,3.3]{MR686042}, we should assign to each $(x_0,\cdots,x_k)\in \Gamma\backslash\widetilde{M}^{k+1}$ with $\mathrm{diam}(x_0,\cdots,x_k)\leqslant\varepsilon$ a singular simplex $\Delta_{(x_0,\cdots,x_k)}$ having vertices $(x_0,\cdots,x_k)$, such that every subset of $(x_0,\cdots,x_k)$ go under this assignment to the corresponding face of $\Delta_{(x_0,\cdots,x_k)}$. Since $\widetilde{M}$ is connected, the construction is straightforward for $k=1$. For $k\geqslant 2$, we proceed inductively using homotopies between lower dimensional simplices. To ensure the existence of these homotopies, it is sufficient that each ball $B^{\widetilde{M}}(x,\varepsilon)$ is contained in a contractible subset of $\widetilde{M}$. This leads us to the following definition.

	\begin{defi}\label{def2.3}
		We define the homological radius $r(M)$ of $M$ by the supremum of all $\varepsilon>0$ with the following property, for any $x\in\widetilde{M}$, there is a contractible open set $U_x\subseteq \widetilde{M}$ such that $B^{\widetilde{M}}(x,\varepsilon)\subseteq U_x$.
	\end{defi}

	Let $(\Omega^*(M),d)$ be the de Rham cochain complex of differential forms and exterior derivative and $H^*_{\mathrm{dR}}(M)$ its cohomology. When $0<\varepsilon< r(M)$, we can define \eqref{2.18} implicitly by
	\begin{equation}\label{2.19'}
		\langle \chi^{-1}_\varepsilon\mu,\alpha\rangle=\int_{\Gamma\backslash\widetilde{M}^{k+1}}\Big(\int_{\Delta_{(x_0,\cdots,x_k)}}\alpha\Big)\mu(x_0,\cdots,x_k)
	\end{equation}
	for any $\alpha\in H^*_{\mathrm{dR}}(M)$. To define the integral, the simplices $\Delta_{(x_0,\cdots,x_k)}$ should be constructed in $C^1$-regularity. This can always be done, since in Definition \ref{def2.3} we assume that the contractible set $U_x$ is open. We now summarize the above discussion in the following main result of this section.

	\begin{prop}\label{prop2.4}
		For any $0<\varepsilon<r(M)$, we have
		\begin{equation}\label{2.20}
			H_*^\varepsilon\big(\widetilde{M},\Gamma\big)\cong H_*^{\mathrm{AS}}\big(\widetilde{M},\Gamma\big).
		\end{equation}
	\end{prop}
	Moreover, we clearly have the following result for $r(M)$.
	\begin{coro}\label{Coro2.5}
		We have $r(M)> \mathrm{Inj}(\widetilde{M})$. Also, if $M$ is aspherical, that is when $\widetilde{M}$ is contractible, then $r(M)=\infty$.
	\end{coro}

	Note that all \eqref{2.16}, \eqref{2.17}, \eqref{2.18}, \eqref{2.19'}, and \eqref{2.20} have a cohomology courtpart. In particular, we can define a cochain map dual to \eqref{2.17}, still denoted by $\chi_\varepsilon$. Replacing continuous functions in $C^{*}_{\varepsilon}(\widetilde{M},\Gamma)$ with smooth functions makes no essential difference at the level of cohomology. Following \cite[Lemma 1.5]{MR1066176}, the map $\chi_\varepsilon$ is given by
	\begin{equation}\label{2.19}
		\begin{split}
			\chi_\varepsilon\colon C^{*}_{\varepsilon}(\widetilde{M},\Gamma)&\to \Omega^*(M)\\
			f&\mapsto \chi_\varepsilon(f)_{x}(v_1,\cdots,v_k)=\frac{1}{k!}\sum_{\sigma\in\mathfrak{G}_k}\mathrm{sgn}(\sigma)(v_{1,x_{\sigma(1)}}\cdots v_{k,x_{\sigma(k)}}f)(x,\cdots,x),
		\end{split}
	\end{equation}
	where $v_{i,x_{\sigma(i)}}$ indicates that the vector $v_i\in T_{x}M$ acts on the $\sigma(i)$-th coordinate of $f$.

	\section{Small scale index theorem and simplicial norms}\label{SSS}
	
	In this section, we formulate a small scale index theorem for manifolds whose universal coverings are spin. In \cref{SS1}, we explain how the spin condition can be weakened, following Rosenberg \cite[\S,3B]{MR720934}. In \cref{SS2}, we introduce the small scale index of the Dirac operator represented by the difference of two idempotents with small propagation following Roe \cite[\S\,7]{MR996446}. In \cref{SS3}, we state our small scale index theorem, which can be viewed as a quantitative version of the higher index theorem of Connes-Moscovici \cite[Theorems 3.9, 5.4]{MR1066176}. Our approach also relies on ideas and techniques developed by Connes \cite{MR823176}, Moscovici–Wu \cite{MR1254310}, Roe \cite{MR1147350}, and Willett-Yu \cite{MR4411373}.


	\subsection{Weakened spin condition}\label{SS1}

	In Rosenberg \cite[\S\,3B)]{MR720934}, a version of Atiyah's $L^2$-index theorem was established under the spin assumption on the universal covering. We shall use the similar trick.

	Now we suppose that $\widetilde{M}$ is spin, while $M$ itself is not necessarily spin. Let $P_{\widetilde{M}}^{x\mathrm{SO}_m}$ be the principal bundle of orthonormal oriented frames on $\widetilde{M}$ and $P_{\widetilde{M}}^{\mathrm{Spin}_m}$ the principal spin bundle. A key observation is that, by the path lifting property, we can define a two covering of $\Gamma$, denoted by $\widetilde{\Gamma}$, using the pullback diagram
	\begin{equation}\label{fig1}
		\begin{tikzcd}[ampersand replacement=\&, column sep=normal,row sep=normal,background color=white!20]
			\widetilde{\Gamma}\arrow[r,""swap]\arrow[d,""swap]\& \mathrm{Aut}(P_{\widetilde{M}}^{\mathrm{Spin}_m})\arrow[d,""swap]\\
			\Gamma \arrow[r,""swap]\& \mathrm{Aut}(P_{\widetilde{M}}^{\mathrm{SO}_m})
		\end{tikzcd}
	\end{equation}
	then $\widetilde{\Gamma}$ acts isometrically on $S_{\widetilde{M}}$ as well as $L^2\big(\widetilde{M},S_{\widetilde{M}}\big)$.

	Let $\mathscr{K}^{\widetilde{\Gamma}}$ be the space of $\widetilde{\Gamma}$-invariant smooth kernel operators on $L^2\big(\widetilde{M},S_{\widetilde{M}}\big)$, the element of which takes the form of $K=K(x,y)\in S_{\widetilde{M},x}\otimes S_{\widetilde{M},y}^*$ with
	\begin{equation}\label{3.2''}
		K(x,y)=\widetilde{\gamma}K\big(\widetilde{\gamma}^{-1} x,\widetilde{\gamma}^{-1} y\big)\widetilde{\gamma}^{-1}
	\end{equation} 
	for any $\widetilde{\gamma}\in \widetilde{\Gamma}$. Although $K$ is not $\gamma$-invariant, its pointwise trace is $\gamma$-invariant in the sense that
	\begin{equation}
		\mathrm{Tr}^{S_{\widetilde{M}}}\big[K(x,x)\big]=\mathrm{Tr}^{S_{\widetilde{M}}}\big[K(\gamma x,\gamma x)\big]
	\end{equation} 
	for any $\gamma\in\Gamma$, due to the fact that $\widetilde{\Gamma}$ is a two-covering of $\Gamma$.

	Let $D^{S_{\widetilde{M}}}$ be the Dirac operator on $\widetilde{M}$, then $D^{S_{\widetilde{M}}}$ is $\widetilde{\Gamma}$-invariant, then for any $\phi\in\mathscr{S}(\mathbb{R})$, the functional calculus $\phi(D^{S_{\widetilde{M}}})$ is a smooth kernel operator and $\phi(D^{S_{\widetilde{M}}})\in\mathscr{K}^{\widetilde{\Gamma}}$.

	\subsection{Finite propagation idempotents}\label{SS2}

	Let $\cos(tD^{S_{\widetilde{M}}})(x,y)$ be the wave kernel of $D^{S_{\widetilde{M}}}$. Let $\Omega\subset \widetilde{M}$ be an open bounded domain and $\cos(tD^{S_{\Omega}})(x,y)$ be the wave kernel of $D^{S_{\Omega}}$. Then we have the following finite propagation speed property of wave equations.
	
	\begin{prop}
		For any $x\in\Omega$, if $\lv t\rv\leqslant d(y,\pa\Omega)$, then we have
		\begin{equation}\label{3.4}
			\cos(tD^{S_{\widetilde{M}}})(x,y)=\cos(tD^{S_{\Omega}})(x,y).
		\end{equation}
	\end{prop}

	Now we follow the construction of Roe \cite[Lemma 7.5]{MR996446}. Choose an odd smooth function $f(\lambda)$ such that $f(-\infty)=-1,f(\infty)=1$, $f'(\lambda)$ is even, $f'(t)\in\mathscr{S}(\mathbb{R})$ and $\mathrm{supp}(\widehat{f'}(\lambda))\subseteq [-\tfrac{1}{6m},\tfrac{1}{6m}]$, here $\widehat{f'}$ 
	is the Fourier transform of $f'$. We then define an even function $g(\lambda)\in \mathscr{S}(\mathbb{R})$ with $\mathrm{supp}(\widehat{g}(\lambda))\subseteq [-\tfrac{1}{3m},\tfrac{1}{3m}]$ by
	\begin{equation}\label{3.5'}
		g(\lambda)=f(\lambda)^2-1.
	\end{equation}

	Let $L_\varepsilon$ be an invertible operator acting on $L^2(\widetilde{M},S_{\widetilde{M}})$ with its inverse $L_\varepsilon^{-1}$ given by
	\begin{equation}
		\begin{split}
			L_\varepsilon&=\begin{pmatrix} -g^+(\varepsilon D^{S_{\widetilde{M}}})&(g^+(\varepsilon D^{S_{\widetilde{M}}})-1)f^-(\varepsilon D^{S_{\widetilde{M}}})\\
				f^+(\varepsilon D^{S_{\widetilde{M}}})&-g^-(\varepsilon D^{S_{\widetilde{M}}})
			\end{pmatrix},\\
			L_\varepsilon^{-1}&=\begin{pmatrix} -g^+(\varepsilon D^{S_{\widetilde{M}}})&-(g^+(\varepsilon D^{S_{\widetilde{M}}})-1)f^-(\varepsilon D^{S_{\widetilde{M}}})\\
				-f^+(\varepsilon D^{S_{\widetilde{M}}})&-g^-(\varepsilon D^{S_{\widetilde{M}}})
			\end{pmatrix},
		\end{split}
	\end{equation}
	and let $P_\varepsilon$ be the corresponding idempotent defined by
	\begin{equation}
		\begin{split}
			P_\varepsilon&=L_\varepsilon\begin{pmatrix} 1&0\\
				0&0
			\end{pmatrix}L_\varepsilon^{-1}\\
			&=\begin{pmatrix} g^+(\varepsilon D^{S_{\widetilde{M}}})^2&-g^+(\varepsilon D^{S_{\widetilde{M}}})\big(1-g^+(\varepsilon D^{S_{\widetilde{M}}})\big)f^-(\varepsilon D^{S_{\widetilde{M}}})\\
				-g^-(\varepsilon D^{S_{\widetilde{M}}})f^+(\varepsilon D^{S_{\widetilde{M}}})&1-g^-(\varepsilon D^{S_{\widetilde{M}}})^2
			\end{pmatrix},
		\end{split}
	\end{equation}
	and $I_\varepsilon$ the ``difference idempotent" defined by
	\begin{equation}\label{3.8'''}
		\begin{split}
			I_\varepsilon&=P_\varepsilon-\begin{pmatrix} 0&0\\
				0&1
			\end{pmatrix}\\
			&=\begin{pmatrix} g^+(\varepsilon D^{S_{\widetilde{M}}})^2&-g^+(\varepsilon D^{S_{\widetilde{M}}})\big(1-g^+(\varepsilon D^{S_{\widetilde{M}}})\big)f^-(\varepsilon D^{S_{\widetilde{M}}})\\
				-g^-(\varepsilon D^{S_{\widetilde{M}}})f^+(\varepsilon D^{S_{\widetilde{M}}})&-g^-(\varepsilon D^{S_{\widetilde{M}}})^2
			\end{pmatrix}.
		\end{split}
	\end{equation}
	In operator K-theory, this difference idempotent $I_\varepsilon$ can be regarded as a small scale representative of the higher index of the Dirac operator $D^{S_{\widetilde{M}}}$. From \eqref{3.4}, we see that
	$I_\varepsilon(x,y)$  has propagation 
	$\tfrac{\varepsilon}{m},$ namely
	\begin{equation}\label{3.8}
		\supp (I_\varepsilon(x,\cdot))\subseteq B^{\widetilde{M}}(x,\tfrac{\varepsilon}{m})
	\end{equation}
	for all  $x\in \widetilde{M}$,
	where $I_\varepsilon(x,y)$ is the kernel of $I_\varepsilon$.

	\subsection{Small scale index theorem}\label{SS3}

	For any $k\in\mathbb{N}$, let us form the small scale index kernel $\mathrm{Ind}_\varepsilon(D^{S_{\widetilde{M}}})(x_{0},x_1,\cdots,x_{2k})$ of the Dirac operator $D^{S_{\widetilde{M}}}$ by
	\begin{equation}\label{3.10}
		\mathrm{Ind}_\varepsilon(D^{S_{\widetilde{M}}})\big(x_{0},x_1,\cdots,x_{2k}\big)=I_\varepsilon\big(x_{0},x_{1}\big)\cdots I_\varepsilon\big(x_{2k-1},x_{2k}\big)I_\varepsilon\big(x_{2k},x_{0}\big).
	\end{equation}
	Since $I_\varepsilon\in \mathscr{K}^{\widetilde{\Gamma}}$, we have
	\begin{equation}
		\mathrm{Ind}_\varepsilon(D^{S_{\widetilde{M}}})\big(x_{0},x_1,\cdots,x_{2k}\big)=\widetilde{\gamma}\mathrm{Ind}_\varepsilon(D^{S_{\widetilde{M}}})\big(\widetilde{\gamma}^{-1} x_{0},\cdots,\widetilde{\gamma}^{-1}x_{2k}\big)\widetilde{\gamma}^{-1},
	\end{equation}
	for any $\widetilde{\gamma}\in\widetilde{\Gamma}$, hence the pointwise trace is $\Gamma$-invariant, that is,
	\begin{equation}
		\mathrm{Tr}^{S_{\widetilde{M}}}\big[\mathrm{Ind}_\varepsilon(D^{S_{\widetilde{M}}})\big(x_{0},\cdots,x_{2k}\big)\big]=\mathrm{Tr}^{S_{\widetilde{M}}}\big[\mathrm{Ind}_\varepsilon(D^{S_{\widetilde{M}}})\big(\gamma x_{0},\cdots,\gamma x_{2k}\big)\big]
	\end{equation}
	for any $\gamma\in\Gamma$. Moreover, by \eqref{3.8}, we have
	\begin{equation}
		\mathrm{supp}\big(\mathrm{Ind}_\varepsilon(D^{S_{\widetilde{M}}})\big)\subseteq \mathrm{diag}_\varepsilon(\widetilde{M}^{2k+1}).
	\end{equation}
	

	\begin{defi}
		We define the Connes-Chern class  $\mathrm{Ch}_{2k}^\varepsilon\big(\mathrm{Ind}_\varepsilon(D^{S_{\widetilde{M}}})\big)\in C_{2k}^\varepsilon(\widetilde{M},\Gamma)$ of $D^{S_{\widetilde{M}}}$ as an antisymmetric measure
		\begin{equation}\label{3.13}
			\begin{split}
				&\mathrm{Ch}_{2k}\big(\mathrm{Ind}_\varepsilon(D^{S_{\widetilde{M}}})\big)\\
				&=\sum_{\sigma\in \mathfrak{G}_{2k+1}}\mathrm{sgn}(\sigma)\mathrm{Tr}^{S_{\widetilde{M}}}\big[\mathrm{Ind}_\varepsilon(D^{S_{\widetilde{M}}})\big(x_{\sigma(0)},\cdots,x_{\sigma(2k)}\big)\big]dv_{\widetilde{M}^{2k+1}}\big(x_{0},\cdots,x_{2k}\big).
			\end{split}
		\end{equation}		
	\end{defi}

	Now we state the small scale index theorem.
	\begin{theo}\label{T3.2}
		Suppose that $\widetilde{M}$ is spin, then for any $\varepsilon>0$, we have $\mathrm{Ch}_{2k}^\varepsilon\big(\mathrm{Ind}_\varepsilon(D^{S_{\widetilde{M}}})\big)\in H_{2k}^\varepsilon(\widetilde{M},\Gamma)$, and
		\begin{equation}\label{3.15..}
			\mathrm{Ch}_{2k}^\varepsilon\big(\mathrm{Ind}_\varepsilon(D^{S_{\widetilde{M}}})\big)=\chi_\varepsilon\Big(\widehat{A}^{[m-2k]}\big(T{M},\nabla^{T{M}}\big)\cap[M]\Big),
		\end{equation}
		where the chain map $\chi_\varepsilon$ is defined in \eqref{2.17}. Equivalently, for any $f\in H^{2k}_{\varepsilon}(\widetilde{M},\Gamma)$,
		\begin{equation}\label{3.14}
			\begin{split}
				\big\langle \mathrm{Ch}_{2k}\big(\mathrm{Ind}_\varepsilon (D^{S_{\widetilde{M}}})\big),f\big\rangle_\Gamma&=\big\langle\chi_\varepsilon(f),\widehat{A}^{[m-2k]}\big(T{M},\nabla^{T{M}}\big)\cap[M]\big\rangle\\
				&=\int_{M}\chi_\varepsilon(f)\wedge\widehat{A}\big(T{M},\nabla^{T{M}}\big),
			\end{split}
		\end{equation}
		where the cochain map $\chi_\varepsilon$ is given in \eqref{2.19}.
	\end{theo}
	
	\begin{remark}
		Theorem \ref{T3.2} expresses the image of $\widehat{A}^{[m-2k]}\big(T{M},\nabla^{T{M}}\big)\cap[M]$ in $\varepsilon$-homology. It is worth noting that this result holds without assuming $0<\varepsilon< r(M)$, where $r(M)$ is the homological radius given in Definition \ref{def2.3}.
	\end{remark}

	For any $f_0,\ldots,f_{2k}\in C_0^\infty(\widetilde{M})$, the space of smooth functions on $\widetilde{M}$ with compact support, we have
	\begin{equation}
		\chi_\varepsilon(f_0\otimes \cdots\otimes f_{2k})=f_0df_1\wedge\cdots\wedge df_{2k},
	\end{equation}
	where $f_0\otimes \cdots\otimes f_{2k}\in C^\infty_0(\widetilde{M}^{2k+1})$ is defined by $(f_0\otimes \cdots\otimes f_{2k})(x_0,  \cdots, x_{2k}) =f_0(x_0)\cdots f(x_{2k})$. Although the support of $f$ in \eqref{3.14} is never a compact subset of $\widetilde{M}^{2k+1}$, by \eqref{2.10''}, the integral is over the fundamental domain of the $\varepsilon$-neighbourhood of the diagonal, which is compact. Using a partition of unity argument, to deduce \eqref{3.14}, it is therefore sufficient to establish the following stronger result.
	\begin{prop}\label{p3.4}
		For any $f_0,\ldots,f_{2k}\in C_0^\infty(\widetilde{M})$, we have
		\begin{equation}\label{3.16}
			\lim_{\varepsilon\to 0}\mathrm{Ch}_{2k}\big(\mathrm{Ind}_\varepsilon (D^{S_{\widetilde{M}}})\big)(f_0\otimes \cdots\otimes f_{2k})=\int_{\widetilde{M}}f_0df_1\wedge\cdots\wedge df_{2k}\wedge\widehat{A}\big(T\widetilde{M},\nabla^{T\widetilde{M}}\big).
		\end{equation}
	\end{prop}

	The detailed proofs of Theorem \ref{T3.2} and Proposition \ref{p3.4} will be provided later, in \cref{SCCC}, \cref{SLIC}, and \cref{FAM}. Also, we remark that an odd dimensional analogue of the small scale index theorem can be formulated and proved using a standard suspension argument.

	\section{Simplicial norm estimates}\label{SNE}

	In this section, we introduce our main estimates on simplicial norms. In \cref{SNE1}, we review the definition of Gromov's simplicial norm \cite[\S\,0.2]{MR686042} and some related properties. In \cref{SNE2}, by combining our small scale index theorem and estimates on the kernel functions, we obtain our main result on simplicial norms.

	\subsection{Simplicial norms}\label{SNE1}
	
	First, let us define the simplicial norm of a homology class.	
	\begin{defi}
		For any singular homology class $\alpha\in H_*(M)$, its simplicial norm is defined by the following semi-norm
		\begin{equation}\label{3.17}
			\lV\alpha \rV_{\ell^1}=\inf_{[\sum_{i}a_i\Delta_i]=\alpha} \sum_i \lv a_i\rv,
		\end{equation}
		where the infimum is taken over all singular cycles $\sum_{i}a_i\Delta_i$ representing $\alpha$. Equivalently, $\lV\alpha\rV_{\ell^1}$ is the infimal number of simplices that we need to represent $\alpha$. In particular, the simplicial norm $\lV[M]\rV_{\ell^1}$ of the fundamental class $[M]\in H_m(M)$ is called the simplicial volume of $M$.
	\end{defi}

	For example, the simplicial volume $\lV[S^1]\rV_{\ell^1}$ of the circle $S^1$ vanishes. To check this, consider a simplex $\Delta_{100}$ that winds $100$ times around $S^1$, then $[S^1]=\frac{1}{100}\Delta_{100}$, and by \eqref{3.17} we have $\lV[S^1]\rV_{\ell^1}\leqslant \frac{1}{100}$. Replacing $100$ by an arbitrarily large number, we see that $\lV[S^1]\rV_{\ell^1}=0$. In contrast, Gromov-Thurston \cite[\S\,1.2]{MR686042} proved that if $M$ is a hyperbolic manifold, then there is a constant $C_m>0$  depending only on $m$ such that
	\begin{equation}
		\lV[M]\rV_{\ell^1}=C_m\mathrm{Vol}(M).
	\end{equation}
	This shows that hyperbolic manifolds carry richer topological complexity.

	Now for any $\mu\in H^\varepsilon_*(M)$, let us consider its simplicial norm $\lV\mu\rV_{\ell^1}$. If $0<\varepsilon<r(M)$, by \eqref{2.18} and \eqref{2.20}, we can formally view $\mu$ as an infinite linear combination of simplices
	\begin{equation}
		\mu=\sum_{(x_0,\cdots,x_k)\in \Gamma\backslash\widetilde{M}^{k+1}}\mu(x_0,\cdots,x_k)\Delta_{(x_0,\cdots,x_k)},
	\end{equation}
	where each $\Delta_{(x_0,\cdots,x_k)}$ having vertices $(x_0,\cdots,x_k)$. From \eqref{3.17}, it follows that
	\begin{equation}
		\lV\mu\rV_{\ell^1}\leqslant \sum_{(x_0,\cdots,x_k)\in \Gamma\backslash\widetilde{M}^{k+1}}\lv\mu(x_0,\cdots,x_k)\rv.
	\end{equation}
	To make sense of this infinite sum, it is natural to interpret it as an integral, leading to the following estimate.
	\begin{lemma}
		Suppose that $0<\varepsilon<r(M)$, then for any $\mu\in H^\varepsilon_*(M)\cong H_*(M)$, we have
		\begin{equation}\label{4.4..}
			\lV\mu\rV_{\ell^1}\leqslant \lv\mu \rv\big(\Gamma\backslash\widetilde{M}^{k+1}\big),
		\end{equation}
		where $\lv\mu\rv$ denotes the total variation measure of $\mu$.
	\end{lemma}
	This intuitive argument can be made rigorous using Gromov’s duality principle \cite[\S\,1.1]{MR686042}. Indeed, Thurston \cite[\S\,6]{MR4554426} suggested an alternative description of the simplicial norm in terms of measures on the total space of simplices, defining a seminorm as the infimum of total variations. Löh \cite[Theorem 1.1]{MR2206643} showed that these two constructions result in the same simplicial norm. For our purposes, the weaker upper bound \eqref{4.4..} is sufficient.

	\subsection{Kernel function estimates}\label{SNE2}

	Combining \eqref{3.13}, \eqref{3.15..} and \eqref{4.4..}, we get the following estimate on the simplicial norm of $\widehat{A}^{[m-2k]}\big(TM,\nabla^{TM}\big)\cap[M]$.
	\begin{theo}\label{T3.4}
		Suppose that $\widetilde{M}$ is spin. For any $0<\varepsilon<r(M)$, we have
		\begin{equation}\label{3.18}
			\begin{split}
				&\lV\widehat{A}^{[m-2k]}\big(TM,\nabla^{TM}\big)\cap[M]\rV_{\ell^1}\\
				&\leqslant\sum_{\sigma\in\mathfrak{G}_{2k+1}} \int_{\Gamma\backslash\widetilde{M}^{2k+1}} \lv\mathrm{Tr}^{S_{\widetilde{M}}}\big[\mathrm{Ind}^\varepsilon(D^{S_{\widetilde{M}}})\big(x_{\sigma(0)},\cdots,x_{\sigma(2k)}\big)\big]\rv dv_{\Gamma\backslash\widetilde{M}^{2k+1}}\big(x_{0},\cdots,x_{2k}\big).
			\end{split}
		\end{equation}
	\end{theo}

	To proceed, it remains to control the integrand in \eqref{3.18}. In \cref{SEE}, we establish the following bound for the index kernel, see Theorem \ref{T7.6} and Lemma \eqref{L8.7} for details.
	\begin{prop}
		We have for $0<\varepsilon<r(M)$
		\begin{equation}\label{4.7''}
			\lV I_\varepsilon(x_0,x_1)\rV_{S_{\widetilde{M},x_0}\otimes S_{\widetilde{M},x_1}^*}\leqslant C\Big(1+\frac{1}{\varepsilon^{m+1}}\Big)\Big(\frac{\sup_{x\in M} \mathrm{Sc}^-_{g^{TM}}(x)^{\tfrac{m}{2}}+1}{\mathrm{ID}\big(B^{\widetilde{M}}(x_0,\varepsilon)\big)^{m}}
			\Big),
		\end{equation}
		where $\mathrm{ID}\big(B^{\widetilde{M}}(x_0,\varepsilon)\big)$ is the Dirichlet isoperimetric constant of $B^{\widetilde{M}}(x_0,\varepsilon)$ later defined in \eqref{2.5''}.
	\end{prop}

	By Corollary \ref{Coro2.5}, \eqref{3.10}, \eqref{3.18} and \eqref{4.7''}, when $M$ is aspherical, we can take $\varepsilon=1$ to obtain the $C^0$-stable bound \eqref{1.4.'},
	\begin{equation}\label{4.8n}
		\begin{split}
			&\lV\widehat{A}^{[i]}(TM)\cap[M]\rV_{\ell^1}\\
			&\leqslant C_m\bigg(\frac{\sup_{x\in M} \mathrm{Sc}^-_{g^{TM}}(x)^{\tfrac{m}{2}}+1}{\inf_{x\in\widetilde{M}}\mathrm{ID}\big(B^{\widetilde{M}}(x,1)\big)^m}\bigg)^{m-i+1}\sup_{x\in\widetilde{M}}\mathrm{Vol}\big(B^{\widetilde{M}}(x,1)\big)^{m-i}\mathrm{Vol}(M).
		\end{split}
	\end{equation}

	In \cite[Lemma 6, Theorem 11]{MR608287}, Croke showed that if a region $\Omega\subset \widetilde{M}$ has diameter less than the injective radius $\mathrm{Inj}(\widetilde{M})$, then the isoperimetric constant $ID(\Omega)$ has a positive lower bound only depending on the dimension $m$.
	\begin{prop}
		There is a constant $C_m>0$ such that for any $i_0\leqslant \mathrm{Inj}(M)$
		\begin{equation}\label{4.8..}
			\mathrm{ID}\big(B^{\widetilde{M}}(x_0,i_0)\big)\geqslant C_m.
		\end{equation}
	\end{prop}
	
	From Corollary \ref{Coro2.5}, \eqref{3.10}, \eqref{3.18} \eqref{4.7''} and \eqref{4.8..}, we can take $\varepsilon=i_0$ to get
	\begin{equation}
		\begin{split}
			&\lV\widehat{A}^{[m-2k]}(TM)\cap[M]\rV_{\ell^1}\\
			&\leqslant C_m\Big(\big(1+i_0^{-m-1}\big)\big(\sup_{x\in M}\mathrm{Sc}^-_{g^{TM}}(x)^{\tfrac{m}{2}}+1\big)\Big)^{2k+1}\mathrm{Vol}(M)\sup_{\widetilde{x}\in \widetilde{M}}\mathrm{Vol}\big(B^{\widetilde{M}}(\widetilde{x},i_0)\big)^{2k}.
		\end{split}
	\end{equation}
	Finally, by first rescaling $g^{TM}$ to $i_0^{-2} g^{TM}$, and then to $t g^{TM}$ as $t \to \infty$, we obtain
	\begin{equation}\label{4.11n}
		\begin{split}
			&\lV\widehat{A}^{[i]}(TM)\cap[M]\rV_{\ell^1}\\
			&\leqslant C_m\bigg(\sup_{x\in M} \mathrm{Sc}^-_{g^{TM}}(x)^{\tfrac{m}{2}}+i_0^{-m}\bigg)^{m-i+1}\cdot\mathrm{Vol}(M)\cdot\sup_{\widetilde{x}\in \widetilde{M}}\mathrm{Vol}\big(B^{\widetilde{M}}\big(\widetilde{x},i_0\big)\big)^{m-i},
		\end{split}
	\end{equation}
	which is equivalent to \eqref{1.11111}.

	\begin{remark}
		Taking $\varepsilon=\mathrm{Inj}(M)/2$, we obtain
		\begin{equation}\label{4.10}
			\lV\widehat{A}^{[i]}(TM)\cap[M]\rV_{\ell^1}\leqslant C_m\bigg(\sup_{x\in M} \mathrm{Sc}^-_{g^{TM}}(x)^{\tfrac{m}{2}}+\mathrm{Inj}(M)^{-m}\bigg)^{m-i+1}\mathrm{Vol}(M)^{m-i+1}.
		\end{equation}
	Grove-Petersen-Wu \cite[Theorem C]{MR1029396} proved a curvature free generalization of Cheeger's finiteness theorem, that the class of closed manifolds with injectivity radius lower bound and volume upper bound contains at most finitely many homeomorphism types when $m\neq 3$, and only finitely many diffeomorphism types if, in addition, $m\neq 4$. Therefore, the estimate \eqref{4.10} is of less interest compared with the bounds in \eqref{4.8n} and \eqref{4.11n}.
	\end{remark}

	\section{Connes-Chern character}\label{SCCC}

	In this section, we discuss properties of the Connes-Chern character $\mathrm{Ch}_{2k}\big(\mathrm{Ind}_\varepsilon(D^{S_{\widetilde{M}}})\big)\in C_{2k}^\varepsilon(\widetilde{M},\Gamma)$ defined in \eqref{3.13} following Connes-Moscovici \cite[\S\,2]{MR1066176} and Moscovici-Wu \cite[\S\,1]{MR1292020}. In \cref{CCC1}, we introduce general trace homomorphism maps. In \cref{CCC2}, we show that $\mathrm{Ch}_{2k}\big(\mathrm{Ind}_\varepsilon(D^{S_{\widetilde{M}}})\big)\in C_{2k}^\varepsilon(\widetilde{M},\Gamma)$ is closed, thus $\mathrm{Ch}_{2k}\big(\mathrm{Ind}_\varepsilon(D^{S_{\widetilde{M}}})\big)\in H_{2k}^\varepsilon(\widetilde{M},\Gamma)$, and given  $\varepsilon_0>0$ we prove that this homology class
	$\mathrm{Ch}_{2k}\big(\mathrm{Ind}_\varepsilon(D^{S_{\widetilde{M}}})\big)$ is independent of  $\varepsilon\leqslant \varepsilon_0$ in $H_{2k}^{\varepsilon_0}(\widetilde{M},\Gamma)$. In \cref{CCC3}, we give a different form of $\mathrm{Ch}_{2k}\big(\mathrm{Ind}_\varepsilon(D^{S_{\widetilde{M}}})\big)$, which is useful in later local index computations.

	\subsection{Trace homomorphism}\label{CCC1}
	
	For $(K_0,\cdots,K_k)\in (\mathscr{K}^{\widetilde{\Gamma}})^{k+1}$, we associate a measure $\mu_k(K)\in C_{k}(\widetilde{M},\Gamma)$ such that for $f_0,\cdots,f_k\in C_0^\infty(\widetilde{M})$,
	\begin{equation}
		\begin{split}
			&\mu_k(K_0,\cdots,K_k)(f_0\otimes\cdots\otimes f_{k})\\
			&=\mathrm{Tr}^{L^2(\widetilde{M},S_{\widetilde{M}})}\big[f_0K_0\cdots f_kK_k\big]\\
			&=\int_{\widetilde{M}^{k+1}}\mathrm{Tr}^{S_{\widetilde{M}}}\big[f_0(x_0)K_0(x_0,x_1)\cdots f_{k}(x_k)K_{k}(x_k,x_0)\big]dv_{\widetilde{M}}(x_0)\cdots dv_{\widetilde{M}}(x_k),
		\end{split}
	\end{equation}
	in other words,
	\begin{equation}\label{4.2}
		\mu_k(K_0,\cdots,K_k)=\mathrm{Tr}^{S_{\widetilde{M}}}\big[K_0(x_0,x_1)K_1(x_1,x_2)\cdots K_{k}(x_k,x_0)\big]dv_{\widetilde{M}^{k+1}}\big(x_{0},\cdots,x_{k}\big).
	\end{equation}
	By the property of trace, we have
	\begin{equation}\label{4.3}
		\mu_k(K_0,\cdots,K_k)(f_0\otimes\cdots\otimes f_{k})=\mu_k(K_1,\cdots,K_k,K_0)(f_1\otimes\cdots\otimes f_{k}\otimes f_0).
	\end{equation}

	We define the antisymmetrization measure $\mathrm{Alt}(\mu_k)(K_0,\cdots,K_k)$ of $\mu_k(K_0,\cdots,K_k)$ by
	\begin{equation}\label{4.4}
		\begin{split}
			&\mathrm{Alt}(\mu_k)(K_0,\cdots,K_k)(f_0\otimes\cdots\otimes f_{k})\\
			&=\frac{1}{(k+1)!}\sum_{\sigma\in\mathfrak{G}_{k+1}}\mathrm{sgn}(\sigma)\mu_k(K_0,\cdots,K_k)\big(f_{\sigma(0)}\otimes\cdots\otimes f_{\sigma(k)}\big).
		\end{split}
	\end{equation}

	In general, for $(K_0+Z_0,\cdots,K_k+Z_k)$, where $K_i\in \mathscr{K}^{\widetilde{\Gamma}}$ and $Z_i\in \mathbb{C}\left(\begin{smallmatrix}
		1&0\\
		0&1
	\end{smallmatrix}\right)+\mathbb{C}\left(\begin{smallmatrix}
		1&0\\
		0&-1
	\end{smallmatrix}\right)$, we define
	\begin{equation}\label{4.5}
		\begin{split}
			&\mu_k(K_0+Z_0,\cdots,K_k+Z_k)(f_0\otimes\cdots\otimes f_{k})\\
			&=\mathrm{Tr}^{L^2(\widetilde{M},S_{\widetilde{M}})}\big[f_0(K_0+Z_0)\cdots f_k(K_k+Z_k)-f_0Z_0\cdots f_kZ_k\big],
		\end{split}
	\end{equation}
	and its antisymmetrization $\mathrm{Alt}(\mu_k)(K_0+Z_0,\cdots,K_k+Z_k)$ is defined similar to \eqref{4.4} by
	\begin{equation}
		\begin{split}
			&\mathrm{Alt}(\mu_k)(K_0+Z_0,\cdots,K_k+Z_k)(f_0\otimes\cdots\otimes f_{k})\\
			&=\frac{1}{(k+1)!}\sum_{\sigma\in\mathfrak{G}_{k+1}}\mathrm{sgn}(\sigma)\mu_k(K_0+Z_0,\cdots,K_k+Z_k)\big(f_{\sigma(0)}\otimes\cdots\otimes f_{\sigma(k)}\big)
		\end{split}
	\end{equation}

	
	\begin{prop}
		We have the following identity
		\begin{equation}\label{4.7}
			\mathrm{Alt}(\mu_k)(K+Z)=\mathrm{Alt}(\mu_k)(K).
		\end{equation}	
		Moreover, the asymmetrization measure is cyclic, 
		\begin{equation}\label{6.14'}
			\mathrm{Alt}(\mu_k)(K_1,\cdots,K_k,K_0)=(-1)^k\mathrm{Alt}(\mu_k)(K_0,K_1,\cdots,K_k).
		\end{equation}	
		Also, the boundary operator can be expressed by
		\begin{equation}\label{6.15}
			\begin{split}
				\pa\mathrm{Alt}(\mu_k)(K_0,K_1,\cdots,K_k)=\sum_{j=0}^{k}&(-1)^j \mathrm{Alt}(\mu_{k-1})(K_0,\cdots,K_jK_{j+1},\cdots,K_{k+1})\\
				&+(-1)^{k+1}\mathrm{Alt}(\mu_{k-1})(K_{k+1}K_0,\cdots,K_k).
			\end{split}
		\end{equation}
	\end{prop}
	
	\begin{pro}
		For any term inside \eqref{4.5}, there exist a product $f_iZ_if_{i+1}$, and we have $f_iZ_if_{i+1}=f_{i+1}Z_if_{i}$, and \eqref{4.7} follows from antisymmetrization.	
		
		We compute that
		\begin{equation}
			\begin{split}
				\mu_k(K_1,\cdots,K_k,K_0)(f_0\otimes\cdots\otimes f_{k})&=\mathrm{Tr}^{L^2(\widetilde{M},S_{\widetilde{M}})}\big[f_0K_1f_1K_2\cdots f_kK_0\big]\\
				&=\mathrm{Tr}^{L^2(\widetilde{M},S_{\widetilde{M}})}\big[f_kK_0f_0K_1f_1K_2\cdots f_{k-1}K_k\big]\\
				&=\mu_k(K_0,\cdots,K_k)(f_k\otimes f_0\cdots\otimes f_{k-1}),
			\end{split}
		\end{equation}
		and we also have $\mathrm{sgn}\binom{0,\cdots,k}{k,0,\cdots,k-1}=(-1)^k$, then we get \eqref{6.14'}.
		
		Recall that the coboundary operator satisfies
		\begin{equation}\label{4.11}
			\begin{split}
				\pa(f_0\otimes\cdots\otimes f_k)(x_0,\cdots,x_{k+1})&=\sum_{j=0}^{k+1}(-1)^j(f_0\otimes\cdots\otimes f_k)(x_0,\cdots,\widehat{x}_j,\cdots,x_{k+1})\\
				&=\sum_{j=0}^{k+1}(-1)^jf_0(x_0)\cdots f_{j-1}(x_{j-1})f_j(x_{j+1})\cdots f_k(x_{k+1}).
			\end{split}
		\end{equation}
		Combining \eqref{4.2} and the fact that
		\begin{equation}
			\int_{\widetilde{M}}K_{j-1}(x_{j-1},x_j)K_j(x_j,x_{j+1})dv_{\widetilde{M}}(x_j)=(K_{j-1}K_j)(x_{j-1},x_{j+1}),
		\end{equation}
		we get
		\begin{equation}\label{4.13}
			\begin{split}
				&\mu_{k+1}(K_0,\cdots,K_{k+1})(f_0\otimes\cdots f_{j-1}\otimes 1\otimes f_j\otimes\cdots f_{k})\\
				&=\mu_{k}(K_0,\cdots,K_{j-1}K_{j},\cdots,K_{k+1})(f_0\otimes\cdots\otimes f_k).
			\end{split}
		\end{equation}
		From \eqref{4.11} and \eqref{4.13} we obtain
		\begin{equation}
			\begin{split}
				\pa\mu_k(K_0,K_1,\cdots,K_k)=\sum_{j=0}^{k}&(-1)^j \mu_{k}(K_0,\cdots,K_jK_{j+1},\cdots,K_{k+1})\\
				&+(-1)^{k+1}\mu_k(K_{k+1}K_0,\cdots,K_k),
			\end{split}
		\end{equation}
		and antisymmetrization implies \eqref{6.15}.\qed 
	\end{pro}

	\subsection{Connes-Chern character}\label{CCC2}
	
	From \eqref{3.13}, \eqref{4.4} and \eqref{4.7}, Connes-Chern character can be expressed by
	\begin{equation}\label{4.15}
		\begin{split}
			\mathrm{Ch}_{2k}\big(\mathrm{Ind}_\varepsilon(D^{S_{\widetilde{M}}})\big)&=(2\pi\sqrt{-1})^k\frac{(2k)!}{k!}\mathrm{Alt}(\mu_{2k})\big(I_\varepsilon,\cdots,I_\varepsilon\big)\\
			&=(2\pi\sqrt{-1})^k\frac{(2k)!}{k!}\mathrm{Alt}(\mu_{2k})\big(P_\varepsilon,\cdots,P_\varepsilon\big)
		\end{split}
	\end{equation}
	We note that
	\begin{equation}\label{6.25.}
		\mathrm{Alt}(\mu_{2k-1})\big(I_\varepsilon,\cdots,I_\varepsilon\big)=0
	\end{equation}
	by \eqref{4.3} and antisymmetization.

	\begin{theo}
		We have
		\begin{equation}\label{6.26.}
			\mathrm{Ch}_{2k}\big(\mathrm{Ind}_\varepsilon(D^{S_{\widetilde{M}}})\big)\in H^\varepsilon_{2k}\big(\widetilde{M},\Gamma\big).
		\end{equation}
		Moreover, we have a transgression formula
		\begin{equation}\label{6.27.}
			\frac{\pa}{\pa\varepsilon}\mathrm{Ch}_{2k}\big(\mathrm{Ind}_\varepsilon(D^{S_{\widetilde{M}}})\big)=(2\pi\sqrt{-1})^k\frac{(2k+1)}{k!}\pa\mathrm{Alt}(\mu_{2k-1})\big(T_\varepsilon,P_\varepsilon,\cdots,P_\varepsilon\big)
		\end{equation}
		where $T_\varepsilon$ is defined by
		\begin{equation}
			T_\varepsilon=\big(1-2P_\varepsilon\big)\frac{\pa P_\varepsilon}{\pa\varepsilon}.
		\end{equation}
		Also, given $\varepsilon_0>0$ the homology classes
		\begin{equation}
			\{\mathrm{Ch}_{2k}\big(\mathrm{Ind}_{\varepsilon}(D^{S_{\widetilde{M}}})\big)\}_{\varepsilon\leqslant \varepsilon_0 }\in H^{\varepsilon_0}_{2k}\big(\widetilde{M},\Gamma\big)
		\end{equation}is independent on $\varepsilon\leqslant \varepsilon_0$.
	\end{theo}
	
	\begin{pro}
		To get \eqref{6.26.}, we compute by \eqref{4.15}, \eqref{6.15} that
		\begin{equation}
			\begin{split}
				\pa \mathrm{Alt}(\mu_{2k})\big(I_\varepsilon,\cdots,I_\varepsilon\big)&=\pa \mathrm{Alt}(\mu_{2k})\big(P_\varepsilon,\cdots,P_\varepsilon\big)\\
				&=\sum_{j}(-1)^j\mathrm{Alt}(\mu_{2k-1})\big(P_\varepsilon,\cdots,P_\varepsilon^2,\cdots,P_\varepsilon\big),
			\end{split}
		\end{equation}
		and this vanishes by \eqref{6.25.}.
		
		Since
		\begin{equation}
			\begin{split}
				\frac{\pa P_\varepsilon}{\pa\varepsilon} P_\varepsilon+P_\varepsilon\frac{\pa P_\varepsilon}{\pa\varepsilon}=\frac{\pa P_\varepsilon}{\pa\varepsilon},
			\end{split}
		\end{equation}
		we get
		\begin{equation}
			P_\varepsilon\frac{\pa P_\varepsilon}{\pa\varepsilon} P_\varepsilon=0.
		\end{equation}
		It follows that
		\begin{equation}
			\begin{split}
				\big[T_\varepsilon,P_\varepsilon\big]&=\big(1-2P_\varepsilon\big)\frac{\pa P_\varepsilon}{\pa\varepsilon}P_\varepsilon-P_\varepsilon\big(1-2P_\varepsilon\big)\frac{\pa P_\varepsilon}{\pa\varepsilon}\\
				&=\big(1-2P_\varepsilon\big)\frac{\pa P_\varepsilon}{\pa\varepsilon} P_\varepsilon+P_\varepsilon\frac{\pa P_\varepsilon}{\pa\varepsilon}\\
				&=\frac{\pa P_\varepsilon}{\pa\varepsilon}
			\end{split},
		\end{equation}
		then we have
		\begin{equation}\label{6.33.}
			\begin{split}
				\frac{\pa}{\pa\varepsilon}\mathrm{Alt}(\mu_{2k})\big(P_\varepsilon,\cdots,P_\varepsilon\big)&=\sum_{}^{}\mathrm{Alt}(\mu_{2k})\big(P_\varepsilon,\cdots,\big[T_\varepsilon,P_\varepsilon\big],\cdots,P_\varepsilon\big)\\
				&=(2k+1)\mathrm{Alt}(\mu_{2k})\big(\big[T_\varepsilon,P_\varepsilon\big],P_\varepsilon,\cdots,P_\varepsilon\big).
			\end{split}
		\end{equation}
		From \eqref{6.15}, we have
		\begin{equation}\label{6.34.}
			\begin{split}
				\pa&\mathrm{Alt}(\mu_{2k+1})\big(T_\varepsilon,P_\varepsilon,\cdots,P_\varepsilon\big)\\
				=&\mathrm{Alt}(\mu_{2k})\big(T_\varepsilon P_\varepsilon,P_\varepsilon,\cdots,P_\varepsilon\big)-\mathrm{Alt}(\mu_{2k})\big(P_\varepsilon T_\varepsilon,P_\varepsilon,\cdots,P_\varepsilon\big)\\
				&+\sum_{j=1}^{2k}(-1)^j\mathrm{Alt}(\mu_{2k})\big(T_\varepsilon,P_\varepsilon,\cdots,P_\varepsilon\big)\\
				=&\mathrm{Alt}(\mu_{2k})\big([T_\varepsilon,P_\varepsilon],P_\varepsilon,\cdots,P_\varepsilon\big).
			\end{split}
		\end{equation}
		Combining \eqref{6.33.} and \eqref{6.34.}, we obtain \eqref{6.27.}.\qed
	\end{pro}

	\subsection{A equivalent form of Connes-Chern character}\label{CCC3}

	The following equivalent expression for the Connes-Chern character is crucial for the subsequent local index computation.
	\begin{prop}
		We have
		\begin{equation}\label{4.27}
			\begin{split}
				&\mathrm{Ch}_{2k}\big(\mathrm{Ind}_\varepsilon D^{S_{\widetilde{M}}}\big)\big(f_0\otimes\cdots\otimes f_{2k}\big)\\
				&=\frac{(2\pi\sqrt{-1})^k}{k!}\sum_{\sigma\in\mathfrak{G}_{2k}}\mathrm{Tr}^{L^2(\widetilde{M},S_{\widetilde{M}})}\Big(I_\varepsilon f_0\big[I_\varepsilon,f_{\sigma(1)}\big]\cdots\big[I_\varepsilon,f_{\sigma(2k)}\big]\Big).
			\end{split}
		\end{equation}
	\end{prop}

	\begin{pro}
		We need three auxiliary identities
		\begin{equation}\label{4.28}
			\begin{split}
				&\mathrm{Tr}^{L^2(\widetilde{M},S_{\widetilde{M}})}\Big(\sum_{\sigma\in\mathfrak{G}_{2k+1}}\mathrm{sgn}(\sigma)f_{\sigma(0)}I_\varepsilon f_{\sigma(1)}I_\varepsilon\cdots f_{\sigma(2k)}I_\varepsilon\Big)\\
				&=\mathrm{Tr}^{L^2(\widetilde{M},S_{\widetilde{M}})}\Big(\sum_{\sigma\in\mathfrak{G}_{2k}}\mathrm{sgn}(\sigma)f_0I_\varepsilon f_{\sigma(1)}I_\varepsilon\cdots f_{\sigma(2k)}I_\varepsilon\Big),\\
				&\mathrm{Tr}^{L^2(\widetilde{M},S_{\widetilde{M}})}\Big(\sum_{\sigma\in\mathfrak{G}_{2k}}\mathrm{sgn}(\sigma)f_0I_\varepsilon f_{\sigma(1)}I_\varepsilon\cdots f_{\sigma(2k)}I_\varepsilon\Big)\\
				&=\mathrm{Tr}^{L^2(\widetilde{M},S_{\widetilde{M}})}\Big(\sum_{\sigma\in\mathfrak{G}_{2k}}\mathrm{sgn}(\sigma)f_0P_\varepsilon f_{\sigma(1)}P_\varepsilon\cdots f_{\sigma(2k)}P_\varepsilon\Big),\\
				&\mathrm{Tr}^{L^2(\widetilde{M},S_{\widetilde{M}})}\Big(\sum_{\sigma\in\mathfrak{G}_{2k}}\mathrm{sgn}(\sigma)I_\varepsilon f_0\big[I_\varepsilon,f_{\sigma(1)}\big]\cdots\big[I_\varepsilon,f_{\sigma(2k)}\big]\Big)\\
				&=\mathrm{Tr}^{L^2(\widetilde{M},S_{\widetilde{M}})}\Big(\sum_{\sigma\in\mathfrak{G}_{2k}}\mathrm{sgn}(\sigma)P_\varepsilon f_0\big[P_\varepsilon,f_{\sigma(1)}\big]\cdots\big[P_\varepsilon,f_{\sigma(2k)}\big]\Big).
			\end{split}
		\end{equation}
		
		The first identity follows immediately from \eqref{6.14} and \eqref{4.15}.

		For the second identity, we expand the product $f_0I_\varepsilon f_{\sigma(1)}I_\varepsilon\cdots f_{\sigma(2k)}I_\varepsilon$. Notice that if the following term appears,
		\begin{equation}
			f_{\sigma(i)}\left(\begin{smallmatrix} 0&0\\
				0&1
			\end{smallmatrix}\right)f_{\sigma(i+1)},
		\end{equation}
		the sum vanishes because of antisymmetrization. Therefore, we only need to consider the following two terms
		\begin{equation}\label{4.30}
			\begin{split}
				&\mathrm{Tr}^{L^2(\widetilde{M},S_{\widetilde{M}})}\Big(f_0\left(\begin{smallmatrix} 0&0\\
					0&1
				\end{smallmatrix}\right)f_{\sigma(1)}I_\varepsilon\cdots f_{\sigma(2k)}I_\varepsilon\Big),\\ &\mathrm{Tr}^{L^2(\widetilde{M},S_{\widetilde{M}})}\Big(f_0I_\varepsilon f_{\sigma(1)}I_\varepsilon\cdots f_{\sigma(2k)}\left(\begin{smallmatrix} 0&0\\
					0&1
				\end{smallmatrix}\right)\Big).
			\end{split}
		\end{equation}
		By the property of trace, we have
		\begin{equation}
			\begin{split}
				&\mathrm{Tr}^{L^2(\widetilde{M},S_{\widetilde{M}})}\Big(f_0\left(\begin{smallmatrix} 0&0\\
					0&1
				\end{smallmatrix}\right)f_{\sigma(1)}I_\varepsilon\cdots f_{\sigma(2k)}I_\varepsilon\Big)\\
				&=\mathrm{Tr}^{L^2(\widetilde{M},S_{\widetilde{M}})}\Big(f_0I_\varepsilon f_{\sigma(2)}I_\varepsilon\cdots f_{\sigma(2k)}I_\varepsilon f_{\sigma(1)}\left(\begin{smallmatrix} 0&0\\
					0&1
				\end{smallmatrix}\right)\Big)
			\end{split},
		\end{equation}
		therefore, the antisymmetric sum of the two terms in \eqref{4.30} cancels out.
		
		For the third identity, since $[I_\varepsilon,f_i]=[P_\varepsilon,f_0]$, we only need to show that
		\begin{equation}
			\sum_{\sigma\in\mathfrak{G}_{2k}}\mathrm{sgn}(\sigma)\mathrm{Tr}^{L^2(\widetilde{M},S_{\widetilde{M}})}\Big(\left(\begin{smallmatrix} 0&0\\
				0&1
			\end{smallmatrix}\right)f_0\big[I_\varepsilon,f_{\sigma(1)}\big]\cdots\big[I_\varepsilon,f_{\sigma(2k)}\big]\Big)=0.
		\end{equation}
		Let us expand $\left(\begin{smallmatrix} 0&0\\
			0&1
		\end{smallmatrix}\right)f_0[I_\varepsilon,f_{\sigma(1)}]\cdots[I_\varepsilon,f_{\sigma(2k)}]$, and observe that if the following term appears in the product
		\begin{equation}
			I_\varepsilon f_{\sigma(i)}f_{\sigma(i+1)}I_\varepsilon,
		\end{equation}
		due to antisymmetrization, the sum vanishes. Therefore, we only need to consider the following terms
		\begin{equation}\label{4.34}
			\begin{split}
				&\mathrm{Tr}^{L^2(\widetilde{M},S_{\widetilde{M}})}\Big(\left(\begin{smallmatrix} 0&0\\
					0&1
				\end{smallmatrix}\right)f_0f_{\sigma(1)}I_\varepsilon\cdots f_{\sigma(i)}I_\varepsilon I_\varepsilon f_{\sigma(i)}\cdots I_\varepsilon f_{\sigma(2k)}\Big),\\
				&\mathrm{Tr}^{L^2(\widetilde{M},S_{\widetilde{M}})}\Big(\left(\begin{smallmatrix} 0&0\\
					0&1
				\end{smallmatrix}\right)f_0I_\varepsilon f_{\sigma(1)}\cdots I_\varepsilon f_{\sigma(2k)}\Big),\\ &\mathrm{Tr}^{L^2(\widetilde{M},S_{\widetilde{M}})}\Big(\left(\begin{smallmatrix} 0&0\\
					0&1
				\end{smallmatrix}\right)f_0f_{\sigma(1)}I_\varepsilon\cdots f_{\sigma(2k)}I_\varepsilon\Big).
			\end{split}
		\end{equation}
		By the property of trace, the first term is equal to
		\begin{equation}
			\begin{split}
				&\mathrm{Tr}^{L^2(\widetilde{M},S_{\widetilde{M}})}\Big(I_\varepsilon f_{\sigma(2k)}\left(\begin{smallmatrix} 0&0\\
					0&1
				\end{smallmatrix}\right)f_0f_{\sigma(1)}I_\varepsilon\cdots f_{\sigma(i)}I_\varepsilon I_\varepsilon f_{\sigma(i)}\cdots I_\varepsilon f_{\sigma(2k-1)}\Big)\\
				&=\mathrm{Tr}^{L^2(\widetilde{M},S_{\widetilde{M}})}\Big(\left(\begin{smallmatrix} 0&0\\
					0&1
				\end{smallmatrix}\right)f_0f_{\sigma(2k)}I_\varepsilon f_{\sigma(2)}I_\varepsilon\cdots f_{\sigma(i)}I_\varepsilon I_\varepsilon f_{\sigma(i)}\cdots I_\varepsilon f_{\sigma(2k-1)}I_\varepsilon f_{\sigma(1)}\Big)\\
			\end{split},
		\end{equation}
		due to antisymmetrization, the sum vanishes. By the property of trace, we have
		\begin{equation}
			\begin{split}
				&\mathrm{Tr}^{L^2(\widetilde{M},S_{\widetilde{M}})}\Big(\left(\begin{smallmatrix} 0&0\\
					0&1
				\end{smallmatrix}\right)f_0I_\varepsilon f_{\sigma(1)}\cdots I_\varepsilon f_{\sigma(2k)}\Big)\\
				&=\mathrm{Tr}^{L^2(\widetilde{M},S_{\widetilde{M}})}\Big(\left(\begin{smallmatrix} 0&0\\
					0&1
				\end{smallmatrix}\right)f_0f_{\sigma(2k)}I_\varepsilon f_{\sigma(1)}I_\varepsilon\cdots f_{\sigma(2k-1)}I_\varepsilon\Big)
			\end{split},
		\end{equation}
		therefore, because of antisymmetrization, the sum of the last two terms in \eqref{4.34} cancel out.

		Now let us return to prove \eqref{4.27}. By \eqref{4.28}, $P_\varepsilon^2=P_\varepsilon$, and the property of trace, it suffices to show that
		\begin{equation}\label{4.37}
			\begin{split}
				&\mathrm{Tr}^{L^2(\widetilde{M},S_{\widetilde{M}})}\Big(\sum_{\sigma\in\mathfrak{G}_{2k}}\mathrm{sgn}(\sigma)P_\varepsilon f_0\big[P_\varepsilon,f_{\sigma(1)}\big]\cdots\big[P_\varepsilon,f_{\sigma(2k)}\big]\Big)\\
				&=\mathrm{Tr}^{L^2(\widetilde{M},S_{\widetilde{M}})}\Big(\sum_{\sigma\in\mathfrak{G}_{2k}}\mathrm{sgn}(\sigma)P_\varepsilon f_0P_\varepsilon f_{\sigma(1)}\cdots P_\varepsilon f_{\sigma(2k)}P_\varepsilon\Big).
			\end{split}
		\end{equation}
		We compute that
		\begin{equation}\label{4.38}
			\begin{split}
				P_\varepsilon f_{\sigma(i)}P_\varepsilon f_{\sigma(i+1)}P_\varepsilon=&\big[P_\varepsilon ,f_{\sigma(i)}\big]P_\varepsilon f_{\sigma(i+1)}P_\varepsilon+f_{\sigma(i)}P_\varepsilon f_{\sigma(i+1)}P_\varepsilon \\
				=&\big[P_\varepsilon ,f_{\sigma(i)}\big]\big[P_\varepsilon ,f_{\sigma(i+1)}\big]P_\varepsilon +\big[P_\varepsilon ,f_{\sigma(i)}\big]f_{\sigma(i+1)}P_\varepsilon +f_{\sigma(i)}P_\varepsilon f_{\sigma(i+1)}P_\varepsilon \\
				=&\big[P_\varepsilon ,f_{\sigma(i)}\big]\big[P_\varepsilon ,f_{\sigma(i+1)}\big]P_\varepsilon +P_\varepsilon f_{\sigma(i)}f_{\sigma(i+1)}P_\varepsilon,
			\end{split}
		\end{equation}
		and similarly,
		\begin{equation}\label{4.39}
			\begin{split}
				P_\varepsilon f_{\sigma(i)}P_\varepsilon f_{\sigma(i+1)}P_\varepsilon=&P_\varepsilon f_{\sigma(i)}P_\varepsilon \big[f_{\sigma(i+1)},P_\varepsilon \big]+P_\varepsilon f_{\sigma(i)}P_\varepsilon f_{\sigma(i+1)}\\
				=&P_\varepsilon \big[f_{\sigma(i)},P_\varepsilon \big]\big[f_{\sigma(i+1)},P_\varepsilon \big]+P_\varepsilon f_{\sigma(i)}\big[f_{\sigma(i+1)},P_\varepsilon \big]+P_\varepsilon f_{\sigma(i)}P_\varepsilon f_{\sigma(i+1)}\\
				=&\big[P_\varepsilon ,f_{\sigma(i)}\big]\big[P_\varepsilon ,f_{\sigma(i+1)}\big]P_\varepsilon +P_\varepsilon f_{\sigma(i)}f_{\sigma(i+1)}P_\varepsilon .
			\end{split}
		\end{equation}
		By \eqref{4.38} and \eqref{4.39}, we get
		\begin{equation}\label{4.40}
			\begin{split}
				\big[P_\varepsilon ,f_{\sigma(i)}\big]\big[P_\varepsilon ,f_{\sigma(i+1)}\big]P_\varepsilon&=P_\varepsilon \big[P_\varepsilon ,f_{\sigma(i)}\big]\big[P_\varepsilon ,f_{\sigma(i+1)}\big]\\
				&=P_\varepsilon \big[P_\varepsilon ,f_{\sigma(i)}\big]\big[P_\varepsilon ,f_{\sigma(i+1)}\big]P_\varepsilon .
			\end{split}
		\end{equation}
		Combining \eqref{4.38} and \eqref{4.40}, due to antisymmetrization, we obtain \eqref{4.37}.\qed
	\end{pro}

	\section{Local index computation}\label{SLIC}

	In this section, we derive the local index formula \eqref{3.16}. In \cref{LIC1}, we introduce a symmetrized version of the Connes-Chern character $\mathrm{Ch}_{2k}\big(\mathrm{Ind}_\varepsilon D^{S_{\widetilde{M}}}\big)$ following Moscovici-Wu \cite[\S\,3.1]{MR1254310}. In \cref{LIC2}, we trivialize the local geometry to reduce the index computation on a manifold to one on the Euclidean space. In \cref{LIC3}, we apply Getzler's rescaling \cite{MR836727} to the localized operators. In \cref{LIC4}, we compute the limit of the rescaled operators. In \cref{LIC5}, we we evaluate the kernel of the resulting limit operator.


	\subsection{Symmetrized cycle}\label{LIC1}

	Let us recall the functions $f(t),g(t)$ used in \eqref{3.5'}, then let $\mathcal{L}_\varepsilon$ be an invertible operator acting on $L^2(\widetilde{M},S_{\widetilde{M}})^{\oplus 2}$ with its inverse $\mathcal{L}_\varepsilon^{-1}$ given by
	\begin{equation}
		\begin{split}
			\mathcal{L}_\varepsilon&=\begin{pmatrix} -g(\varepsilon D^{S_{\widetilde{M}}})&(g(\varepsilon D^{S_{\widetilde{M}}})-1)f(\varepsilon D^{S_{\widetilde{M}}})\\
				f(\varepsilon D^{S_{\widetilde{M}}})&-g(\varepsilon D^{S_{\widetilde{M}}})
			\end{pmatrix},\\
			\mathcal{L}_\varepsilon^{-1}&=\begin{pmatrix} -g(\varepsilon D^{S_{\widetilde{M}}})&-(g(\varepsilon D^{S_{\widetilde{M}}})-1)f(\varepsilon D^{S_M})\\
				-f(\varepsilon D^{S_{\widetilde{M}}})&-g(\varepsilon D^{S_{\widetilde{M}}})
			\end{pmatrix}.
		\end{split}
	\end{equation}
	Let $\mathcal{P}_\varepsilon$ be the corresponding idempotent defined by
	\begin{equation}
		\mathcal{P}_\varepsilon=\mathcal{L}_\varepsilon\begin{pmatrix} (\begin{smallmatrix}
				1&0\\
				0&0
			\end{smallmatrix})&0\\
			0&(\begin{smallmatrix}
				1&0\\
				0&0
			\end{smallmatrix})
		\end{pmatrix}\mathcal{L}_\varepsilon^{-1},
	\end{equation}
	and $\mathcal{I}_\varepsilon$ the ``difference idempotent" given by
	\begin{equation}\label{5.3}
		\begin{split}
			\mathcal{I}_\varepsilon&=\mathcal{P}_\varepsilon-\begin{pmatrix} (\begin{smallmatrix}
					0&0\\
					0&1
				\end{smallmatrix})&0\\
				0&(\begin{smallmatrix}
					0&0\\
					0&1
				\end{smallmatrix})
			\end{pmatrix}\\
			&=\begin{pmatrix} g(\varepsilon D^{S_{\widetilde{M}}})^2\tau&-g(\varepsilon D^{S_{\widetilde{M}}})\big(1-g(\varepsilon D^{S_{\widetilde{M}}})\big)f(\varepsilon D^{S_{\widetilde{M}}})\tau\\
				g(\varepsilon D^{S_{\widetilde{M}}})f(\varepsilon D^{S_{\widetilde{M}}})\tau&g(\varepsilon D^{S_{\widetilde{M}}})^2\tau
			\end{pmatrix}.
		\end{split}
	\end{equation}
	where $\tau$ the grading operator on $S_{\widetilde{M}}$ give by
	\begin{equation}\label{5.4}
		\tau=\begin{pmatrix}
			1&0\\
			0&-1
		\end{pmatrix}.
	\end{equation}
	
	For $0\leqslant s\leqslant1$, set
	\begin{equation}
		a(s)=\begin{pmatrix} 1&0&0&0\\
			0&\cos(\tfrac{\pi}{2}s)&0&-\sin(\tfrac{\pi}{2}s)\\
			0&0&1&0\\
			0&\sin(\tfrac{\pi}{2}s)&0&\cos(\tfrac{\pi}{2}s)
		\end{pmatrix},
	\end{equation}
	then we have a homotopy
	\begin{equation}
		\mathcal{P}_\varepsilon\xrightarrow{\mathrm{Ad}_{a(s)}}\begin{pmatrix} P_\varepsilon&0\\
			0&P_\varepsilon^*
		\end{pmatrix},
	\end{equation}
	where $P_\varepsilon^*$ is the adjoint of $P_\varepsilon$. Since $D^{S_{\widetilde{M}}}$ is self-adjoint, if we replace $P_\varepsilon$ with its adjoint $P_\varepsilon^*$ in \eqref{4.15}, the resulting homology class remains the same. Therefore, by \eqref{4.27}, we get the required symmetrized Connes-Chern character
	\begin{equation}
		\begin{split}
			&\mathrm{Ch}_{2k}\big(\mathrm{Ind}_\varepsilon D^{S_{\widetilde{M}}}\big)\big(f_0\otimes\cdots\otimes f_{2k}\big)\\
			&=\frac{(2\pi\sqrt{-1})^k}{2\cdot k!}\sum_{\sigma\in\mathfrak{G}_{2k}}\mathrm{Tr}^{L^2(\widetilde{M},S_{\widetilde{M}})^{\oplus 2}}\Big(\mathcal{I}_\varepsilon f_0\big[\mathcal{I}_\varepsilon,f_{\sigma(1)}\big]\cdots\big[\mathcal{I}_\varepsilon,f_{\sigma(2k)}\big]\Big),
		\end{split}
	\end{equation}
	also, we can express it in terms of super-traces by
	\begin{equation}\label{5.8}
		\begin{split}
			&\mathrm{Ch}_{2k}\big(\mathrm{Ind}_\varepsilon D^{S_{\widetilde{M}}}\big)\big(f_0\otimes\cdots\otimes f_{2k}\big)\\
			&=\frac{(2\pi\sqrt{-1})^k}{2\cdot k!}\sum_{\sigma\in\mathfrak{G}_{2k}}\int_{\widetilde{M}}\mathrm{Tr}^{S_{\widetilde{M}}^{\oplus2}}_s\Big(\tau\mathcal{I}_\varepsilon f_0\big[\mathcal{I}_\varepsilon,f_{\sigma(1)}\big]\cdots\big[\mathcal{I}_\varepsilon,f_{\sigma(2k)}\big](x,x)\Big)dv_{\widetilde{M}}(x),
		\end{split}
	\end{equation}
	where $\tau$ the grading operator given in \eqref{5.4}.
	
	Let us denote the integrand term in \eqref{5.8} by
	\begin{equation}\label{5.9}
		\mathcal{K}_\varepsilon=\tau\mathcal{I}_\varepsilon f_0\big[\mathcal{I}_\varepsilon,f_{1}\big]\cdots\big[\mathcal{I}_\varepsilon,f_{2k}\big]
	\end{equation}
	and analyze its structure. We define functions $v(\lambda),w(\lambda)$ such that
	\begin{equation}\label{5.10'}
		v(\lambda^2)=\frac{f(\lambda)}{\lambda},\ \ w(\lambda^2)=g(\lambda),
	\end{equation}
	where $f(t),g(t)$ are used in \eqref{3.5'}, and put
	\begin{equation}\label{5.10}
		\begin{split}
			\mathcal{A}_\varepsilon&=\begin{pmatrix} w^2(\varepsilon^2 D^{S_{\widetilde{M}},2})&0\\
				0&w^2(\varepsilon^2 D^{S_{\widetilde{M}},2})
			\end{pmatrix},\\
			\mathcal{B}_\varepsilon&=\begin{pmatrix} 0&-\big(w(1-w)v\big)(\varepsilon^2 D^{S_{\widetilde{M}},2})\varepsilon D^{S_{\widetilde{M}}}\\
				(wv)(\varepsilon^2 D^{S_{\widetilde{M}},2})\varepsilon D^{S_{\widetilde{M}}}&0
			\end{pmatrix},\\
			\mathcal{A}_{\varepsilon,i}&=\begin{pmatrix} [w^2(\varepsilon^2 D^{S_{\widetilde{M}},2}),f_i]&0\\
				0&[w^2(\varepsilon^2 D^{S_{\widetilde{M}},2}),f_i]
			\end{pmatrix},\\
			\mathcal{B}_{\varepsilon,i}&=\begin{pmatrix} 0&-\big[\big(w(1-w)v\big)(\varepsilon^2 D^{S_{\widetilde{M}},2})\varepsilon D^{S_{\widetilde{M}}},f_i\big]\\
				\big[(wv)(\varepsilon^2 D^{S_{\widetilde{M}},2})\varepsilon D^{S_{\widetilde{M}}},f_i\big]&0
			\end{pmatrix},\\
			\mathcal{B}_{\varepsilon,i}'&=\begin{pmatrix} 0&-\big(w(1-w)v\big)(\varepsilon^2 D^{S_{\widetilde{M}},2})\big[\varepsilon D^{S_{\widetilde{M}}},f_i\big]\\
				(wv)(\varepsilon^2 D^{S_{\widetilde{M}},2})\big[\varepsilon D^{S_{\widetilde{M}}},f_i\big]&0
			\end{pmatrix},\\
			\mathcal{B}_{\varepsilon,i}''&=\begin{pmatrix} 0&-\big[\big(w(1-w)v\big)(\varepsilon^2 D^{S_{\widetilde{M}},2}),f_i\big]\varepsilon D^{S_{\widetilde{M}}}\\
				\big[(wv)(\varepsilon^2 D^{S_{\widetilde{M}},2}),f_i\big]\varepsilon D^{S_{\widetilde{M}}}&0
			\end{pmatrix},
		\end{split}
	\end{equation}
	where $D^{S_{\widetilde{M}},2}$ denotes the square of the Dirac operator $D^{S_{\widetilde{M}}}$. Then by \eqref{5.3} and \eqref{5.9}, we have
	\begin{equation}\label{5.11}
		\begin{split}
			\mathcal{K}_\varepsilon&=\tau\big(\mathcal{A}_\varepsilon+\mathcal{B}_\varepsilon\big)\tau f_0\prod_{i=1}^{2k}\big(\mathcal{A}_{\varepsilon,i}+\mathcal{B}_{\varepsilon,i}\big)\tau\\
			&=\tau\big(\mathcal{A}_\varepsilon+\mathcal{B}_\varepsilon\big)\tau f_0\prod_{i=1}^{2k}\big(\mathcal{A}_{\varepsilon,i}+\mathcal{B}_{\varepsilon,i}'+\mathcal{B}_{\varepsilon,i}''\big)\tau.
		\end{split}
	\end{equation}

	\subsection{Trivialization}\label{LIC2}

	Let us fix $x_0\in \widetilde{M}$. For ${\varepsilon}>0$, let $B^{\widetilde{M}}(x_0,\mathrm{inj}(\widetilde{M}))$ be the open balls in $\widetilde{M}$ with center $x_0$ and radius $\mathrm{inj}(\widetilde{M})$, and $B^{T_{x_0}\widetilde{M}}(0,\mathrm{inj}(\widetilde{M}))$ the open balls in $T_{{x_0}}\widetilde{M}$ with center $0$ and radius $\mathrm{inj}(\widetilde{M})$. The exponential map $\exp_{x_0}^{\widetilde{M}}\colon B^{T_{x_0}\widetilde{M}}(0,\mathrm{inj}(\widetilde{M}))\to B^{\widetilde{M}}(x_0,\mathrm{inj}(\widetilde{M}))$ is a diffeomorphism, and this gives local coordinates by identifying $T_{x_0}\widetilde{M}$ with $\mathbb{R}^{m}$ via an orthonormal basis $\{e_{i}\}_{i=1}^m$ of $T_{x_0}\widetilde{M}$, in other words,
	\begin{equation}\label{fb1}
		X=(X_1,\cdots,X_m)\in\mathbb{R}^{m} \mapsto  X_i{e}_{i}\in T_{x_0}\widetilde{M} \longmapsto \exp^{\widetilde{M}}_{{x_0}}(X_i{e}_{i}).
	\end{equation}
	We shall always identify $B^{\widetilde{M}}(x_0,\mathrm{inj}(\widetilde{M}))$ with $B^{T_{x_0}\widetilde{M}}(0,\mathrm{inj}(\widetilde{M}))$ through \eqref{fb1}. Using the parallel transport associated with $\nabla^{T\widetilde{M}}$ along the curve $s\in[0,1]\mapsto sX$, we extend $\{e_i\}_{i=1}^m$ to a local orthonormal frame $\{\widetilde{e}_i\}_{i=1}^m$ of $T\widetilde{M}$ over $B^{T_{x_0}\widetilde{M}}(0,\mathrm{inj}(\widetilde{M}))$.

	For $X\in B^{T_{x_0}\widetilde{M}}(0,\mathrm{inj}(\widetilde{M}))$, we identify $S_{\widetilde{M},X}$ with $S_{\widetilde{M},x_0}$ by parallel transport corresponding to $\nabla^{S_{\widetilde{M}}}$ along the curve $s\in[0,1]\mapsto sX$. This identification also preserves the metric $h^{S_{\widetilde{M}}}$ and the Clifford multiplication $c(\widetilde{e}^i)$, therefore,
	\begin{equation}
		\Big(S_{\widetilde{M}},h^{S_{\widetilde{M}}},c(\widetilde{e}^i)\Big)\Big|_{B^{T_{x_0}\widetilde{M}}(0,\mathrm{inj}(\widetilde{M}))}\cong \Big(S_{\widetilde{M},x_0},h^{S_{\widetilde{M},x_0}},c(e^i)\Big)\Big|_{B^{T_{x_0}\widetilde{M}}(0,\mathrm{inj}(\widetilde{M}))}.
	\end{equation}
	Let $\theta^{S_{\widetilde{M}}}$ be the associated connection forms of $\nabla^{S_{\widetilde{M}}}$ on $B^{T_{x_0}\widetilde{M}}(0,\mathrm{inj}(\widetilde{M}))$, then
	\begin{equation}
		\begin{split}
			\nabla^{S_{\widetilde{M}}}&=d+\theta^{S_{\widetilde{M}}}\\
			&=d+\frac{1}{4}\langle \theta^{S_{\widetilde{M}}}\widetilde{e}_k,\widetilde{e}_\ell\rangle c(e^k)c(e^\ell).
		\end{split}
	\end{equation}

	We denote $\widetilde{M}_0=T_{x_0}\widetilde{M}\cong\mathbb{R}^m$, and we construct a metric $g^{T\widetilde{M}_0}$ on $T\widetilde{M}_0$ such that
	\begin{equation}\label{fb.73}
		g^{T\widetilde{M}_0}=\begin{cases}
			g^{T\widetilde{M}}, &\text{ over } B^{T_{x_0}\widetilde{M}}(0,\mathrm{inj}(X)/2),\\
			g^{T_{x_0}\widetilde{M}}, &\text{ outside } B^{T_{x_0}\widetilde{M}}(0,\mathrm{inj}(X)),
		\end{cases}
	\end{equation}
	and let $dv_{\widetilde{M}_0}$ be the associated volume form. Let $dv_{T_{x_0}\widetilde{M}}$ be the Riemannian volume form of $(T_{x_0}\widetilde{M},g^{T_{x_0}\widetilde{M}})$ and $\kappa$ the smooth positive function on $\widetilde{M}_0$ defined by
	\begin{equation}
		dv_{\widetilde{M}_0}(X)=\kappa(X)dv_{T_{x_0}\widetilde{M}},
	\end{equation}
	where $\kappa\lk0\rk=1$.
	
	Let $\rho(X)$ be a radially increasing smooth function on $T_{x_0}\widetilde{M}$ such that
	\begin{equation}\label{5.18'}
		\rho(X)=\begin{cases}
			1, &\text{ over } B^{T_{x_0}\widetilde{M}}(0,\mathrm{inj}(X)/2),\\
			0, &\text{ outside } B^{T_{x_0}\widetilde{M}}(0,\mathrm{inj}(X)).
		\end{cases}
	\end{equation}
	Let us extend $\nabla^{S_{\widetilde{M}}}$ on $S_{\widetilde{M}}$ over $B^{T_{x_0}\widetilde{M}}(0,r)$ to $\nabla^{S_{\widetilde{M},x_0}}$ on $S_{\widetilde{M},x_0}$ over $T_{x_0}\widetilde{M}$ by
	\begin{equation}
		\begin{split}
			\nabla^{S_{\widetilde{M},x_0}}&=d+\rho(X)\theta^{S_{\widetilde{M}}}_X\\
			&=d+\frac{1}{4}\rho(X)\big\langle \theta^{S_{\widetilde{M}}}_X\widetilde{e}_k,\widetilde{e}_\ell\big\rangle c(e^k)c(e^\ell)
		\end{split}
	\end{equation}

	Now we set
	\begin{equation}
		\begin{split}
			f_i^{x_0}&=\rho(X)f_i(X),\ \ \ 	D^{x_0}=c(\widetilde{e}^i)\nabla^{S_{\widetilde{M},x_0}}_{\widetilde{e}_i},\\
			D^{2,x_0}&=g^{T\widetilde{M}_0,ij}\Big(\nabla^{S_{\widetilde{M},x_0}}_{\pa_{X_i}}\nabla^{S_{\widetilde{M},x_0}}_{\pa_{X_j}}-\nabla^{S_{\widetilde{M},x_0}}_{\nabla^{T\widetilde{M}_0}_{\pa_{X_i}}\pa_{X_j}}\Big)+\frac{1}{4}\rho(X)\mathrm{Sc}_{g^{T\widetilde{M}}}(X),
		\end{split}
	\end{equation}
	then by the following expression for $D^{S_{\widetilde{M}}}$ and the classical Lichnerowicz formula for $D^{S_{\widetilde{M}},2}$
	\begin{equation}
		D^{S_{\widetilde{M}}}=c(\widetilde{e}^i)\nabla^{S_{\widetilde{M}}}_{\widetilde{e}_i},\ \ \ 	D^{S_{\widetilde{M}},2}=\Delta^{S_{\widetilde{M}}}+\frac{1}{4}\mathrm{Sc}_{g^{T\widetilde{M}}},
	\end{equation}
	we have
	\begin{equation}\label{5.20}
		\Big(f_i^{x_0},D^{x_0},D^{2,x_0}\Big)\Big|_{B^{T_{x_0}\widetilde{M}}(0,\mathrm{inj}(\widetilde{M})/2)}=\Big(f_i,D^{S_{\widetilde{M}}},D^{S_{\widetilde{M}},2},\Big)\Big|_{B^{T_{x_0}\widetilde{M}}(0,\mathrm{inj}(\widetilde{M})/2)}.
	\end{equation}

	We put operators
	\begin{equation}
		\Big(\mathcal{A}_\varepsilon^{x_0},\mathcal{B}_\varepsilon^{x_0},\mathcal{A}_{\varepsilon,i}^{x_0},\mathcal{B}_{\varepsilon,i}^{x_0},\mathcal{B}_{\varepsilon,i}^{x_0}{}',\mathcal{B}_{\varepsilon,i}^{x_0}{}''\Big)
	\end{equation}
	acting on $C_0^\infty\big(T_{x_0}\widetilde{M},S_{\widetilde{M},x_0}\big)$ by replacing $(f_i,D^{S_{\widetilde{M}}},D^{S_{\widetilde{M}},2})$ with $(f_i^{x_0},D^{x_0},D^{2,x_0})$ in
	\begin{equation}
		\Big(\mathcal{A}_\varepsilon,\mathcal{B}_\varepsilon,\mathcal{A}_{\varepsilon,i},\mathcal{B}_{\varepsilon,i},\mathcal{B}_{\varepsilon,i}',\mathcal{B}_{\varepsilon,i}''\Big),
	\end{equation}
	then from \eqref{5.10}, we have
	\begin{equation}\label{5.23}
		\begin{split}
			\mathcal{A}_\varepsilon^{x_0}&=\begin{pmatrix} w^2\big(\varepsilon^2D^{2,x_0}\big)&0\\
				0&w^2\big(\varepsilon^2D^{2,x_0}\big)
			\end{pmatrix},\\
			\mathcal{B}_\varepsilon^{x_0}&=\begin{pmatrix} 0&-\big(w(1-w)v\big)\big(\varepsilon^2D^{2,x_0}\big)\varepsilon D^{x_0}\\
				(wv)\big(\varepsilon^2D^{2,x_0}\big)\varepsilon D^{x_0}&0
			\end{pmatrix},\\
			\mathcal{A}_{\varepsilon,i}^{x_0}&=\begin{pmatrix} \big[w^2(\varepsilon^2D^{2,x_0}),f_i^{x_0}\big]&0\\
				0&\big[w^2(\varepsilon^2D^{2,x_0}),f_i^{x_0}\big]
			\end{pmatrix},\\
			\mathcal{B}_{\varepsilon,i}^{x_0}&=\begin{pmatrix} 0&-\big[\big(w(1-w)v\big)(\varepsilon^2D^{2,x_0})\varepsilon D^{x_0},f_i^{x_0}\big]\\
				\big[(wv)(\varepsilon^2D^{2,x_0})\varepsilon D^{x_0},f_i^{x_0}\big]&0
			\end{pmatrix},\\
			\mathcal{B}_{\varepsilon,i}^{x_0}{}'&=\begin{pmatrix} 0&-\big(w(1-w)v\big)(\varepsilon^2 D^{2,x_0})\big[\varepsilon D^{x_0},f_i^{x_0}\big]\\
				(wv)(\varepsilon^2 D^{2,x_0})\big[\varepsilon D^{x_0},f_i^{x_0}\big]&0
			\end{pmatrix},\\
			\mathcal{B}_{\varepsilon,i}^{x_0}{}''&=\begin{pmatrix} 0&-\big[\big(w(1-w)v\big)(\varepsilon^2 D^{2,x_0}),f_i^{x_0}\big]\varepsilon D^{x_0}\\
				\big[(wv)(\varepsilon^2 D^{2,x_0}),f_i^{x_0}\big]\varepsilon D^{x_0}&0
			\end{pmatrix},
		\end{split}
	\end{equation}
	where $v(\lambda),w(\lambda)$ are given in \eqref{5.10'}.

	Let us denote
	\begin{equation}\label{5.24}
		\begin{split}
			\mathcal{K}^{x_0}_\varepsilon=&\tau\big(\mathcal{A}_\varepsilon^{x_0}+\mathcal{B}_\varepsilon^{x_0}\big)\tau f_0^{x_0}\prod_{i=1}^{2k}\big(\mathcal{A}_{\varepsilon,i}^{x_0}+\mathcal{B}_{\varepsilon,i}^{x_0}\big)\tau\\
			=&\tau\big(\mathcal{A}_\varepsilon^{x_0}+\mathcal{B}_\varepsilon^{x_0}\big)\tau f_0^{x_0}\prod_{i=1}^{2k}\big(\mathcal{A}_{\varepsilon,i}^{x_0}+\mathcal{B}_{\varepsilon,i}^{x_0}{}'+\mathcal{B}_{\varepsilon,i}^{x_0}{}''\big)\tau,
		\end{split}
	\end{equation}
	and $\mathcal{K}^{x_0}_\varepsilon(X,X')$ the kernel of $\mathcal{K}^{x_0}_\varepsilon$ with respect to $dv_{M_0}$. By \eqref{3.4}, \eqref{5.3}, \eqref{5.9}, \eqref{5.11}, \eqref{5.20}, \eqref{5.23} and \eqref{5.24}, when $0<\varepsilon\leqslant\mathrm{inj}(\widetilde{M})/2$ we have
	\begin{equation}
		\mathcal{K}^{x_0}_\varepsilon(0,\cdot)=\mathcal{K}_\varepsilon(x_0,\cdot),
	\end{equation}
	and in particular,
	\begin{equation}\label{5.26}
		\mathcal{K}^{x_0}_\varepsilon(0,0)=\mathcal{K}_\varepsilon(x_0,x_0).
	\end{equation}

	\subsection{Rescaling}\label{LIC3}

	Let us define the space rescaling operator $S_\varepsilon$ by
	\begin{equation}\label{5.29'}
		S_\varepsilon(u)(X)=u(\varepsilon X),
	\end{equation}
	then we easily check that $S_{\varepsilon}\mathcal{K}^{x_0}_\varepsilon S_{\varepsilon^{-1}}$ is also a kernel operator with
	\begin{equation}\label{5.28}
		S_{\varepsilon}\mathcal{K}^{x_0}_\varepsilon S_{\varepsilon^{-1}}(0,0)=\varepsilon^m\mathcal{K}_\varepsilon^{x_0}(0,0).
	\end{equation}

	We now introduce an additional Clifford rescaling on $\mathcal{K}_\varepsilon^{x_0}$ following Getzler \cite{MR836727}.

	By \cite[Proposition 3.21]{MR2273508}, among all the monoids $c(e^{i_1})\cdots c(e^{i_j})$, only $c(e^1)\cdots c(e^m)$ has a nonzero supertrace and
	\begin{equation}\label{5.29}
		\mathrm{Tr}_s^{S_{\widetilde{M}}}\big[c(e^1)\cdots c(e^m)\big]=(-2\sqrt{-1})^{m/2}.
	\end{equation}
	We set
	\begin{equation}\label{5.32'}
		c_\varepsilon(e^i)=\varepsilon^{-1}(e^i\wedge)-\varepsilon\iota_{e_i}
	\end{equation}
	acting on $\Lambda(T_{x_0}^*\widetilde{M})$, where $\wedge$ means exterior product and $\iota$ denotes the interior product.. For any $Q\in \mathrm{End}\big(\Lambda(T_{x_0}^*\widetilde{M})\big)$, we can express it in the form of
	\begin{equation}
		Q=\sum_{
			1 \leqslant i_1<\cdots< i_k\leqslant m,\ 
			1 \leqslant j_1<\cdots< j_\ell\leqslant m}Q^{j_1,\cdots,j_\ell}_{i_1,\cdots,i_k}e^{i_1}\wedge\cdots\wedge e^{i_k}\wedge i_{e_{j_1}}\cdots i_{e_{j_\ell}},
	\end{equation}
	and we set
	\begin{equation}
		[Q]^{\max}=(-2\sqrt{-1})^{m/2}Q_{1,\cdots,m}.
	\end{equation}
	Then from \eqref{5.29} we have for any monoid,
	\begin{equation}\label{5.33}
		\mathrm{Tr}_s^{S_{\widetilde{M}}}\big[c(e^{i_1})\cdots c(e^{i_k})\big]=\varepsilon^m\big[c_\varepsilon(e^{i_1})\cdots c_\varepsilon(e^{i_k})\big]^{\max}.
	\end{equation}
	

	Getzler rescaling is the combination of the space rescaling \eqref{5.29'} and the Clifford rescaling \eqref{5.32'}. Let us define the Getzler rescaling operators
	\begin{equation}\label{5.34}
		\Big(f_\varepsilon^{x_0},d_\varepsilon^{x_0},\Delta_\varepsilon^{x_0},\mathscr{A}_{\varepsilon}^{x_0},\mathscr{B}_{\varepsilon}^{x_0},\mathscr{A}_{\varepsilon,i}^{x_0},\mathscr{B}_{\varepsilon,i}^{x_0},\mathscr{B}_{\varepsilon,i}^{x_0}{}',\mathscr{B}_{\varepsilon,i}^{x_0}{}''\Big)
	\end{equation}
	acting on $C_0^\infty\big(T_{x_0}\widetilde{M},\Lambda\big(T^*_{x_0}\widetilde{M}\big)\big)$, by replacing each $c(e^i)$ with $c_\varepsilon(e^i)$ in
	\begin{equation}
		S_{\varepsilon}\Big(f^{x_0},\varepsilon^2 D^{x_0},\varepsilon^2 D^{2,x_0},\mathcal{A}_{\varepsilon}^{x_0},\mathcal{B}_{\varepsilon}^{x_0},\mathcal{A}_{\varepsilon,i}^{x_0},\mathcal{B}_{\varepsilon,i}^{x_0},\mathcal{B}_{\varepsilon,i}^{x_0}{}',\mathcal{B}_{\varepsilon,i}^{x_0}{}''\Big)S_{\varepsilon^{-1}},
	\end{equation}
	where the adjoint of $S_\varepsilon$ acts on each term. Then by evaluating the tensorial terms on the left sides at $X$ and right sides at $\varepsilon X$, we have
	\begin{equation}\label{5.36}
		\begin{split}
			f^{x_0}_\varepsilon=f^{x_0},\ \ \ \ \ \ &\\
			d_\varepsilon^{x_0}=\varepsilon c_\varepsilon(e^i)\Big(&\widetilde{e}_i+\frac{1}{4}\varepsilon\rho(\varepsilon X)\langle \theta^{T\widetilde{M}}_{\varepsilon X}(\widetilde{e}_i)\widetilde{e}_k,\widetilde{e}_\ell\rangle c_\varepsilon(e^k)c_\varepsilon(e^\ell)\Big),\\
			\Delta_\varepsilon^{x_0}=g^{T\widetilde{M}_0,ij}_{\varepsilon X}\bigg[&\Big(\pa_{X_i}+\frac{1}{4}\varepsilon\rho(\varepsilon X)\langle \theta^{T\widetilde{M}}_{\varepsilon X}(\pa_{X_i})\widetilde{e}_k,\widetilde{e}_\ell\rangle c_\varepsilon(e^k)c_\varepsilon(e^\ell)\Big)\\
			&\cdot\Big(\pa_{X_j}+\frac{1}{4}\varepsilon\rho(\varepsilon X)\langle \theta^{T\widetilde{M}}_{\varepsilon X}(\pa_{X_j})\widetilde{e}_k,\widetilde{e}_\ell\rangle c_\varepsilon(e^k)c_\varepsilon(e^\ell)\Big)\\
			&-\varepsilon\Big(\nabla^{T_x\widetilde{M}}_{\pa_{X_i}}\pa_{X_j}+\frac{1}{4}\varepsilon\rho(\varepsilon X)\Big\langle \theta^{T\widetilde{M}}_{\varepsilon X}\big(\nabla^{T_xM}_{\pa_{X_i}}\pa_{X_j}\big)\widetilde{e}_k,\widetilde{e}_\ell\Big\rangle c_\varepsilon(e^k)c_\varepsilon(e^\ell)\Big)\bigg]\\
			+\frac{1}{4}\varepsilon^2&\rho(\varepsilon X)\mathrm{Sc}_{g^{T\widetilde{M}}}(\varepsilon X).
		\end{split}
	\end{equation}
	Moreover,
	\begin{equation}\label{5.37}
		\begin{split}
			\mathscr{A}_\varepsilon^{x_0}&=\begin{pmatrix} w^2(\Delta_\varepsilon^{x_0})&0\\
				0&w^2(\Delta_\varepsilon^{x_0})
			\end{pmatrix},\\
			\mathscr{B}_\varepsilon^{x_0}&=\frac{1}{\varepsilon}\begin{pmatrix} 0&-\big(w(1-w)v\big)(\Delta_\varepsilon^{x_0})d^{x_0}_\varepsilon\\
				(wv)(\Delta_\varepsilon^{x_0})d^{x_0}_\varepsilon&0
			\end{pmatrix},\\
			\mathscr{A}_{\varepsilon,i}^{x_0}&=\begin{pmatrix} \big[w^2(\Delta_\varepsilon^{x_0}), f_{\varepsilon,i}^{x_0}\big]&0\\
				0&\big[w^2(\Delta_0^{x_0}),f_{\varepsilon,i}^{x_0}\big]
			\end{pmatrix},\\
			\mathscr{B}_{\varepsilon,i}^{x_0}&=\begin{pmatrix} 0&-\big[\big(w(1-w)v\big)(\Delta_\varepsilon^{x_0})\varepsilon^{-1}d_\varepsilon^{x_0},f_{\varepsilon,i}^{x_0}\big]\\
				\big[(wv)(\Delta_\varepsilon^{x_0})\varepsilon^{-1}d_\varepsilon^{x_0},f_{\varepsilon,i}^{x_0}\big]&0
			\end{pmatrix},\\
			\mathscr{B}_{\varepsilon,i}^{x_0}{}'&=\begin{pmatrix} 0&-\big(w(1-w)v\big)(\Delta_\varepsilon^{x_0})[\varepsilon^{-1}d^{x_0}_\varepsilon,f_{\varepsilon,i}^{x_0}]\\
				(wv)(\Delta_\varepsilon^{x_0})[\varepsilon^{-1}d^{x_0}_\varepsilon,f_{\varepsilon,i}^{x_0}]&0
			\end{pmatrix},\\
			\mathscr{B}_{\varepsilon,i}^{x_0}{}''&=\begin{pmatrix} 0&-\big[\big(w(1-w)v\big)(\Delta_\varepsilon^{x_0}),f_{\varepsilon,i}^{x_0}\big]\varepsilon^{-1}d^{x_0}_\varepsilon\\
				\big[(wv)(\Delta_\varepsilon^{x_0}),f_{\varepsilon,i}^{x_0}\big]\varepsilon^{-1}d^{x_0}_\varepsilon&0
			\end{pmatrix}.
		\end{split}
	\end{equation}

	Similar to \eqref{5.24}, we denote
	\begin{equation}\label{5.38}
		\begin{split}
			\mathscr{K}_\varepsilon^{x_0}&=\tau\big(\mathscr{A}_\varepsilon^{x_0}+\mathscr{B}_\varepsilon^{x_0}\big)\tau f_0^{x_0}\prod_{i=1}^{2k}\big(\mathscr{A}_{\varepsilon,i}^{x_0}+\mathscr{B}_{\varepsilon,i}^{x_0}\big)\tau\\
			&=\tau\big(\mathscr{A}_\varepsilon^{x_0}+\mathscr{B}_\varepsilon^{x_0}\big)\tau f_0^{x_0}\prod_{i=1}^{2k}\big(\mathscr{A}_{\varepsilon,i}^{x_0}+\mathscr{B}_{\varepsilon,i}^{x_0}{}'+\mathscr{B}_{\varepsilon,i}^{x_0}{}''\big)\tau.
		\end{split}
	\end{equation}
	From \eqref{5.33} we obtain
	\begin{equation}
		\begin{split}
			\mathrm{Tr}^{S_{\widetilde{M}}^{\oplus2}}\Big(S_{\varepsilon}\mathcal{K}_\varepsilon^{x_0}S_{\varepsilon^{-1}}(0,0)\Big)=\varepsilon^m\big(\mathscr{K}_\varepsilon^{x_0}(0,0)\big)^{\text{max}},
		\end{split}
	\end{equation}
	which, together with \eqref{5.38} and \eqref{5.28}, implies
	\begin{equation}\label{5.40}
		\mathrm{Tr}^{S_{\widetilde{M}}^{\oplus2}}\big(\mathcal{K}_\varepsilon^{}(x_0,x_0)\big)=\big(\mathscr{K}_\varepsilon^{x_0}\big(0,0)\big)^{\text{max}}.
	\end{equation}
	As indicated by \eqref{5.8} \eqref{5.9} and \eqref{5.40}, the Getzler rescaling transfers the Connes-Chern character to the kernel of the associated rescaled operator.

	\subsection{Limit operators}\label{LIC4}

	We now analyze the asymptotic behavior of the operators in \eqref{5.36}, \eqref{5.37} and \eqref{5.38} as $\varepsilon \to 0$. We note that the asymptotic expansions in this section are relatively intuitive, and a rigorous treatment, including precise estimates for all remainder terms, will be given in \cref{FAM}.

	\begin{lemma}
		We have
		\begin{equation}\label{5.41}
			\begin{split}
				f_\varepsilon^{x_0}&=f(x_0)+\varepsilon \nabla f(x_0)\cdot X+O(\varepsilon^2),\\
				d_\varepsilon^{x_0}&=d+O(\varepsilon),\\
				\Delta_\varepsilon^{x_0}&=-\Big(\pa_{X_i}+\frac{1}{4}X_j\big(R^{T\widetilde{M}}_{x_0}\pa_{X_j},\pa_{X_i}\big)\Big)^2+O(\varepsilon).
			\end{split}
		\end{equation}
	\end{lemma}

	\begin{pro}
		By \eqref{5.36}, we have
		\begin{equation}
			f_\varepsilon^{x_0}(X)=f^{x_0}(\varepsilon X),
		\end{equation}
		which implies the first equality in \eqref{5.41}
		
		By \cite[Proposition 1.18]{MR2273508}, we have
		\begin{equation}
			\theta^{T\widetilde{M}}_X(\cdot)=\frac{1}{2}X_jR^{T\widetilde{M}}_{x_0}(\pa_{X_j},\cdot)+O(X^2),
		\end{equation}
		which, together with \eqref{5.36}, implies
		\begin{equation}
			\begin{split}
				\varepsilon\rho(\varepsilon X)\langle \theta^{T\widetilde{M}}_{\varepsilon X}(\widetilde{e}_i)\widetilde{e}_k,\widetilde{e}_\ell\rangle c_\varepsilon(e^k)c_\varepsilon(e^\ell)&=\frac{1}{2}X_j\big\langle R^{T\widetilde{M}}_{x_0}(\pa_{X_{j}},\pa_{X_i})\pa_{X_k},\pa_{X_\ell}\big\rangle dX_k\wedge dX_\ell+O(\varepsilon)\\
				&=X_jR^{T\widetilde{M}}_{x_0}(\pa_{X_{j}},\pa_{X_i})+O(\varepsilon).
			\end{split}
		\end{equation}
		Then by \eqref{5.36}, the identity
		\begin{equation}\label{5.45}
			\begin{split}
				X_jdX_i\wedge\big(R^{T\widetilde{M}}\pa_{X_j},\pa_{X_i}\big)=-X_j\big(R^{T\widetilde{M}}(\pa_{X_k},\pa_{X_\ell})\pa_{X_i},\pa_{X_j}\big)dX_i\wedge dX_k\wedge dX_\ell
			\end{split}
		\end{equation}
		and the Bianchi identity
		\begin{equation}\label{5.46}
			R^{T\widetilde{M}}(\pa_{X_k},\pa_{X_\ell})\pa_{X_i}+R^{T\widetilde{M}}(\pa_{X_i},\pa_{X_k})\pa_{X_\ell}+R^{T\widetilde{M}}(\pa_{X_\ell},\pa_{X_i})\pa_{X_k}=0,
		\end{equation}
		we get the second and third identity in \eqref{5.41}.\qed
	\end{pro}
	
	By \eqref{5.41}, let us denote
	\begin{equation}
		d_0^{x_0}=d,\ \ 
		\Delta_0^{x_0}=-\Big(\pa_{X_i}+\frac{1}{4}X_j\big(R^{T\widetilde{M}}_{x_0}\pa_{X_j},\pa_{X_i}\big)\Big)^2.
	\end{equation}
	Then combining \eqref{5.37} and \eqref{5.41}, we have
	\begin{equation}\label{5.47}
		\begin{split}
			\mathscr{A}_\varepsilon^{x_0}&=\begin{pmatrix} w^2(\Delta_0^{x_0})&0\\
				0&w^2(\Delta_0^{x_0})
			\end{pmatrix}+O(\varepsilon),\\
			\mathscr{B}_\varepsilon^{x_0}&=\varepsilon^{-1}\begin{pmatrix} 0&-\big(w(1-w)v\big)(\Delta_0^{x_0})d\\
				(wv)(\Delta_0^{x_0})d&0
			\end{pmatrix}+O(1),\\
			\mathscr{A}_{\varepsilon,i}^{x_0}&=\varepsilon\begin{pmatrix} \big[w^2(\Delta_0^{x_0}),\nabla f_i(0)\cdot X\big]&0\\
				0&\big[w^2(\Delta_0^{x_0}),\nabla f_i(0)\cdot X\big]
			\end{pmatrix}+O(\varepsilon^2),\\
			\mathscr{B}_{\varepsilon,i}^{x_0}&=\begin{pmatrix} 0&-\big[\big(w(1-w)v\big)(\Delta_0^{x_0})d,\nabla f_i(0)\cdot X\big]\\
				\big[(wv)(\Delta_0^{x_0})d,\nabla f_i(0)\cdot X\big]&0
			\end{pmatrix}+O(\varepsilon),\\
			\mathscr{B}_{\varepsilon,i}^{x_0}{}'&=\begin{pmatrix} 0&-\big(w(1-w)v\big)(\Delta_0^{x_0})df_i(0)\wedge\\
				(wv)(\Delta_0^{x_0})df_i(0)\wedge&0
			\end{pmatrix}+O(\varepsilon),\\
			\mathscr{B}_{\varepsilon,i}^{x_0}{}''&=\begin{pmatrix} 0&-\big[\big(w(1-w)v\big)(\Delta_0^{x_0}),\nabla f_i(0)\cdot X\big]d\\
				\big[(wv)(\Delta_0^{x_0}),\nabla f_i(0)\cdot X\big]d&0
			\end{pmatrix}+O(\varepsilon),
		\end{split}
	\end{equation}
	where we use the fact that $[d,\nabla f(0)\cdot X]=df(0)\wedge$. Use \eqref{5.47} to count the order of $\varepsilon$ in $\mathscr{K}_\varepsilon^{x_0}$ as defined in \eqref{5.38}, the only divergent term with order $\varepsilon^{-1}$ is
	\begin{equation}
		\tau\mathscr{B}_\varepsilon^{x_0}\tau f_0(0)\prod_{i=1}^{2k}\big(\mathscr{B}_{\varepsilon,i}^{x_0}\tau\big),
	\end{equation}
	but this operator is off-diagonal, so it does not contribute to the super-trace. Hence, we only need to compute the constant terms with order $0$, which are given by
	\begin{equation}\label{5.49}
		\tau\mathscr{A}_0^{x_0}\tau f_0(0)\prod_{i=1}^{2k}\mathscr{B}_{0,i}^{x_0}\tau+\sum_{i=1}^{2k}\tau\mathscr{B}_0^{x_0}\tau f_0(0)\mathscr{B}_{0,1}^{x_0}\tau\cdots\mathscr{A}_{0,i}^{x_0}\tau\cdots\tau\mathscr{B}_{0,2k}^{x_0}\tau.
	\end{equation}
	Since $\mathscr{A}_{0}^{x_0},\mathscr{A}_{0,i}^{x_0}$ are even and $\mathscr{B}_{0}^{x_0},\mathscr{B}_{0,i}^{x_0}$ are odd, we can write \eqref{5.49} as
	\begin{equation}
		(-1)^k\mathscr{A}_0^{x_0} f_0(0)\prod_{i=1}^{2k}\mathscr{B}_{0,i}^{x_0}+\sum_{i=1}^{2k}(-1)^{k+i}\mathscr{B}_0^{x_0}f_0(0)\mathscr{B}_{0,1}^{x_0}\cdots\mathscr{A}_{0,i}^{x_0}\cdots\mathscr{B}_{0,2k}^{x_0}.
	\end{equation}
	Then from \eqref{5.40} we obtain
	\begin{equation}\label{5.51}
		\begin{split}
			&\lim_{\varepsilon\to0}\mathrm{Tr}^{S_{\widetilde{M}}^{\oplus2}}\big(\mathcal{K}_\varepsilon(x_0,x_0)\big)\\
			&=(-1)^k\Big[\mathscr{A}_0^{x_0} f_0(0)\prod_{i=1}^{2k}\big(\mathscr{B}_{0,i}^{x_0}{}'+\mathscr{B}_{0,i}^{x_0}{}''\big)(0,0)\Big]^{\text{max}}\\
			&\ \ \ \ +\sum_{i=1}^{2k}(-1)^{k+i}\Big[\mathscr{B}_0^{x_0}f_0(0)\big(\mathscr{B}_{0,1}^{x_0}{}'+\mathscr{B}_{0,1}^{x_0}{}''\big)\cdots\mathscr{A}_{0,i}^{x_0}\cdots\big(\mathscr{B}_{0,2k}^{x_0}{'}+\mathscr{B}_{0,2k}^{x_0}{''}\big)(0,0)\Big]^{\text{max}}.
		\end{split}
	\end{equation}

	
	
	We proceed to further simplify the product in \eqref{5.51}.
	
	\begin{lemma}
		We have
		\begin{equation}\label{5.52}
			\begin{split} \big[d,\Delta_0^{x_0}\big]&=0,\\
				\big[d,\big[\Delta_0^{x_0},\nabla f(0)\cdot X\big]\big]&=0.
			\end{split}
		\end{equation}
	\end{lemma}
	
	\begin{pro}
		In addition to verifying \eqref{5.52} by straightforward computation, we may alternatively apply the rescaling procedure described in \eqref{5.34} to the following identities
		\begin{equation}
			\begin{split}
				\big[\varepsilon^2 D^{S_{\widetilde{M}},2},\varepsilon^2D^{S_{\widetilde{M}}}\big]&=0,\\
				\big[\varepsilon^2D^{S_{\widetilde{M}}},\big[\varepsilon^2D^{S_{\widetilde{M}},2},\varepsilon^{-1}f\big]\big]&=\big[\varepsilon^2D^{S_{\widetilde{M}},2},\big[\varepsilon^2D^{S_{\widetilde{M}}},\varepsilon^{-1}f\big]\big],
			\end{split}
		\end{equation}
		from which \eqref{5.52} follows.\qed
	\end{pro}

	Let us list a series of vanishing terms in the expansion of \eqref{5.51}. If a term starts with $\mathscr{A}_0^{x_0}$ and contains at least two factors $\mathscr{B}_{0,i}^{x_0}{''}$ and $\mathscr{B}_{0,j}^{x_0}{''}$, since, we can move the two $d$'s together by \eqref{5.52}, and since $d^2=0$, it vanishes. If a term starts with $\mathscr{B}_{0}^{x_0}$ and contains at least one factor $\mathscr{B}_{0,i}^{x_0}{''}$, similarly, we can bring the two $d$'s together and it vanishes. Summarizing the above results, we obtain
	\begin{equation}\label{5.54}
		\begin{split}
			\lim_{\varepsilon\to0}\mathrm{Tr}^{S_{\widetilde{M}}^{\oplus2}}\big(\mathcal{K}_\varepsilon(x_0,x_0)\big)=&(-1)^k\Big[\mathscr{A}_0^{x_0} f_0(0)\prod_{i=1}^{2k}\mathscr{B}_{0,i}^{x_0}{}'(0,0)\Big]^{\text{max}}\\
			&+\sum_{i=1}^{2k}(-1)^k\Big[\mathscr{A}_0^{x_0} f_0(0)\mathscr{B}_{0,1}^{x_0}{}'\cdots \mathscr{B}_{0,i}^{x_0}{}''\cdots \mathscr{B}_{0,2k}^{x_0}{}'(0,0)\Big]^{\text{max}}\\
			& +\sum_{i=1}^{2k}(-1)^{k+i}\Big[\mathscr{B}_0^{x_0}f_0(0)\mathscr{B}_{0,1}^{x_0}{}'\cdots\mathscr{A}_{0,i}^{x_0}\cdots\mathscr{B}_{0,2k}^{x_0}{'}(0,0)\Big]^{\text{max}}.
		\end{split}
	\end{equation}

	\subsection{Explicit kernel computations}\label{LIC5}
	
	We now explicitly calculate the limit kernel \eqref{5.54}.

	Let us start with the first term in \eqref{5.54}. By \eqref{5.47}, we compute that
	\begin{equation}
		\begin{split}
			&(-1)^k\mathscr{A}_0^{x_0} f_0(0)\prod_{i=1}^{2k}\mathscr{B}_{0,i}^{x_0}{}'\\
			&=(-1)^k\begin{pmatrix} w^2(\Delta_0^{x_0})&0\\
				0&w^2(\Delta_0^{x_0})
			\end{pmatrix}f_0(0)
			\prod_{i=1}^{2k}\begin{pmatrix} 0&-\big(w(1-w)v\big)(\Delta_0^{x_0})\\
				(wv)(\Delta_0^{x_0})&0
			\end{pmatrix}df_i(0)\\
			&=\begin{pmatrix} T_{1,11}&0\\
				0&T_{1,22}
			\end{pmatrix}(f_0df_1\wedge\cdots \wedge df_{2k})(0),
		\end{split}
	\end{equation}
	where
	\begin{equation}
		T_{1,11}=T_{1,22}=\big(w^2(w^2(1-w)v^2)^{k}\big)(\Delta_0^{x_0}).
	\end{equation}
	Therefore, we may reduce to considering
	\begin{equation}\label{5.58}
		\big[2\big(w^2(w^2(1-w)v^2)^{k}\big)(\Delta_0^{x_0})(0,0)\cdot(f_0df_1\wedge\cdots \wedge df_{2k})(0)\big]^{\max}.
	\end{equation}

	We then compute the second term in \eqref{5.54}. By \eqref{5.47} and the Cauchy formula, we have
	\begin{equation}
		\begin{split}
			\mathscr{B}_{0,i}^{x_0}{}''=\frac{1}{2\pi\sqrt{-1}}\int&\begin{pmatrix} 0&-\big(w(1-w)v\big)(\lambda)\\
				(wv)(\lambda)&0
			\end{pmatrix}\\
			&(\lambda-\Delta_0^{x_0})^{-1}\big[\Delta_0^{x_0},\nabla f_i(0)\cdot X\big]d(\lambda-\Delta_0^{x_0})^{-1}d\lambda,
		\end{split}
	\end{equation}
	where the integral is over a suitable contour. Since we need the integral of supertrace in \eqref{5.8}, therefore, we can replace the second term in \eqref{5.54} with
	\begin{equation}
		\begin{split}
			\frac{(-1)^{k+i-1}}{2\pi\sqrt{-1}}\int&\big[\Delta_0^{x_0},\nabla f_i(0)\cdot X\big]d(\lambda-\Delta_0^{x_0})^{-1}\mathscr{B}_{0,i+1}^{x_0}{}'\cdots \mathscr{B}_{0,2k}^{x_0}{}'\\
			&\mathscr{A}_0^{x_0}f_0(0)\mathscr{B}_{0,1}^{x_0}{}'\cdots \mathscr{B}_{0,i-1}^{x_0}{}'\begin{pmatrix} 0&-\big(w(1-w)v\big)(\lambda)\\
				(wv)(\lambda)&0
			\end{pmatrix}(\lambda-\Delta_0^{x_0})^{-1}d\lambda\\
			=&(-1)^{k+i-1}\big[\Delta_0^{x_0},\nabla f_i(0)\cdot X\big]d\mathscr{B}_{0,i+1}^{x_0}{}'\cdots \mathscr{B}_{0,2k}^{x_0}{}'\\
			&\mathscr{A}_0^{x_0}f_0(0)\mathscr{B}_{0,1}^{x_0}{}'\cdots \mathscr{B}_{0,i-1}^{x_0}{}'\begin{pmatrix} 0&-\pa_\lambda\big(w(1-w)v\big)(\Delta_0^{x_0})\\
				\pa_\lambda(wv)(\Delta_0^{x_0})&0
			\end{pmatrix},
		\end{split}
	\end{equation}
	where $(-1)^{(i-1)(2k-(i-1))}=(-1)^{i-1}$. Applying Cauchy formula again, this expression is equal to
	\begin{equation}
		(-1)^{i-1}\big[\Delta_0^{x_0},\nabla f_i(0)\cdot X\big]d\begin{pmatrix} T_{2,11}&0\\
			0&T_{2,22}
		\end{pmatrix}\ \big(f_0df_{i+1}\wedge\cdots\wedge df_{2k}\wedge  df_1\wedge\cdots\wedge df_{i-1}\big)(0),
	\end{equation}
	where 
	\begin{equation}
		\begin{split}
			T_{2,11}&=\big(w(1-w)v\pa_\lambda(wv)w^2(w^2(1-w)v^2)^{k-1}\big)(\Delta_0^{x_0}),\\
			T_{2,22}&=\big(wv\pa_\lambda(w(1-w)v)w^2(w^2(1-w)v^2)^{k-1}\big)(\Delta_0^{x_0}),
		\end{split}
	\end{equation}
	or with $T_{2,11},T_{2,22}$ interchanged, depending on the parity of $i$. Thus, it suffices to consider
	\begin{equation}\label{5.63}
		\begin{split}
			(-1)^{i-1}\Big[\big[\Delta_0^{x_0},\nabla f_i(0)\cdot X\big]d\pa_\lambda \big(w^2(1-w)v^2w^2(w^2(1-w)v^2)^{k-1}\big)(\Delta_0^{x_0})(0,0)&\\
			\big(f_0df_1\wedge\cdots\widehat{df_i}\cdots\wedge df_{2k}\big)(0)&\Big]^{\max}.
		\end{split}
	\end{equation}

	To evaluate \eqref{5.58} and \eqref{5.63}, we will use two key tools, Mehler formula
	\begin{equation}\label{5.64}
		e^{-t\Delta_0^{x_0}}(X,0)=\frac{1}{(4\pi t)^{m/2}}\det{}^{1/2}\bigg(\frac{tR^{T\widetilde{M}}_{x_0}/2}{\sinh \big(tR^{T\widetilde{M}}_{x_0}/2\big)}\bigg)e^{-\frac{1}{4t}\big\langle tR^{T\widetilde{M}}_{x_0}/2\coth\big(tR^{T\widetilde{M}}_{x_0}/2\big)X,X\big\rangle},
	\end{equation}
	and Euler integral
	\begin{equation}\label{5.65}
		\frac{1}{\Gamma(n)}\int_0^\infty e^{-t\lambda}\lambda^{n-1}d\lambda=t^{-n}.
	\end{equation}

	By \eqref{5.64} and \eqref{5.65}, we have
	\begin{equation}\label{5.66}
		\begin{split}
			&\big[e^{-t\Delta_0^{x_0}}(0,0)(f_0df_1\wedge\cdots \wedge df_{2k})(0)\big]^{\max}\\
			&=\bigg[\frac{1}{(4\pi )^{m/2}}\det{}^{1/2}\bigg(\frac{R^{T\widetilde{M}}_{x_0}/2}{\sinh \big(R^{T\widetilde{M}}_{x_0}/2\big)}\bigg)(f_0df_1\wedge\cdots \wedge df_{2k})_{x_0}\bigg]^{\max}t^{-k}\\
			&=\bigg[\frac{1}{(4\pi )^{m/2}}\det{}^{1/2}\bigg(\frac{R^{T\widetilde{M}}_{x_0}/2}{\sinh \big(R^{T\widetilde{M}}_{x_0}/2\big)}\bigg)(f_0df_1\wedge\cdots \wedge df_{2k})_{x_0}\bigg]^{\max}\frac{1}{\Gamma(k)}\int_0^\infty e^{-t\lambda}\lambda^{k-1}d\lambda
		\end{split}
	\end{equation}
	From \eqref{5.66}, we see that \eqref{5.58} can be expressed by
	\begin{equation}\label{5.67}
		\begin{split}
			\big[2&\big(w^2(w^2(1-w)v^2)^{k}\big)(\Delta_0^{x_0})(0,0)\cdot(f_0df_1\wedge\cdots \wedge df_{2k})(0)\big]^{\max}\\
			=&\bigg[\frac{1}{(4\pi )^{m/2}}\det{}^{1/2}\bigg(\frac{R^{T\widetilde{M}}_{x_0}/2}{\sinh \big(R^{T\widetilde{M}}_{x_0}/2\big)}\bigg)(f_0df_1\wedge\cdots \wedge df_{2k})_{x_0}\bigg]^{\max}\\
			&\cdot\frac{1}{\Gamma(k)}\int_{0}^\infty2w^2(w^2(1-w)v^2)^{k}(\lambda)\lambda^{k-1}d\lambda.
		\end{split}
	\end{equation}

	Now we turn to find an integral expression for \eqref{5.63}, similar to \eqref{5.67}. By \eqref{5.41}, we have
	\begin{equation}
		\big[\Delta_0^{x_0},\nabla f_i(0)\cdot X\big]d=-2(\pa_{X_j}f_i)(0)\Big(\pa_{X_j}+\frac{1}{4}X_k\big(R^{T\widetilde{M}}_{x_0}\pa_{X_k},\pa_{X_j}\big)\Big)dX_\ell\wedge\pa_{X_\ell}.
	\end{equation}
	Then, by evaluating at $(0,0)$, we have
	\begin{equation}\label{5.69}
		\begin{split}
			&\big[\Delta_0^{x_0},\nabla f_i(0)\cdot X\big]d e^{-\frac{1}{4t}\big\langle tR^{T\widetilde{M}}_{x_0}/2\coth\big(tR^{T\widetilde{M}}_{x_0}/2\big)X,X\big\rangle}\\
			&=-2(\pa_{X_j}f_i)(0)dX_\ell\wedge\pa_{X_jX_\ell}^2e^{-\frac{1}{4t}\big\langle tR^{T\widetilde{M}}_{x_0}/2\coth\big(tR^{T\widetilde{M}}_{x_0}/2\big)X,X\big\rangle}\\
			&=2(\pa_{X_j}f_i)(0)dX_\ell\wedge\pa_{X_jX_\ell}^2\bigg(\frac{1}{4t}\Big\langle tR^{T\widetilde{M}}_{x_0}/2\coth \big(tR^{T\widetilde{M}}_{x_0}/2\big)X,X\Big\rangle\bigg)\\
			&=\frac{1}{t}(\pa_{X_j}f_i)(0)dX_\ell\wedge\Big\langle tR^{T\widetilde{M}}_{x_0}/2\coth \big(tR^{T\widetilde{M}}_{x_0}/2\big)\pa_{X_{\ell}},\pa_{X_j}\Big\rangle.
		\end{split}
	\end{equation}
	Similar to \eqref{5.45} and \eqref{5.46}, we get
	\begin{equation}
		dX_\ell\wedge\Big\langle tR^{T\widetilde{M}}_{x_0}/2\coth \big(tR^{T\widetilde{M}}_{x_0}/2\big)\pa_{X_{\ell}},\pa_{X_j}\Big\rangle=dX_j,
	\end{equation}
	hence \eqref{5.69} is equal to
	\begin{equation}\label{5.71}
		\frac{1}{t}(\pa_{X_j}f_i)(0)dX_j\wedge=\frac{1}{t}df_i(0)\wedge.
	\end{equation}
	By \eqref{5.64}, \eqref{5.69} and \eqref{5.71}, we get
	\begin{equation}\label{5.72}
		\begin{split}
			&(-1)^{i-1}\Big[\big[\Delta_0^{x_0},\nabla f_i(0)\cdot X\big]de^{-t\Delta_0^{x_0}}(0,0)\big(f_0df_1\wedge\cdots\widehat{df_i}\cdots\wedge df_{2k}\big)(0)\Big]^{\max}\\
			&=\bigg[\frac{1}{(4\pi )^{m/2}}\det{}^{1/2}\bigg(\frac{R^{T\widetilde{M}}_{x_0}/2}{\sinh \big(R^{T\widetilde{M}}_{x_0}/2\big)}\bigg)(f_0df_1\wedge\cdots \wedge df_{2k})_{x_0}\bigg]^{\max}t^{-k-1}.
		\end{split}
	\end{equation}
	By \eqref{5.65}, \eqref{5.72}, and the fact that
	$\Gamma(k+1)=k\Gamma(k)$, the summing of \eqref{5.63} for $i$ can be expressed by
	\begin{equation}\label{5.73}
		\begin{split}
			\sum_i(-1)^{i-1}\Big[\big[\Delta_0^{x_0},\nabla f_i(0)\cdot X\big]d\pa_\lambda \big(w^2(1-w)v^2w^2(w^2(1-w)v^2)^{k-1}\big)(\Delta_0^{x_0})(0,0)&\\
			\big(f_0df_1\wedge\cdots\widehat{df_i}\cdots\wedge df_{2k}\big)(0)&\Big]^{\max}\\
			=\bigg[\frac{1}{(4\pi )^{m/2}}\det{}^{1/2}\bigg(\frac{R^{T\widetilde{M}}_{x_0}/2}{\sinh\big(R^{T\widetilde{M}}_{x_0}/2\big)}\bigg)(f_0df_1\wedge\cdots \wedge df_{2k})_{x_0}&\bigg]^{\max}\\
			\cdot\frac{1}{\Gamma(k)}\int_{0}^\infty2\pa_\lambda (w^2(1-w)v^2)w^2(w^2(1-w)v^2)^{k-1}(\lambda)\lambda^k&d\lambda.
		\end{split}
	\end{equation}

	Now we put together \eqref{5.67} and \eqref{5.73}. By \eqref{3.5'} and \eqref{5.10'}, we have
	\begin{equation}
		v(\lambda)^2=\frac{1+w(\lambda)}{\lambda},
	\end{equation}
	and thus
	\begin{equation}
		\begin{split}
			w^2(w^2(1-w)v^2)^{k}\lambda^{k-1}=&w^2(w^2(1-w^2))^{k}/\lambda,\\
			\pa_\lambda (w^2(1-w)v^2)w^2(w^2(1-w)v^2)^{k-1}\lambda^k=&\pa_\lambda (w^2(1-w^2))w^2(w^2(1-w^2))^{k-1}\\
			&-w^2(w^2(1-w^2))^{k}/\lambda.
		\end{split}
	\end{equation}
	Then, up to the complicated coefficient, the sum of \eqref{5.67} and \eqref{5.73} is equal to
	\begin{equation}
		\begin{split}
			\frac{1}{\Gamma(k)}\int_{0}^\infty\pa_\lambda (w^2(1-w^2))w^2(w^2(1-w^2))^{k-1}d\lambda.
		\end{split}
	\end{equation}
	Since $w(0)=-1,w(\infty)=0$, make the change of variable $a=(1-w^2)$, we get
	\begin{equation}\label{5.77}
		\begin{split}
			&\frac{1}{\Gamma(k)}\int_0^1(1-a)^2(a(1-a))^{k-1}-(a(1-a))^{k}da\\
			&=\frac{1}{\Gamma(k)}\big(B(k,k+2)-B(k+1,k+1)\big)=\frac{k!}{(2k+1)!}.
		\end{split}
	\end{equation}
	
	By \eqref{5.66}, \eqref{5.77}, the first two terms in \eqref{5.54} contribute
	\begin{equation}\label{5.78}
		\begin{split}
			\frac{k!}{(2k+1)!}	\bigg[\frac{1}{(4\pi )^{m/2}}\det{}^{1/2}\bigg(\frac{R^{T\widetilde{M}}_{x_0}/2}{\sinh \big(R^{T\widetilde{M}}_{x_0}/2\big)}\bigg)(f_0df_1\wedge\cdots \wedge df_{2k})_{x_0}\bigg]^{\max}.
		\end{split}
	\end{equation}
	The computation of the third term in \eqref{5.54} is analogous to that of the second term and yields the same contribution as \eqref{5.78}. 
	
	Finally, combining \eqref{5.8}, \eqref{5.9}, \eqref{5.33}, \eqref{5.54} and \eqref{5.78}, we obtain \eqref{3.14} and \eqref{3.16}.

	\section{A rigorous treatment of the asymptotic expansion}\label{FAM}

	In this section, we provide detailed estimates to justify the asymptotic expansions \eqref{5.47}, \eqref{5.51} and \eqref{5.54}, following the analytic technique of Bismut-Lebeau \cite[\S\,11]{MR1188532} and Ma-Marinescu \cite[\S\,4]{MR2339952} with certain
	adaptations. In \cref{ART1}, we introduce a family of weighted Sobolev spaces tailored for the structure of the rescaled operator $\Delta_\varepsilon^{x_0}$ given in \eqref{5.36} and \eqref{5.41}. In \cref{ART2}, we give regularity estimates for the resolvent of $\Delta_\varepsilon^{x_0}$ following \cite[\S\,4]{MR4665497}. In \cref{ART3}, we define holomorphic functional calculus of $\Delta_\varepsilon^{x_0}$. In \cref{ART4}, we prove the convergence of functional calculus.

	\subsection{A family of Sobolev spaces with weights}\label{ART1}

	For any $0<\varepsilon\leqslant 1$ and $u\in\mathscr{C}^\infty_0\big(T_{x_0}\widetilde{M},\Lambda^\ell(T_{x_0}^*\widetilde{M})\big)$, the space of sections with compact support, we define Sobolev norms
	\begin{equation}\label{6.1}
		\begin{split}
			\Vert u\Vert_{0,\varepsilon}^2&=\int_{T_{x_0}\widetilde{M}}\Vert u(X)\Vert_{S_{\widetilde{M},x}}^2\lk1+(1+\vert X\vert^2)^{1/2}\rho(\varepsilon X)\rk^{n-\ell}\de X,\\
			\Vert u\Vert_{k,\varepsilon}^2&=\sum_{1\leqslant j\leqslant k}\sum_{1\leqslant i_1,\cdots,i_j\leqslant m}\Vert\pa_{X_{i_1}}\cdots\pa_{X_{i_j}}u\Vert_{0,\varepsilon}^2,
		\end{split}
	\end{equation}
	where $\rho(X)$ is definded in \eqref{5.18'}. Let $\mathscr{H}^{k,\varepsilon}$ be the Sobolev space with respect to $\lV\cdot\rV_{k,\varepsilon}$ and $\mathscr{H}^{-k,\varepsilon}$ the Sobolev space of negative order with norm defined by
	\begin{equation}
		\Vert u\Vert_{-k,\varepsilon}=\sup_{s'\neq0,s'\in \mathscr{H}_{k,\varepsilon}}\frac{\lv\langle u,u'\rangle_{0,\varepsilon}\rv}{\lV u'\rV_{k,\varepsilon}}.
	\end{equation}
	Clearly $\lV\cdot\rV_{k,\varepsilon}$ is increasing on $k$. Recall that $\rho(X)$ is radically decreasing, hence $\lV\cdot\rV_{k,\varepsilon}$ is decreasing on $\varepsilon$, and $\lV\cdot\rV_{k,\infty}$ is the usual Sobolev norm.

	These norms with weights are tailored for the following elliptic estimate.
	\begin{prop}
		There exists $C>0$  such that for any $0<\varepsilon\leqslant 1$ and $u,u'\in\mathscr{C}^\infty_0\big(T_{x_0}\widetilde{M},\Lambda (T_{x_0}\widetilde{M})\big)$, we have
		\begin{equation}\label{6.4}
			\begin{split}
				\Re\big\langle \Delta_\varepsilon^{x_0} u,u\big\rangle_{0,\varepsilon}&\geqslant C\lV u\rV_{1,\varepsilon}^2-C\lV u\rV_{0,\varepsilon}^2,\\
				\lv\Im\big\langle \Delta_\varepsilon^{x_0} u,u\big\rangle_{0,\varepsilon}\rv&\leqslant C\lV u\rV_{1,\varepsilon}\lV u\rV_{0,\varepsilon},\\
				\lv\big\langle \Delta_\varepsilon^{x_0} u,u'\big\rangle_{0,\varepsilon}\rv&\leqslant C\lV u\rV_{1,\varepsilon}\lV u'\rV_{1,\varepsilon}.
			\end{split}
		\end{equation}
	\end{prop}
	
	\begin{pro}
		Recall the expansion of $\Delta_\varepsilon^{x_0}$ given in \eqref{5.36}. All the seemingly unbounded terms are compensated for by weights. To see this, we easily check that the following terms
		\begin{equation}\label{6.5}
			\begin{split}
				&\lk\rho(\varepsilon X)\lv X\rv\rk^{1/2}\lk1+(1+\vert X\vert^2)^{1/2}\rho(\varepsilon X)\rk^{-1/2},\\ &\varepsilon\lk\rho(\varepsilon X)\lv X\rv\rk^{1/2}\lk1+(1+\vert X\vert^2)^{1/2}\rho(\varepsilon X)\rk^{1/2},\\
				&\pa_{X_i}\lk1+(1+\vert X\vert^2)^{1/2}\rho(\varepsilon X)\rk,\\
				&\pa_{X_i}\lk\varepsilon\rho(\varepsilon X)\langle \theta^{T\widetilde{M}}_{\varepsilon X}(\pa_{X_i})\widetilde{e}_k,\widetilde{e}_\ell\rangle c_\varepsilon(e^k)c_\varepsilon(e^\ell)\rk,
			\end{split}
		\end{equation}
		are uniformly bounded for $0<\varepsilon\leqslant 1,X\in T_{x_0}\widetilde{M}$. Then \eqref{6.4} is reduced to a standard elliptic estimate.\qed
	\end{pro}

	\begin{lemma}
		For any $k\in\mathbb{N}$, there exists $C>0$ such that for any $0<\varepsilon\leqslant 1$, $Q_1,\cdots Q_k\in\{\pa_{X_i},X_i\}_{i=1}^{m}$ and $u\in \mathscr{C}^\infty_0\big(T_{x_0}\widetilde{M},\Lambda(T_{x_0}^*\widetilde{M})\big)$, we have
		\begin{equation}\label{6.6}
			\lV\big[Q_{1},\big[Q_{2},\cdots\big[Q_{k},\Delta_\varepsilon^{x_0}\big]\cdots\big]\big]u\rV_{-1,\varepsilon}\leqslant C\lV u\rV_{1,\varepsilon}.
		\end{equation}
	\end{lemma}
	
	\begin{pro}
		It will be suﬀicient to prove that
		\begin{equation}\label{6.7}
			\lv\big\langle\big[Q_{1},\big[Q_{2},\cdots\big[Q_{k},\Delta_\varepsilon^{x_0}\big]\cdots\big]\big]u,u'\big\rangle_{0,\varepsilon}\rv\leqslant C\lV u\rV_{1,\varepsilon}\lV u'\rV_{1,\varepsilon}
		\end{equation}
		Since $[\pa_{X_i},X_j]=\delta_{ij}$, we know that $[X_j,\Delta_\varepsilon^{x_0}]$ verifies \eqref{6.7}. By the fourth inequality in \eqref{6.5}, $[\pa_{X_i},\Delta_\varepsilon^{x_0}]$ verifies \eqref{6.7}. Iterating this process, we see that $\big[Q_{1},\big[Q_{2},\cdots\big[Q_{k},\Delta_\varepsilon^{x_0}\big]\cdots\big]\big]$ has a similar structure as  $\Delta_\varepsilon^{x_0}$. By the third inequality in \eqref{6.4}, we get \eqref{6.7}.\qed
	\end{pro}

	\begin{lemma}
		For any $k\in\mathbb{N}$, there exists $C>0$ such that for any $0<\varepsilon\leqslant 1$, $Q_1,\cdots Q_k\in\{\pa_{X_i},X_i\}_{i=1}^{m}$ and $u\in \mathscr{C}^\infty_0\big(T_{x_0}\widetilde{M},\Lambda(T_{x_0}^*\widetilde{M})\big)$, we have
		\begin{equation}\label{6.8}
			\lV\big[Q_{1},\big[Q_{2},\cdots\big[Q_{k},\Delta_\varepsilon^{x_0}-\Delta_0^{x_0}\big]\cdots\big]\big]u\rV_{-1,\varepsilon}\leqslant C\varepsilon\sum_{\lv\alpha\rv\leqslant 2}\lV X^\alpha u\rV_{1,0}.
		\end{equation}
	\end{lemma}

	\begin{pro}
		By applying the Taylor expansion to \eqref{5.36} with respect to $\varepsilon$, we can analyze $(\frac{\pa\Delta_\varepsilon^{x_0}}{\pa\varepsilon})_{\varepsilon=\varepsilon_0}$. This operator has a similar structure as $\Delta_\varepsilon^{x_0}$, except for a possible extra term of order $O(\lv X\rv)$. Moreover, since $\varepsilon_0$ is arbitrary, we cannot use the trick in \eqref{6.5} of pairing $\rho(\varepsilon X)$ with $\lv X\rv$, which introduces another extra order $O(\lv X\rv)$. As a result, we have
		\begin{equation}
			\lv\big\langle(\Delta_\varepsilon^{x_0}-\Delta_0^{x_0})u,u'\big\rangle_{0,\varepsilon}\rv\leqslant C\varepsilon \sum_{\lv\alpha\rv\leqslant2}\lV X^\alpha u\rV_{1,\varepsilon}\cdot\lV u'\rV_{1,\varepsilon},
		\end{equation}
		and since $\lV\cdot\rV_{\ell,\varepsilon}$ is decreasing on $\varepsilon$, we obtain \eqref{6.8} for $k=0$. For $k\geqslant 1$, the argument leading to \eqref{6.6} also applies here. \qed
	\end{pro}

	\subsection{Sobolev inequalities for resolvents}\label{ART2}

	For any $c_1,c_2>0$, we define
	\begin{equation}\label{6.10}
		V_{c_1,c_2}=\{\lambda\in\mathbb{C}\mid \Re\lambda\leqslant c_1(\Im\lambda)^2-c_2\}.
	\end{equation}

	We have the following regularity property for the resolvent of $\Delta_\varepsilon^{x_0}$.
	\begin{prop}
		There exist $c_1,c_2>0$ such that for any $\lambda\in V_{c_1,c_2}$ the resolvent $(\lambda-\Delta_\varepsilon^{x_0})^{-1}$ exists. Moreover, there is $C>0$ such that for any for any $\lambda\in V_{c_1,c_2}$, $0<\varepsilon\leqslant 1$ and $u\in\mathscr{C}^\infty_0\big(T_{x_0}\widetilde{M},\Lambda(T_{x_0}^*\widetilde{M})\big)$, we have
		\begin{equation}\label{6.11}
			\begin{split}
				\lV(\lambda-\Delta_\varepsilon^{x_0})^{-1}u\rV_{0,\varepsilon}&\leqslant C\lV u\rV_{0,\varepsilon},\\
				\lV(\lambda-\Delta_\varepsilon^{x_0})^{-1}u\rV_{1,\varepsilon}&\leqslant C\lk1+\lv\lambda\rv\rk^2\lV u\rV_{-1,\varepsilon}.
			\end{split}
		\end{equation}
	\end{prop}
	
	\begin{pro}
		By \eqref{6.4}, we have
		\begin{equation}
			\begin{split}
				&\lv\left\langle\blk\lambda-\Delta_\varepsilon^{x_0}\brk u,u\right\rangle_{0,\varepsilon}\rv\\
				&\geqslant\sup\Big\{\Re\bli \Delta_\varepsilon^{x_0} u,u\bri_{0,\varepsilon}-\Re\lambda\lV u\rV_{0,\varepsilon}^2,\Bv\Im\bli\Delta_\varepsilon^{x_0} u,u\bri_{0,\varepsilon}-\Im\lambda\lV u\rV_{0,p,\mu}^2\Bv\Big\}\\
				&\geqslant \sup\Big\{C\lV u\rV_{1,\varepsilon}^2-(\Re\lambda+C)\lV u\rV_{0,\varepsilon}^2,\lv \Im\lambda\rv\lV u\rV_{0,\varepsilon}^2-C\lV u\rV_{0,\varepsilon}\lV u\rV_{1,\varepsilon}\Big\}.
			\end{split}
		\end{equation}
		Then by the obvious inequality $\lV u\rV_{0,\varepsilon}\leqslant \lV u\rV_{1,\varepsilon}$ we get
		\begin{equation}\label{6.13}
			\bv\big\langle\blk\lambda-\Delta_\varepsilon^{x_0}\brk u,u\big\rangle_{0,\varepsilon}\bv\geqslant C(\lambda) \lV u\rV_{0,\varepsilon}^2,
		\end{equation}
		where $C(\lambda)$ is given by
		\begin{equation}\label{6.14}
			C(\lambda)=\inf_{a\geqslant 1}\sup\big\{Ca^2-(\Re\lambda+C), \lv \Im\lambda\rv-Ca\big\}.
		\end{equation}
		For any $c>0$, by \eqref{6.14}, $C(\lambda)\geqslant c$ if and only if $\{a\geqslant1 \}$ is contained in
		\begin{equation}
			\{a\in\mathbb{R}\mid Ca^2-(\Re\lambda+C)\geqslant c\}\bigcup\{a\in\mathbb{R}\mid\lv \Im\lambda\rv-Ca\geqslant c\},
		\end{equation}
		which means that at least one of the following conditions holds
		\begin{equation}\label{6.16}
			\begin{cases}
				C-(\Re\lambda+C)\geqslant c,\\
				(C^{-1}({\Re\lambda+C+c}))^{1/2}\leqslant C^{-1}({\lv \Im\lambda\rv-c}).
			\end{cases}
		\end{equation}
		Simplify this, we see that $C(\lambda)\geqslant c$ if and only if at least one of the following inequalities holds 
		\begin{equation}
			\Re\lambda\leqslant\begin{cases}-c,\\
				C^{-1}\lk\lv \Im\lambda\rv-c\rk^2-C-c.
			\end{cases}
		\end{equation}
		Therefore, we can choose $0<c_1<C^{-1}$ and $c_2>C+c$ such that every point in $V_{c_1,c_2}$ satisfies one of these two inequalities, hence
		\begin{equation}\label{6.18}
			\inf_{\lambda\in V_{c_1,c_2}}C(\lambda)\geqslant  c>0.
		\end{equation}
		Then by \eqref{6.13}, if $\lambda\in V_{c_1,c_2}$, the 
		the resolvent  $\blk\lambda-\Delta_\varepsilon^{x_0}\brk^{-1}$ exists and 
		\begin{equation}\label{6.19}
			\lV(\lambda-\Delta_\varepsilon^{x_0})^{-1}u\rV_{0,\varepsilon}\leqslant c^{-1}\lV u\rV_{0,\varepsilon}.
		\end{equation}
		
		By \eqref{6.4}, it is clear that $\Delta_\varepsilon^{x_0}$ extends to a continuous linear map from $H_{1,\varepsilon}$ into $H_{-1,\varepsilon}$, and if we take $\lambda_0=-C$, then
		\begin{equation}\label{6.20}
			\bV\blk\lambda_0-\Delta_\varepsilon^{x_0}\brk^{-1}s\bV_{1,\varepsilon}\leqslant C^{-1}\lV u\rV_{-1,\varepsilon}.
		\end{equation}
		For $\lambda\in V_{c_1,c_2}$, we have
		\begin{equation}\label{6.21}
			(\lambda-\Delta_\varepsilon^{x_0})^{-1}=(\lambda_0-\Delta_\varepsilon^{x_0})^{-1}+\lk\lambda_0-\lambda\rk(\lambda-\Delta_\varepsilon^{x_0})^{-1}(\lambda_0-\Delta_\varepsilon^{x_0})^{-1}.
		\end{equation}
		By \eqref{6.19}, \eqref{6.20} and \eqref{6.21}, we get
		\begin{equation}\label{6.22}
			\begin{split}
				\bV(\lambda-\Delta_\varepsilon^{x_0})^{-1}u\bV_{0,\varepsilon}&\leqslant \lk C^{-1}+C^{-1}c^{-1}\lv\lambda-\lambda_0\rv\rk\lV u\rV_{-1,\varepsilon}.
			\end{split}
		\end{equation}
		On the other hand, we also have
		\begin{equation}
			(\lambda-\Delta_\varepsilon^{x_0})^{-1}=(\lambda_0-\Delta_\varepsilon^{x_0})^{-1}+\lk\lambda_0-\lambda\rk(\lambda_0-\Delta_\varepsilon^{x_0})^{-1}(\lambda-\Delta_\varepsilon^{x_0})^{-1},
		\end{equation}
		which, together with \eqref{6.20} and \eqref{6.22}, implies
		\begin{equation}
			\begin{split}
				\bV(\lambda-\Delta_\varepsilon^{x_0})^{-1}u\bV_{1,\varepsilon}&\leqslant C^{-1}+C^{-1}\lv\lambda_0-\lambda\rv\lk C^{-1}+C^{-1}c^{-1}\lv\lambda-\lambda_0\rv\rk\\
				&\leqslant C\lk1+\lv\lambda\rv\rk^2.
			\end{split}
		\end{equation}
		We complete the proof of \eqref{6.11}.\qed
	\end{pro}
	
	From now on, we fix the constants $c_1,c_2$ used in \eqref{6.18}. We now establish the higher order regularity for the resolvent of $\Delta_\varepsilon^{x_0}$.
	\begin{prop}
		For any $\lambda\in V_{c_1,c_2}$, $0<\varepsilon\leqslant 1$ and $k\in\mathbb{N}$, the resolvent $(\lambda-\Delta_\varepsilon^{x_0})^{-1}$ maps $\mathscr{H}_{k,\varepsilon}$ to $\mathscr{H}_{k+1,\varepsilon}$. Indeed, for any $\alpha\in\mathbb{N}^m$, there exists $C>0$ such that for any $u\in\mathscr{C}^\infty_0\big(T_{x_0}\widetilde{M},\Lambda(T_{x_0}^*\widetilde{M})\big)$, we have
		\begin{equation}\label{6.25}
			\lV X^\alpha(\lambda-\Delta_\varepsilon^{x_0})^{-1}u\rV_{k+1,\varepsilon}\leqslant C\lk1+\lv\lambda\rv\rk^{2(k+\lv\alpha\rv+1)}\sum_{\beta\leqslant \alpha}\lV X^\beta u\rV_{k,\varepsilon}.
		\end{equation}
	\end{prop}

	\begin{pro} 
		
		First, for any $Q_1,\cdots, Q_k\in\{\pa_{X_i}\}_{i=1}^{m}$, $Q_{k+1},\cdots, Q_{k+\lv\alpha\rv}\in \{X_i\}_{i=1}^{m}$, we can express $Q_1\cdots Q_{k+{\vert}\alpha{\vert}}(\lambda-\Delta_\varepsilon^{x_0})^{-1}$ as a linear combination of operators of the type 
		\begin{equation}\label{6.26}
			\big[Q_{i_1},\big[Q_{i_2},\cdots\big[Q_{i_{\ell}},(\lambda-\Delta_\varepsilon^{x_0})^{-1}\big]\cdots\big]\big]Q_{i_{\ell+1}}\cdots Q_{i_{k+\lv\alpha\rv}},
		\end{equation}
		where $1\leqslant i_1<\cdots<i_\ell\leqslant k+\lv\alpha\rv, 0\leqslant i_{\ell+1}<\cdots< i_{k+\lv\alpha\rv}\leqslant k+\lv\alpha\rv$ and $\{i_1,\cdots,i_{k+\lv\alpha\rv}\}=\{1,\cdots,k+\lv\alpha\rv\}$.
		
		Secondly, any commutator
		\begin{equation}
			\big[Q_{i_1},\big[Q_{i_2},\cdots\big[Q_{i_{\ell}},(\lambda-\Delta_\varepsilon^{x_0})^{-1}\big]\cdots\big]\big]
		\end{equation}
		is a linear combination of operators of the form
		\begin{equation}\label{6.28}
			(\lambda-\Delta_\varepsilon^{x_0})^{-1}C_{1,\varepsilon}(\lambda-\Delta_\varepsilon^{x_0})^{-1}C_{2,\varepsilon}\cdots C_{\ell',\varepsilon}(\lambda-\Delta_\varepsilon^{x_0})^{-1},
		\end{equation}
		where $C_{j,\varepsilon}$ takes the form $\big[Q_{j_1},\cdots\big[Q_{j_{\ell''}},\Delta_\varepsilon^{x_0}\big]\cdots\big]$ with $\ell''\leqslant\ell$.
		
		Combining \eqref{6.5} and \eqref{6.11}, we obtain
		\begin{equation}
			\begin{split}
				\lV(\lambda-\Delta_\varepsilon^{x_0})^{-1}C_{1,\varepsilon}(\lambda-\Delta_\varepsilon^{x_0})^{-1}C_{2,\varepsilon}\cdots C_{\ell',\varepsilon}(\lambda-\Delta_\varepsilon^{x_0})^{-1}Q_{i_{\ell+1}}\cdots Q_{i_{k+\lv\alpha\rv}}u\rV_{1,\varepsilon}&\\
				\leqslant (1+\lv\lambda\rv)^{2(k+\lv\alpha\rv+1)}\sum_{\beta\leqslant\alpha}\lV X^\beta u\rV_{k,\varepsilon}&,
			\end{split}
		\end{equation}
		which, together with \eqref{6.26}, yields \eqref{6.25}.\qed
	\end{pro}

	\begin{prop}
		For $k\in\mathbb{N}$ and $\alpha\in \mathbb{N}^m$, there exists $C>0$ and $\ell\in\mathbb{N}$ such that for any $0<\varepsilon\leqslant 1$ and $u\in\mathscr{C}^\infty_0\big(T_{x_0}\widetilde{M},\Lambda(T_{x_0}^*\widetilde{M})\big)$, we have
		\begin{equation}\label{6.30}
			\lV\big((\lambda-\Delta_\varepsilon^{x_0})^{-k-\lv\alpha\rv}-(\lambda-\Delta_0^{x_0})^{-k-\lv\alpha\rv}\big)\pa^{\alpha}u\rV_{k,\varepsilon}\leqslant C\varepsilon\big(1+\lv\lambda\rv\big)^\ell\sum_{\lv\beta\rv\leqslant \ell}\lV X^\beta u\rV_{0,0}.
		\end{equation}
	\end{prop}

	\begin{pro}
		First, let us consider the case $\alpha=0$. We have
		\begin{equation}
			(\lambda-\Delta_\varepsilon^{x_0})^{-1}-(\lambda-\Delta_0^{x_0})^{-1}=-(\lambda-\Delta_\varepsilon^{x_0})^{-1}(\Delta_\varepsilon^{x_0}-\Delta_0^{x_0})(\lambda-\Delta_0^{x_0})^{-1},
		\end{equation}
		and in general,
		\begin{equation}
			\begin{split}
				&(\lambda-\Delta_\varepsilon^{x_0})^{-k}-(\lambda-\Delta_0^{x_0})^{-k}\\
				&=\sum_{i=0}^{k-1}(\lambda-\Delta_\varepsilon^{x_0})^{-(k-i-1)}\big((\lambda-\Delta_\varepsilon^{x_0})^{-1}-(\lambda-\Delta_0^{x_0})^{-1}\big)(\lambda-\Delta_0^{x_0})^{-i}\\
				&=-\sum_{i=0}^{k-1}(\lambda-\Delta_\varepsilon^{x_0})^{-(k-i)}\big(\Delta_\varepsilon^{x_0}-\Delta_0^{x_0}\big)(\lambda-\Delta_0^{x_0})^{-i-1}.
			\end{split}
		\end{equation}
		Then, we play the commutator game again for $Q_1\cdots Q_{k}(\lambda-\Delta_\varepsilon^{x_0})^{-(k-i)}\big(\Delta_\varepsilon^{x_0}-\Delta_0^{x_0}\big)(\lambda-\Delta_0^{x_0})^{-i-1}$ as in \eqref{6.26}. But here, we leave a $\pa_{X_{j}}$ in front of each operator of the form
		\begin{equation}
			\big[Q_{i_1},\cdots,\big[Q_{i_{\ell}},(\lambda-\Delta_\varepsilon^{x_0})^{-1}\big]\cdots\big],\ \ \big[Q_{i_1},\cdots,\big[Q_{i_{\ell}},(\lambda-\Delta_0^{x_0})^{-1}\big]\cdots\big]
		\end{equation}
		except the one after $\big(\Delta_\varepsilon^{x_0}-\Delta_0^{x_0}\big)$. Since $(k-i)+(i+1)-1=k$, after this procedure, we get a product of operators with the following form
		\begin{equation}\label{6.34}
			\begin{split}
				&\pa_{X_{j}}\big[Q_{i_1},\cdots,\big[Q_{i_{\ell}},(\lambda-\Delta_\varepsilon^{x_0})^{-1}\big]\cdots\big],\\
				&\big[Q_{i_1},\cdots,\big[Q_{i_{\ell}},\Delta_\varepsilon^{x_0}-\Delta_0^{x_0}\big]\cdots\big]\cdot \big[Q_{i_1},\cdots,\big[Q_{i_{\ell}},(\lambda-\Delta_0^{x_0})^{-1}\big]\cdots\big],\\
				&\pa_{X_{j}}\big[Q_{i_1},\cdots,\big[Q_{i_{\ell}},(\lambda-\Delta_0^{x_0})^{-1}\big]\cdots\big].
			\end{split}
		\end{equation}
		Combining this with \eqref{6.6}, \eqref{6.8}, \eqref{6.11}, \eqref{6.25} and \eqref{6.28}, we get \eqref{6.30}. 
		
		For the general case $\alpha\neq0$, the argument proceeds similarly, except that the commutator trick is also performed from right to left. \qed
	\end{pro}

	\subsection{Holomorphic functional calculus}\label{ART3}

	Now we state Paley-Wiener theorem, which follows from Fourier inversion theorem and integration by parts, see \cite[Theorem 7.22]{MR1157815} for instance.

	\begin{lemma}
		Let $\psi\in\mathscr{S}(\mathbb{R})$ with $\mathrm{supp}(\widehat{\psi})\subseteq [-r,r]$, then $\psi$ admits an analytic continuation to an entire function. Moreover, for any $k\in\mathbb{N}$, there is $C>0$ such that for any $\lambda\in\mathbb{C}$,
		\begin{equation}\label{6.35}
			\lv \psi(\lambda)\rv\leqslant C(1+\lv \lambda\rv)^{-k}e^{r\lv\Im \lambda\rv}.
		\end{equation}
	\end{lemma}

	Recall $f'(\lambda),g(\lambda)$ used in \eqref{3.5'}, clearly \eqref{6.35} applies to them. Therefore, $f(\lambda)$ also has a continuation to an entire function by
	\begin{equation}\label{6.36}
		f(\lambda)=-\int_{\lambda}^{\infty}f'(\mu)d\mu,
	\end{equation}
	where the integration is over the horizontal ray from $\lambda=\Re\lambda+i\Im\lambda$ to $+\infty+i\Im\lambda$. Moreover, using \eqref{6.35}, \eqref{6.36} and the fact that $\lv\lambda\rv\leqslant 1$ is compact, we see that for any $k\in\mathbb{N}$, there are $C,r>0$ such that for any $\Re \lambda\geqslant 0$,
	\begin{equation}\label{6.38}
		\lv\frac{f(\lambda)}{\lambda}\rv\leqslant C(1+\lv \lambda\rv)^{-k}e^{r\lv\Im\lambda\rv}.
	\end{equation}

	Let us now use $\phi(\lambda)$ to denote any function appearing in \eqref{5.11}, \eqref{5.24}, \eqref{5.38} for use in functional calculus. By \eqref{5.10'}, \eqref{5.10}, \eqref{5.23} and \eqref{5.37}, $\phi(\lambda)$ is the product of several $v(\lambda)$ and $w(\lambda)$. From \eqref{6.35}, \eqref{6.38} and the following identity 
	\begin{equation}
		\lambda^{1/2}=\pm\lk\lk\frac{\lv\lambda\rv+\Re\lambda}{2}\rk^{1/2}+\sqrt{-1}\frac{\Im\lambda}{\lv\Im\lambda\rv}\lk\frac{\lv\lambda\rv-\Re\lambda}{2}\rk^{1/2}\rk,
	\end{equation}
	we obtain that for any $k\in\mathbb{N}$, there are $C,r>0$ such that for any $\Re \lambda\geqslant 0$,
	\begin{equation}\label{6.39}
		\lv \phi(\lambda)\rv\leqslant C(1+\lv \lambda\rv)^{-k}e^{r(\lv\lambda\rv-\Re\lambda)}.
	\end{equation}
	For any $n\in\mathbb{N}$, similar to \eqref{6.36}, we can define a $n$-th holomorphic antiderivative of $\phi$ inductively by
	\begin{equation}
		\phi_n^{}(\lambda)=-\int_{\lambda}^{\infty}\phi_{n-1}(\mu)d\mu,
	\end{equation}
	and $\phi_n(\lambda)$ also verifies the property \eqref{6.36}.

	\begin{lemma}
		For any $n,k\in\mathbb{N}$, there is $C>0$ such that for any $\lambda\in\pa V_{c_1,c_2}$, the boundary of $V_{c_1,c_2}$ defined in \eqref{6.10} and \eqref{6.18}, we have
		\begin{equation}\label{6.40'}
			\lv\phi_n(\lambda)\rv\leqslant C(1+\lv\lambda\rv)^{-k}.
		\end{equation}
		Also, the Cauchy integral formula applies to $\phi_n$ and the infinite contour $\pa V_{a,b}$, that is, for any $\mu\notin V_{c_1,c_2}$, we have
		\begin{equation}\label{6.42}
			\phi(\mu)=\frac{n!}{2\pi \sqrt{-1}}\int_{\pa V_{c_1,c_2}}\phi_n(\lambda)(\lambda-\mu)^{-n-1}d\lambda.
		\end{equation}	
	\end{lemma}
	
	\begin{pro}
		By \eqref{6.39}, we only need to check that $\sup_{\lambda\in \pa V_{c_1,c_2}}(\lv\lambda\rv-\Re\lambda)<\infty$. For $\lambda\in \pa V_{c_1,c_2}$, we have
		\begin{equation}\label{6.43}
			\begin{split}
				\lv\lambda\rv-\Re\lambda&=\blk(c_1(\Im\lambda)^2-c_2)^2+(\Im\lambda)^2\brk^{1/2}-\big(c_1(\Im\lambda)^2-c_2\big)\\
				&=(\Im\lambda)^2\Big/\Big(\blk(c_1(\Im\lambda)^2-c_2)^2+(\Im\lambda)^2\brk^{1/2}+\big(c_1(\Im\lambda)^2-c_2\big)\Big).
			\end{split}
		\end{equation}
		It is a fortunate that in \eqref{6.43}, the highest order of $\Im\lambda$ in both the numerator and the denominator is $2$!\qed
	\end{pro}

	By \eqref{6.11} and \eqref{6.42}, we can justify the holomorphic functional calculus of $\Delta_\varepsilon^{x_0}$.
	\begin{coro}
		The holomorphic functional calculus of $\Delta_\varepsilon^{x_0}$ is legal, and for any $n\in\mathbb{N}$,
		\begin{equation}\label{6.44}
			\phi(\Delta_\varepsilon^{x_0})=\frac{n!}{2\pi \sqrt{-1}}\int_{\pa V_{c_1,c_2}}\phi_{n}(\lambda)(\lambda-\Delta_\varepsilon^{x_0})^{-n-1}d\lambda.
		\end{equation}
	\end{coro}

	\subsection{Convergence of functional calculus}\label{ART4}

	By \eqref{5.47}, the limit in \eqref{5.54} takes the form of 
	\begin{equation}
		\lim_{\varepsilon\to0} \big(L_\varepsilon^{x_0}\phi(\Delta_\varepsilon^{x_0})\big)(0,0)=\big(L_0^{x_0}\phi(\Delta_0^{x_0})\big)(0,0),
	\end{equation}
	where $L_\varepsilon^{x_0}$ is a family of differential operators, and the coefficients of $L_\varepsilon^{x_0}$ converge to those of $L_0^{x_0}$ as $\varepsilon\to 0$ uniformly over any compact subset of $T_{x_0}\widetilde{M}$. Therefore, 	it is enough to establish the following stronger limit
	\begin{equation}\label{6.45}
		\lim_{\varepsilon\to 0}\phi(\Delta_\varepsilon^{x_0})(X,X')=\phi(\Delta_0^{x_0})(X,X')
	\end{equation}
	in the $C^\infty$-sense over any compact subset of $T_{x_0}\widetilde{M}\times T_{x_0}\widetilde{M}$.

	\begin{prop}
		There exist $C>0,k\in\mathbb{N}$ such that
		\begin{equation}\label{6.46}
			\lV\pa^\alpha_X\pa^\beta_{X'}\psi(\Delta_\varepsilon^{x_0})(X,X')\rV_{\mathrm{End}(\Lambda(T_{x_0}^*\widetilde{M}))}\leqslant C\lk1+\lv X\rv+\lv X'\rv\rk^k.
		\end{equation}
	\end{prop}
	
	\begin{pro}
		By \eqref{6.25}, we have
		\begin{equation}\label{6.47}
			\lV(\Delta^{T_{x_0}\widetilde{M}})^k(\lambda-\Delta_\varepsilon^{x_0})^{-2k}u\rV_{0,\varepsilon}\leqslant C\lk1+\lv\lambda\rv\rk^\ell\lV u\rV_{0,\varepsilon},
		\end{equation}
		where $\Delta^{T_{x_0}\widetilde{M}}$ is the standard Laplacian.
		
		Let $\Delta_\varepsilon^{x_0,*}$ be the formal adjoint of $\Delta_\varepsilon^{x_0}$ with respect to the usual $L^2$-inner product $\langle\cdot,\cdot\rangle_{0,\infty}$, then $\Delta_\varepsilon^{x_0,*}$ has the same structure as $\Delta_\varepsilon^{x_0}$ except that we interchange the operators $e^i\wedge$ and $\iota_{e_i}$. For any $0<\varepsilon\leqslant 1$, we define Sobolev norms $\lV\cdot\rV_{k,\varepsilon,*}$ by replacing $(n-q)$ with $(q-n)$ in \eqref{6.1}. From above analysis, we get the same inequalities for $\Delta_\varepsilon^{x_0,*}$ with respect to $\lV\cdot\rV_{k,\varepsilon,*}$. Also, the weights of $\lV\cdot\rV_{k,\varepsilon}$ and $\lV\cdot\rV_{k,\varepsilon,*}$ cancel out to get $\lV\cdot\rV_{k,\infty}$, hence
		\begin{equation}\label{6.48}
			\begin{split}
				\lV(\lambda-\Delta_\varepsilon^{x_0})^{-2k}(\Delta^{T_{x_0}\widetilde{M}})^ku\rV_{0,\varepsilon}&=\sup_{\lV s'\rV_{0,\varepsilon,*}\leqslant1} \lv\big\langle s',(\lambda-\Delta_\varepsilon^{x_0})^{-2k}(\Delta^{T_{x_0}\widetilde{M}})^ks\big\rangle_{0,\infty}\rv\\
				&=\sup_{\lV s'\rV_{0,\varepsilon,*}\leqslant1} \lv\big\langle (\Delta^{T_{x_0}\widetilde{M}})^k(\lambda-\Delta_\varepsilon^{x_0,*})^{-2k} s',s\big\rangle_{0,\infty}\rv\\
				&\leqslant \sup_{\lV s'\rV_{0,\varepsilon,*}\leqslant1}\lV(\Delta^{T_{x_0}\widetilde{M}})^k(\lambda-\Delta_\varepsilon^{x_0,*})^{-2k} s'\rV_{0,\varepsilon,*}\lV u\rV_{0,\varepsilon}\\
				&\leqslant  C\lk1+\lv\lambda\rv\rk^\ell\lV u\rV_{0,\varepsilon}.
			\end{split}
		\end{equation}
		Note that \eqref{6.48} can also be derived using the commutator trick applied from right to left, as in the proof of \eqref{6.30}.
		
		Applying \eqref{6.44} for $n=4k-1$, we get
		\begin{equation}
			\psi(\Delta_\varepsilon^{x_0})=\frac{(4k)!}{2\pi \sqrt{-1}}\int_{\pa V_{c_1,c_2}}\psi_{4k-1}(\lambda)(\lambda-\Delta_\varepsilon^{x_0})^{-4k}d\lambda,
		\end{equation}
		then by \eqref{6.40'}, \eqref{6.47} and \eqref{6.48}, we get
		\begin{equation}
			\lV(\Delta^{T_{x_0}\widetilde{M}})^k\psi(\Delta_\varepsilon^{x_0})(\Delta^{T_{x_0}\widetilde{M}})^ku\rV_{0,\varepsilon}\leqslant C\lV u\rV_{0,\varepsilon}.
		\end{equation}
		Since $(\Delta^{T_{x_0}\widetilde{M}})^k\psi(\Delta_\varepsilon^{x_0})(\Delta^{T_{x_0}\widetilde{M}})^k$ has kernel $(\Delta^{T_{x_0}\widetilde{M}}_{X})^k(\Delta^{T_{x_0}\widetilde{M}}_{X'})^k\psi(\Delta_\varepsilon^{x_0})(X,X')$, using Sobolev embedding theorem and the fact that $\lV u\rV_{0,\varepsilon}\leqslant \Vert(1+X)^\frac{m}{2}s(X)\Vert_{0,\infty}$, we get \eqref{6.46}.\qed
	\end{pro}

	\begin{prop}
		For any $k\in\mathbb{N}$, we have
		\begin{equation}\label{6.51}
			\lV(\Delta^{T_{x_0}\widetilde{M}})^k	\lk\psi(\Delta_\varepsilon^{x_0})-\psi(\Delta_0^{x_0})\rk (\Delta^{T_{x_0}\widetilde{M}})^ku\rV_{0,\varepsilon}\leqslant C\varepsilon\sum_{\lv\beta\rv\leqslant \ell}\lV X^\alpha u\rV_{0,0}.
		\end{equation}
	\end{prop}

	\begin{pro}
		By \eqref{6.44}, we have
		\begin{equation}\label{6.52}
			\psi(\Delta_\varepsilon^{x_0})-\psi(\Delta_0^{x_0})=\frac{(4k)!}{2\pi \sqrt{-1}}\int_{\pa V_{c_1,c_2}}\psi_{4k-1}(\lambda)\big((\lambda-\Delta_\varepsilon^{x_0})^{-4k}-(\lambda-\Delta_0^{x_0})^{-4k}\big)d\lambda.
		\end{equation}
		By \eqref{6.30}, we get
		\begin{equation}\label{6.53}
			\begin{split}
				\lV\big(\Delta^{T_{x_0}\widetilde{M}})^k((\lambda-\Delta_\varepsilon^{x_0})^{-2k}-(\lambda-\Delta_0^{x_0})^{-2k}\big)(\Delta^{T_{x_0}\widetilde{M}})^ku\rV_{0,\varepsilon}&\\
				\leqslant C\varepsilon\big(1+\lv\lambda\rv\big)^\ell\sum_{\lv\beta\rv\leqslant \ell}\lV X^\alpha u\rV_{0,0}&.
			\end{split}
		\end{equation}
		Using \eqref{6.40'}, \eqref{6.52} and \eqref{6.53}, we get \eqref{6.51}.\qed
	\end{pro}

	Now we are finally ready to prove \eqref{6.45}.
	
	\begin{theo}
		There exist $C>0$ and $k\in \mathbb{N}$ such that
		\begin{equation}\label{6.54}
			\lV\pa^\alpha_X\pa^\beta_{X'}\big(\psi(\Delta_\varepsilon^{x_0})(X,X')-\psi(\Delta_0^{x_0})(X,X')\big)\rV_{\mathrm{End}(\Lambda(T_{x_0}^*\widetilde{M}))}\leqslant C\varepsilon^{\tfrac{1}{m+1}}\lk1+\lv X\rv+\lv X'\rv\rk^k.
		\end{equation}
	\end{theo}

	\begin{pro}
		
		Let us denote
		\begin{equation}
			K_\varepsilon(X,X')=(\Delta^{T_{x_0}\widetilde{M}}_{X})^k(\Delta^{T_{x_0}\widetilde{M}}_{X'})^k\big(\psi(\Delta_\varepsilon^{x_0})(X,X')-\psi(\Delta_0^{x_0})(X,X')\big).
		\end{equation}
		Let $f(X)$ be a smooth function on $T_{x_0}\widetilde{M}$ with compact support and $\int_{T_{x_0}\widetilde{M}}\phi(X)\de X=1$. 
		
		For any $V,V'\in \Lambda(T_{x_0}^*\widetilde{M})$, by \eqref{6.46}, we have
		\begin{equation}
			\begin{split}
				\bbv\langle &K_\varepsilon(X,X')V,V'\rangle_{\Lambda(T_{x_0}^*\widetilde{M})}\\
				&-\int_{T_{x_0}\widetilde{M}\times T_{x_0}\widetilde{M}}\langle K_\varepsilon(X-\delta Y,X'-\delta Y)V,V'\rangle_{\Lambda(T_{x_0}^*\widetilde{M})}f(Y)f(Y')dYdY'\bbv\\
				&\leqslant C\delta\lk1+\lv X\rv+\lv X'\rv\rk^\ell \lV V\rV_{\Lambda(T_{x_0}^*\widetilde{M})}\lV V'\rV_{\Lambda(T_{x_0}^*\widetilde{M})}
			\end{split}
		\end{equation}
		On the other hand, by \eqref{6.51}, we have
		\begin{equation}
			\begin{split}
				&\lv\int_{T_{x_0}\widetilde{M}\times T_{x_0}\widetilde{M}}\langle K_\varepsilon(X-\delta Y,X'-\delta Y)v,v'\rangle_{S_{\widetilde{M},x}}f(Y)f(Y')dYdY'\rv\\
				&=\delta^{-2m}\lv\int_{T_{x_0}\widetilde{M}\times T_{x_0}\widetilde{M}}\langle K_\varepsilon(Y,Y')f((X-Y)/\delta)v,f((X'-Y')/\delta)v'\rangle_{\Lambda(T_{x_0}^*\widetilde{M})}dYdY'\rv\\
				&\leqslant\varepsilon \delta^{-m}\lk1+\lv X\rv+\lv X'\rv\rk^\ell \lV V\rV_{\Lambda(T_{x_0}^*\widetilde{M})}\lV V'\rV_{\Lambda(T_{x_0}^*\widetilde{M})}
			\end{split}
		\end{equation}
		since $\lV Y^\alpha f((X-Y)/\delta)\rV_{0,\infty}^2\leqslant C\delta^{m}\lk1+\lv X\rv\rk^{\lv\alpha\rv}$. By taking $\delta=\varepsilon^{\tfrac{1}{m+1}}$, we get \eqref{6.54}.\qed
	\end{pro}

	\section{Estimates on eigenvalues for manifolds with boundary}\label{SEE}
	
	In this section, we prove the kernel estimate \eqref{4.7''}. In \cref{A}, we recall two fundamental inequalities by Kato and Poincaré-Wirtinger. In \cref{SEB}, we apply Moser iteration \cite{MR159138} to estimate the eigensections of the Dirac operator. In \cref{SEE3}, we establish kernel estimates for certain functional calculus of the Dirac operator.

	\subsection{Kato and Poincaré-Wirtinger inequalities}\label{A}
	
	In this subsection, we review two important inequalities that will be used later.

	We retain the notation that $\Omega$ is a compact spin manifold with a nonempty boundary $\pa\Omega$ and $\dim\Omega=m$. We equip $\Omega$ with a Riemannian metric $g^{T\Omega}$, then let $dv_\Omega$ be the volume form on $\Omega$, $\nabla^{T\Omega}$ the Levi-Civita connection on $T\Omega$, $R^{T\Omega}$ the curvature and $\mathrm{Sc}_{g^{T\Omega}}$ the scalar curvature. Let $D^{S_\Omega}$ be the Dirac operator and $\Delta^{S_\Omega}$ the nonnegative Laplacian, both acting on the space $C^\infty_0(\Omega,S_\Omega)$ of smooth sections of $S_\Omega$ that vanish on $\pa \Omega$.

	Let $H^1_0(\Omega)$ be the Sobolev space on $\Omega$ and $H^{-1}_0(\Omega)$ its dual through the $L^2$-inner product $\langle\cdot,\cdot\rangle_{L^2(\Omega,S_\Omega)}$. Let $\Delta^\Omega$ be the nonnegative Laplacian acting on $C^\infty_0(\Omega)$, which extends to a map $H^1_0(\Omega)\to H^{-1}_0(\Omega)$. 
	
	Now we present the Kato's inequality for vector bundles by Hess-Schrader-Uhlenbrock \cite[Proposition 2.2]{MR602436}, which we apply here specifically to $(\Omega,S_\Omega)$.
	\begin{prop}
		For any $u\in C^\infty_0(\Omega,S_\Omega)$, we have the following inequality in $H^{-1}_0(\Omega)$
		\begin{equation}\label{2.3}
			\Delta^\Omega\Vert u\Vert_{S_\Omega}\leqslant\mathrm{Re}\langle \Delta^{S_\Omega}u,\mathrm{sign}_\xi u\rangle_{S_\Omega},
		\end{equation}
		where $\xi$ is an arbitrary measurable section of the unit sphere bundle of $S_\Omega$ and
		\begin{equation}
			\mathrm{sign}_\xi u=\begin{cases}u/\Vert u\Vert_{S_\Omega}\ \ &\text{on}\ \mathrm{supp}(u),\\
				\xi\ \ &\text{otherwise}.
			\end{cases}
		\end{equation}
	\end{prop}

	\subsubsection{The Dirichlet isoperimetric constant}
	
	The Dirichlet isoperimetric constant of $\Omega$ is defined by
	\begin{equation}\label{2.5''}
		\mathrm{ID}(\Omega)=\inf_{\begin{subarray}{c}
				\pa\Sigma\subset \Omega,\\ \partial \Sigma\cap \partial\Omega=\emptyset
		\end{subarray}}\frac{\mathrm{Area}(\partial\Sigma)}{\mathrm{Vol}(\Sigma)^{(m-1)/m}},
	\end{equation}
	where the infimum is taken over all subdomains $\Sigma\subset\Omega$ with the property that $\partial\Sigma$ is a hypersurface not intersecting $\partial\Omega$.
	
	The Dirichlet isoperimetric constant here corresponds to the one obtained by taking $\alpha=m/(m-1)$ in \cite[Definition 9.1]{MR2962229}. This Dirichlet isoperimetric constant here remains unchanged when the Riemannian metric is rescaled, thanks to the specific choice of the power in the denominator of its definition. Additionally, we note that the Dirichlet isoperimetric constant is continuous with respect to $C^0$-convergence of Riemannian metrics.
	Finally, we mention that intuitively, the Dirichlet isoperimetric constant of a Riemannian manifold $\Omega$ is small if there is a thin neck.

	We have the following Poincaré-Wirtinger inequality, see \cite[Corollary 9.8]{MR2962229}.
	\begin{prop}\label{p2.3'}
		There is $C_m>0$ such that for any $f\in H^1_0(\Omega)$, we have
		\begin{equation}\label{2.6'}
			\begin{split}
				\int_\Omega\vert\nabla f\vert^2dv_\Omega\geqslant C_m\mathrm{ID}(\Omega)^2\Big(\int_X\vert f\vert^{2m/(m-2)}dv_\Omega\Big)^{(m-2)/m}.
			\end{split}
		\end{equation}
	\end{prop}

	From Proposition \ref{p2.3'}, we have the following useful corollary.
	\begin{coro}
		For any $a\geqslant 2$ and nonnegative $f\in H^1_0(\Omega)$, we have
		\begin{equation}\label{2.7'}
			\int_\Omega f^{a-1}\Delta^\Omega fdv_\Omega \geqslant \frac{C_m\mathrm{ID}(\Omega)^2}{a}\Big(\int_\Omega f^{am/(m-2)}dv_\Omega\Big)^{(m-2)/m}.
		\end{equation}    
	\end{coro}
	
	\begin{pro}
		Integrating by parts implies
		\begin{equation}
			\begin{split}
				\int_\Omega f^{a-1}\Delta^\Omega fdv_\Omega&=(a-1)\int_\Omega f^{a-2}\big\vert\nabla f\big\vert^2dv_\Omega\\
				&=\frac{4(a-1)}{a^2}\int_\Omega\big\vert\nabla\big(f^{a/2}\big)\big\vert^2dv_\Omega,
			\end{split}
		\end{equation}    
		from which we deduce \eqref{2.7'} by replacing $f$ with $f^{a/2}$ in \eqref{2.6'}. \qed
	\end{pro}

	\subsection{Moser iteration}\label{SEB}

	Let us list the Dirichlet eigenvalues $\{\lambda_i\}_{i=1}^\infty$ of $D^{S_\Omega,2}$ counted with multiplicities and the corresponding normalized eigensections $\{u_i\}_{i=1}^\infty$, that is,
	\begin{equation}
		D^{S_\Omega,2}u_i=\lambda_iu_i, \ \ \ \lV u_i\rV_{L^2}=1,\ \ \ u_i\in H_0^1(\Omega,S_\Omega).
	\end{equation}
	
	\begin{prop}
		For any $k\in\mathbb{N}^*$ and $u\in\mathrm{span}\{u_1,\cdots,u_k\}$, we have
		\begin{equation}\label{3.2'}
			\lV u\rV_{L^\infty}\leqslant\frac{C\big(\lambda_k+\vert k_0\vert\big)^{m/4}}{\mathrm{ID}(\Omega)^{m/2}}\lV u\rV_{L^2}.
		\end{equation}
	\end{prop}
	
	\begin{pro}
		Taking $b=m/(m-2),\ell\in\mathbb{N}$, we claim that 
		\begin{equation}\label{2.10}
			\begin{split}
				\lV u\rV_{L^{2b^{(\ell+1)}}}\leqslant	\prod_{j=0}^{\ell}\Bigg(\frac{2b^j}{C_m\mathrm{ID}(\Omega)^2}\bigg(\lambda_k+\frac{\lv k_0\rv}{4}\bigg)\Bigg)^{1/(2b^j)}\lV u\rV_{L^{2b^{\ell}}}.
			\end{split}
		\end{equation}    
		To see this, let us choose $u$ such that
		\begin{equation}\label{2.11}
			\frac{\lV u\rV_{L^{2b^{(\ell+1)}}}}{\lV u\rV_{L^{2b^{\ell}}}}=\sup_{u\in\mathrm{span}\{u_1,\cdots,u_k\}}\frac{\lV u\rV_{L^{2b^{(\ell+1)}}}}{\lV u\rV_{L^{2b^{\ell}}}}.
		\end{equation}
		Let us write $\mathrm{Sc}_{g^{T\Omega}}=\mathrm{Sc}^+_{g^{T\Omega}}-\mathrm{Sc}^-_{g^{T\Omega}}$. From \eqref{2.3} and the obvious inequality
		\begin{equation}
			\int_\Omega\lV u\rV_{S_\Omega}^{2b^\ell}\mathrm{Sc}^+_{g^{T\Omega}}dv_\Omega\geqslant0,
		\end{equation}
		we get
		\begin{equation}\label{2.14'}
			\begin{split}
				\int_\Omega\lV u\rV_{S_\Omega}^{2b^\ell-1}\Delta^\Omega\lV u\rV_{S_\Omega} dv_\Omega&\leqslant \int_\Omega\lV u\rV_{S_\Omega}^{2b^\ell-2}\mathrm{Re}\langle \Delta^{S_\Omega}u,u\rangle_{S_\Omega} dv_\Omega\\
				&=\sum_{1\leqslant j\leqslant k} \int_\Omega\lV u\rV_{S_\Omega}^{2b^\ell-2}\mathrm{Re}\Big\langle\big(\lambda_j-{\textstyle\frac{\mathrm{Sc}_{g^{TX}}}{4}}\big)a_ju_j,u\Big\rangle_{S_\Omega} dv_\Omega\\
				&\leqslant \sum_{1\leqslant j\leqslant k} \int_\Omega\lV u\rV_{S_\Omega}^{2b^\ell-2}\mathrm{Re}\Big\langle\big(\lambda_j+{\textstyle\frac{\mathrm{Sc}^-_{g^{T\Omega}}}{4}}\big)a_ju_j,u\Big\rangle_{S_\Omega} dv_\Omega.
			\end{split}
		\end{equation}
		Using the Hölder inequality, we deduce that
		\begin{equation}\label{3.4'}
			\begin{split}
				&\sum_{1\leqslant j\leqslant k} \int_\Omega\lV u\rV_{S_\Omega}^{2b^\ell-2}\mathrm{Re}\Big\langle\big(\lambda_j+{\textstyle\frac{\mathrm{Sc}^-_{g^{T\Omega}}}{4}}\big)a_ju_j,u\Big\rangle_{S_\Omega} dv_\Omega\\
				&\leqslant\bigg(\int_\Omega\BV \sum_{1\leqslant j\leqslant k}\big(\lambda_j+{\textstyle\frac{\mathrm{Sc}^-_{g^{T\Omega}}}{4}}\big)a_ju_j\BV_{S_\Omega}^{2b^{\ell}}dv_\Omega\bigg)^{1/(2b^{\ell})}\bigg(\int_\Omega\lV u\rV_{S_\Omega}^{2b^\ell}dv_\Omega\bigg)^{(2b^\ell-1)/2b^\ell}.
			\end{split}
		\end{equation}
		
		To estimate the first term on the right hand side of \eqref{3.4'}, we observe that the function
		\begin{equation}
			F(\mu_1,\cdots,\mu_k)=	\int_\Omega\BV \sum_{1\leqslant j\leqslant k}\mu_ja_ju_j\BV_{S_\Omega}^{2b^{\ell}}dv_\Omega
		\end{equation}
		is convex as a function of $0\leqslant\mu_j\leqslant \lambda_k+\vert k_0\vert/4$ since
		\begin{equation}
			\frac{\pa^2 F}{\pa \mu_i^2}=2b^\ell(2b^\ell-1)\int_\Omega\vert a_i\vert^2\lV u_i\rV_{S_\Omega}^2\BV\sum_{1\leqslant j\leqslant k}\mu_ja_ju_j\BV_{S_\Omega}^{2b^{\ell}-2}dv_\Omega\geqslant 0.
		\end{equation}
		Therefore, there is $A\subseteq\{1,\cdots,k\}$ such that
		\begin{equation}\label{7.16}
			\int_\Omega\BV \sum_{1\leqslant j\leqslant k}(\lambda_j+\mathrm{Sc}_{g^{T\Omega}}^-)a_ju_j\BV_{S_\Omega}^{2b^{\ell}}dv_\Omega\leqslant  \bigg(\lambda_k+\frac{\vert k_0\vert}{4}\bigg)^{2b^\ell}	\int_\Omega\BV \sum_{j\in A}a_ju_j\BV_{S_\Omega}^{2b^{\ell}}dv_\Omega.
		\end{equation}
		The extremal property \eqref{2.11} of $u$ asserts that
		\begin{equation}
			\frac{\big(\int_\Omega\bV \sum_{j\in A}a_ju_j\bV_{S_\Omega}^{2b^{\ell}}dv_\Omega\big)^{1/(2b^\ell)}}{\big(\sum_{j\in A}\vert a_j\vert^2\big)^{1/2}}\leqslant \frac{\big(\int_\Omega\bV \sum_{0\leqslant j\leqslant k}a_ju_j\bV_{S_\Omega}^{2b^{\ell}}dv_\Omega\big)^{1/(2b^\ell)}}{\big(\sum_{0\leqslant j\leqslant k}\vert a_j\vert^2\big)^{1/2}},
		\end{equation}
		then we get
		\begin{equation}
			\begin{split}
				\int_\Omega\BV \sum_{j\in A}a_ju_j\BV_{S_\Omega}^{2b^{\ell}}dv_\Omega&\leqslant \Bigg(\frac{\sum_{j\in A}\vert a_j\vert^2}{\sum_{0\leqslant j\leqslant k}\vert a_j\vert^2}\Bigg)^{b^\ell}\int_\Omega\BV \sum_{0\leqslant j\leqslant k}a_ju_j\BV_{S_\Omega}^{2b^{\ell}}dv_\Omega\\
				&\leqslant\int_\Omega\BV \sum_{0\leqslant j\leqslant k}a_ju_j\BV_{S_\Omega}^{2b^{\ell}}dv_\Omega.
			\end{split}
		\end{equation}
		Plug this into \eqref{7.16}, we obtain
		\begin{equation}\label{3.5''}
			\begin{split}
				\int_\Omega\BV \sum_{1\leqslant j\leqslant k}(\lambda_j+\mathrm{Sc}_{g^{T\Omega}}^-)a_ju_j\BV_{S_\Omega}^{2b^{\ell}}dv_\Omega\leqslant \bigg(\lambda_k+\frac{\lv k_0\rv}{4}\bigg)^{2b^{\ell}}\int_\Omega\lV u\rV_{S_\Omega}^{2b^{\ell}}dv_\Omega.
			\end{split}
		\end{equation}
		
		Now put $f=\lV u\rV_{S_\Omega}$ and $a=2b^\ell$ in \eqref{2.7'}, then by \eqref{2.14'}, \eqref{3.4'} and \eqref{3.5''}, we get for a $C_m>0$,
		\begin{equation}
			\begin{split}
				\frac{C_m\mathrm{ID}(\Omega)^2}{2b^\ell}\lV u\rV_{L^{2b^{(\ell+1)}}}^{2b^{\ell}}\leqslant \bigg(\lambda_k+\frac{\lv k_0\rv}{4}\bigg)\lV u\rV_{L^{2b^{\ell}}}^{2b^{\ell}},
			\end{split}
		\end{equation}    
		or equivalently,
		\begin{equation}
			\begin{split}
				\lV u\rV_{L^{2b^{(\ell+1)}}}\leqslant	\Bigg(\frac{2b^\ell}{C_m\mathrm{ID}(\Omega)^2}\bigg(\lambda_k+\frac{\lv k_0\rv}{4}\bigg)\Bigg)^{1/(2b^\ell)}\lV u\rV_{L^{2b^{\ell}}},
			\end{split}
		\end{equation}    
		which implies \eqref{2.10}.
		
		Now let $\ell\to \infty$ in \eqref{2.10}, we see that
		\begin{equation}\label{3.8'}
			\begin{split}
				\lV u\rV_{L^{\infty}}&\leqslant	\prod_{\ell=0}^{\infty}	\Bigg(\frac{2b^\ell}{C_m\mathrm{ID}(\Omega)^2}\bigg(\lambda_k+\frac{\lv k_0\rv}{4}\bigg)\Bigg)^{1/(2b^\ell)}\lV u\rV_{L^{2}}\\
				&=\prod_{\ell=0}^{\infty}(2b^\ell)^{1/(2b^\ell)}\Bigg(\frac{\lambda_k+\frac{\lv k_0\rv}{4}}{C_m\mathrm{ID}(\Omega)^2}\Bigg)^{m/4}\lV u\rV_{L^{2}},
			\end{split}
		\end{equation}    
		which implies \eqref{3.2'}. \qed
	\end{pro}

	We shall imitate the proof of \cite[Lemma 7.2]{MR2962229} to get the following useful estimate.
	\begin{prop}
		There is $C_m>0$ such that for any $k\in\mathbb{N}$, we have
		\begin{equation}\label{3.10'}
			\sum_{i=1}^k\Vert u_i(x)\Vert^2_{S_\Omega}\leqslant\frac{C_m\big(\lambda_k+\vert k_0\vert\big)^{m/2}}{\mathrm{ID}(\Omega)^{m}}.
		\end{equation}
	\end{prop}
	
	\begin{pro}
		For any $x\in X$, let $\{v_j\}$ be an orthonormal basis of the fiber $S_{\Omega,x}$ at $x$, then we take $u_{x,v_j}\in C^\infty_0(\Omega,S_\Omega)$ by
		\begin{equation}
			u_{x,v_j}=\sum_{i=1}^k \langle v_j,u_i(x)\rangle_{S_\Omega}u_i.
		\end{equation}
		For any $u\in L^2(\Omega,S_\Omega)$, we denote its orthogonal projection to $\mathrm{span}\{u_1,u_2,\cdots,u_k\}$ by $p_{k}u$, that is,
		\begin{equation}
			p_{k}u=\sum_{i=1}^k\langle u, u_i\rangle_{L^2(\Omega,S_\Omega)}u_i,
		\end{equation}
		then we compute that
		\begin{equation}
			\begin{split}
				\langle u_{x,v_j},u\rangle_{L^2(\Omega,S_\Omega)}&=\sum_{i=1}^k \langle v_j,u_i(x)\rangle_{S_\Omega}\langle u_i,u\rangle_{L^2(\Omega,S_\Omega)}\\
				&=\big\langle v_j,(p_{k}u)(x)\big\rangle_{S_\Omega},
			\end{split}
		\end{equation}
		which, together with \eqref{3.2'}, implies that
		\begin{equation}
			\begin{split}
				\big\vert\langle u_{x,v_j},u\rangle_{L^2(\Omega,S_\Omega)}\big\vert&\leqslant \big\Vert p_{k}u\big\Vert_{C^0(\Omega,S_\Omega)}\\
				&\leqslant C_m\frac{\big(\lambda_k+\vert k_0\vert\big)^{m/4}}{\mathrm{ID}(\Omega)^{m/2}}\big\Vert p_{k}u\big\Vert_{L^2(\Omega,S_\Omega)}\\
				&\leqslant C_m\frac{\big(\lambda_k+\vert k_0\vert\big)^{m/4}}{\mathrm{ID}(\Omega)^{m/2}}\Vert u\Vert_{L^2(\Omega,S_\Omega)}.
			\end{split}
		\end{equation}
		Therefore, we have
		\begin{equation}
			\Vert u_{x,v_j}\Vert_{L^2(\Omega,S_\Omega)}\leqslant C_m\frac{\big(\lambda_k+\vert k_0\vert\big)^{m/4}}{\mathrm{ID}(\Omega)^{m/2}},
		\end{equation}
		in other words,
		\begin{equation}
			\Big(\sum_{i=1}^k\big\vert\langle v_j,u_i(x)\rangle_{S_\Omega}\big\vert^2\Big)^{1/2}\leqslant C_m\frac{\big(\lambda_k+\vert k_0\vert\big)^{m/4}}{\mathrm{ID}(\Omega)^{m/2}}.
		\end{equation}
		Squaring and summing over $j$, we deduce \eqref{3.10'}.  \qed
	\end{pro}

	By taking the integral of \eqref{3.10'}, we derive the following eigenvalue lower bound estimate.
	\begin{coro}
		There is $C_m>0$ such that for any $k\in\mathbb{N}$, we have
		\begin{equation}\label{3.10''}
			\lambda_k\geqslant C_mk^{2/m} \mathrm{ID}(\Omega)^{2}-\vert k_0\vert.
		\end{equation}
	\end{coro}

	\subsection{Functional calculus and estimates}\label{SEE3}

	For any $\phi\in\mathscr{S}(\mathbb{R})$, the operator $\phi\big(D^{S_\Omega}\big)$ has a smooth kernel $\phi\big(D^{S_\Omega}\big)(x,y)\in S_{X,x}\otimes S_{X,y}^*$ with respect to the volume $dv_\Omega$ such that for any $u\in C_0^\infty(M,S_\Omega)$ and $x,y\in M$,
	\begin{equation}
		\big(\phi\big(D^{S_\Omega}\big)u\big)(x)=\int_\Omega\phi\big(D^{S_\Omega}\big)(x,y)u(y)dv_\Omega(y).
	\end{equation}
	We define a norm $\lV\phi\rV_m$ of $\phi$ by
	\begin{equation}
		\lV\phi\rV_m=\sum_{\ell=0}^{\infty}\sup_{\vert t\vert\in[\ell,\ell+1]}\lv \phi(t)\rv(\ell^{m}+1).
	\end{equation}
	
	Now we state the main result of this section.
	\begin{theo}\label{T7.6}
		There exists $C_m>0$ such that for any $\phi\in\mathscr{S}(\mathbb{R})$, we have 
		\begin{equation}\label{3.19'}
			\Vert\phi\big(D^{S_\Omega}\big)(x,y)\Vert_{S_{X,x}\otimes S_{X,y}^*}\leqslant C_m\lV\phi\rV_m\frac{\vert k_0\vert^{m/2}+1}{\mathrm{ID}(\Omega)^{m}}.
		\end{equation}
	\end{theo}
	
	\begin{pro}
		First, we have
		\begin{equation}
			\begin{split}
				\phi\big(D^{S_\Omega}\big)(x,y)&=\sum_{k=0}^{\infty}\phi(\pm\lambda_k^{1/2})u_k(x)u_k(y)^*\\
				&=\sum_{\ell=0}^{\infty}\sum_{\ell \leqslant\lambda_k^{1/2}<\ell+1}\phi(\pm\lambda_k^{1/2})u_k(x)u_k(y)^*.
			\end{split}
		\end{equation}
		Then, by Cauchy-Schwarz inequality, we have
		\begin{equation}\label{3.21'}
			\begin{split}
				&\Vert\phi\big(D^{S_\Omega}\big)(x,y)\Vert_{S_{X,x}\otimes S_{X,y}^*}\\
				&\leqslant \sum_{\ell=0}^{\infty}\Big(\sup_{\vert t\vert\in[\ell,\ell+1]}\lv\phi(t)\rv\Big)\sum_{\ell \leqslant\lambda_k^{1/2}<\ell+1}\lV u_k(x)\rV_{S_\Omega}\lV u_k(y)\rV_{S_\Omega}\\
				&\leqslant \sum_{\ell=0}^{\infty}\Big(\sup_{\vert t\vert\in[\ell,\ell+1]}\lv\phi(t)\rv\Big)\Big(\sum_{\lambda_k^{1/2}<\ell+1}\lV u_k(x)\rV_{S_\Omega}^2\Big)^{1/2}\Big(\sum_{\lambda_k^{1/2}<\ell+1}\lV u_k(y)\rV_{S_\Omega}^2\Big)^{1/2},
			\end{split}
		\end{equation}
		which, together with \eqref{3.10'}, implies \eqref{3.19'}.\qed
	\end{pro}

	\begin{lemma}\label{L8.7}
		There exists $C_{\phi}>0$ such that for any $0<\varepsilon<\infty$, we have
		\begin{equation}\label{3.25}
			\Vert \phi_\varepsilon\Vert_{m}\leqslant C_{\phi}\Big(1+\frac{1}{\varepsilon^{m+1}}\Big),
		\end{equation}
		where $\phi_\varepsilon(t)=\phi(\varepsilon t)$
	\end{lemma}
	
	\begin{pro}
		We define $C_0=\sup_{t\in\mathbb{R}}\lv \phi(t)(1+t^2)\rv$ and compute that
		\begin{equation}\label{3.27}
			\begin{split}
				\sum_{\ell=0}^{\infty}\sup_{\vert t\vert\in[\ell,\ell+1]}\lv \phi_\varepsilon(t)\rv&\leqslant\sup_{t\in\mathbb{R}}\bv \phi(\varepsilon t)\bv+\sum_{\ell=1}^{\infty}\sup_{\vert t\vert\in[\ell,\ell+1]}\bv \phi(\varepsilon t)\bv\\
				&\leqslant C_0+\sum_{\ell=1}^{\infty}\frac{C_0}{1+(\varepsilon\ell)^2}\\
				&\leqslant C_0+\frac{C_0}{\varepsilon}\int_0^\infty\frac{1}{1+t^2}dt.
			\end{split}
		\end{equation}
		We have the obvious inequality
		\begin{equation}\label{3.28}
			\sup_{\vert t\vert\in[\ell,\ell+1]}\bv \phi_\varepsilon(t)\bv\ell^{m}\leqslant \sup_{\vert t\vert\in[\ell,\ell+1]}\bv g_\varepsilon(t)t^{m}\bv=\frac{1}{\varepsilon^{m}}\sup_{\vert t\vert\in[\ell,\ell+1]}\bv g(\varepsilon t)(\varepsilon t)^{m}\bv.
		\end{equation}
		Set $C_1=\sup_{t\in\mathbb{R}}\lv g(t)t^{m}(1+t^2)\rv$, then from \eqref{3.28}, we can repeat \eqref{3.27} by replacing $\phi$ with $t^{m}\phi$ to get
		\begin{equation}\label{3.29}
			\begin{split}
				\sum_{\ell=0}^{\infty}	\sup_{\vert t\vert\in[\ell,\ell+1]}\bv g_\varepsilon(t)\bv\ell^{m}\leqslant\frac{1}{\varepsilon^{m}}\Big(C_1+\frac{C_1}{\varepsilon}\int_0^\infty\frac{1}{1+t^2}dt\Big).
			\end{split}
		\end{equation} 
		By \eqref{3.27} and \eqref{3.29}, we deduce \eqref{3.25}.\qed
	\end{pro}

	\addcontentsline{toc}{section}{References}
	
	\bibliographystyle{abbrv}


\begin{thebibliography}{10}
		
		\bibitem{MR2273508}
		N.~Berline, E.~Getzler, and M.~Vergne.
		\newblock {\em Heat kernels and {D}irac operators}.
		\newblock Grundlehren Text Editions. Springer-Verlag, Berlin, 2004.
		\newblock Corrected reprint of the 1992 original.
		
		\bibitem{MR813584}
		J.-M. Bismut.
		\newblock The {A}tiyah-{S}inger index theorem for families of {D}irac
		operators: two heat equation proofs.
		\newblock {\em Invent. Math.}, 83(1):91--151, 1986.
		
		\bibitem{MR1188532}
		J.-M. Bismut and G.~Lebeau.
		\newblock Complex immersions and {Q}uillen metrics.
		\newblock {\em Inst. Hautes \'{E}tudes Sci. Publ. Math.}, (74):ii+298 pp.
		(1992), 1991.
		
		\bibitem{MR4386411}
		S.~Braun and R.~Sauer.
		\newblock Volume and macroscopic scalar curvature.
		\newblock {\em Geom. Funct. Anal.}, 31(6):1321--1376, 2021.
		
		\bibitem{MR263092}
		J.~Cheeger.
		\newblock Finiteness theorems for {R}iemannian manifolds.
		\newblock {\em Amer. J. Math.}, 92:61--74, 1970.
		
		\bibitem{MR823176}
		A.~Connes.
		\newblock Noncommutative differential geometry.
		\newblock {\em Inst. Hautes \'Etudes Sci. Publ. Math.}, (62):257--360, 1985.
		
		\bibitem{MR866491}
		A.~Connes.
		\newblock Cyclic cohomology and the transverse fundamental class of a
		foliation.
		\newblock In {\em Geometric methods in operator algebras ({K}yoto, 1983)},
		volume 123 of {\em Pitman Res. Notes Math. Ser.}, pages 52--144. Longman Sci.
		Tech., Harlow, 1986.
		
		\bibitem{MR1066176}
		A.~Connes and H.~Moscovici.
		\newblock Cyclic cohomology, the {N}ovikov conjecture and hyperbolic groups.
		\newblock {\em Topology}, 29(3):345--388, 1990.
		
		\bibitem{MR608287}
		C.~B. Croke.
		\newblock Some isoperimetric inequalities and eigenvalue estimates.
		\newblock {\em Ann. Sci. \'Ecole Norm. Sup. (4)}, 13(4):419--435, 1980.
		
		\bibitem{MR836727}
		E.~Getzler.
		\newblock A short proof of the local {A}tiyah-{S}inger index theorem.
		\newblock {\em Topology}, 25(1):111--117, 1986.
		
		\bibitem{MR686042}
		M.~Gromov.
		\newblock Volume and bounded cohomology.
		\newblock {\em Inst. Hautes \'Etudes Sci. Publ. Math.}, (56):5--99, 1982.
		
		\bibitem{MR859578}
		M.~Gromov.
		\newblock Large {R}iemannian manifolds.
		\newblock In {\em Curvature and topology of {R}iemannian manifolds ({K}atata,
			1985)}, volume 1201 of {\em Lecture Notes in Math.}, pages 108--121.
		Springer, Berlin, 1986.
		
		\bibitem{MR3201312}
		M.~Gromov.
		\newblock Dirac and {P}lateau billiards in domains with corners.
		\newblock {\em Cent. Eur. J. Math.}, 12(8):1109--1156, 2014.
		
		\bibitem{MR4577903}
		M.~Gromov.
		\newblock Four lectures on scalar curvature.
		\newblock In {\em Perspectives in scalar curvature. {V}ol. 1}, pages 1--514.
		World Sci. Publ., Hackensack, NJ, [2023] \copyright 2023.
		
		\bibitem{MR577131}
		M.~Gromov and H.~B. Lawson.
		\newblock The classification of simply connected manifolds of positive scalar
		curvature.
		\newblock {\em Ann. of Math. (2)}, 111(3):423--434, 1980.
		
		\bibitem{MR569070}
		M.~Gromov and H.~B. Lawson.
		\newblock Spin and scalar curvature in the presence of a fundamental group.
		{I}.
		\newblock {\em Ann. of Math. (2)}, 111(2):209--230, 1980.
		
		\bibitem{MR1029396}
		K.~Grove, P.~Petersen, V, and J.~Y. Wu.
		\newblock Geometric finiteness theorems via controlled topology.
		\newblock {\em Invent. Math.}, 99(1):205--213, 1990.
		
		\bibitem{MR2827817}
		L.~Guth.
		\newblock Metaphors in systolic geometry.
		\newblock In {\em Proceedings of the {I}nternational {C}ongress of
			{M}athematicians. {V}olume {II}}, pages 745--768. Hindustan Book Agency, New
		Delhi, 2010.
		
		\bibitem{MR2753599}
		L.~Guth.
		\newblock Volumes of balls in large {R}iemannian manifolds.
		\newblock {\em Ann. of Math. (2)}, 173(1):51--76, 2011.
		
		\bibitem{MR1867354}
		A.~Hatcher.
		\newblock {\em Algebraic topology}.
		\newblock Cambridge University Press, Cambridge, 2002.
		
		\bibitem{MR602436}
		H.~Hess, R.~Schrader, and D.~A. Uhlenbrock.
		\newblock Kato's inequality and the spectral distribution of {L}aplacians on
		compact {R}iemannian manifolds.
		\newblock {\em J. Differential Geometry}, 15(1):27--37, 1980.
		
		\bibitem{MR358873}
		N.~Hitchin.
		\newblock Harmonic spinors.
		\newblock {\em Advances in Math.}, 14:1--55, 1974.
		
		\bibitem{MR375153}
		J.~L. Kazdan and F.~W. Warner.
		\newblock Existence and conformal deformation of metrics with prescribed
		{G}aussian and scalar curvatures.
		\newblock {\em Ann. of Math. (2)}, 101:317--331, 1975.
		
		\bibitem{MR2962229}
		P.~Li.
		\newblock {\em Geometric analysis}, volume 134 of {\em Cambridge Studies in
			Advanced Mathematics}.
		\newblock Cambridge University Press, Cambridge, 2012.
		
		\bibitem{MR156292}
		A.~Lichnerowicz.
		\newblock Spineurs harmoniques.
		\newblock {\em C. R. Acad. Sci. Paris}, 257:7--9, 1963.
		
		\bibitem{MR1850748}
		K.~Liu and W.~Zhang.
		\newblock Adiabatic limits and foliations.
		\newblock In {\em Topology, geometry, and algebra: interactions and new
			directions ({S}tanford, {CA}, 1999)}, volume 279 of {\em Contemp. Math.},
		pages 195--208. Amer. Math. Soc., Providence, RI, 2001.
		
		\bibitem{MR2206643}
		C.~L\"oh.
		\newblock Measure homology and singular homology are isometrically isomorphic.
		\newblock {\em Math. Z.}, 253(1):197--218, 2006.
		
		\bibitem{MR4665497}
		Q.~Ma.
		\newblock Toeplitz operators and the full asymptotic torsion forms.
		\newblock {\em J. Funct. Anal.}, 286(3):Paper No. 110210, 74pp, 2024.
		
		\bibitem{ma2024boundingahatgenususing}
		Q.~Ma, J.~Wang, G.~Yu, and B.~Zhu.
		\newblock Bounding the a-hat genus using scalar curvature lower bounds and
		isoperimetric constants.
		\newblock {\em arXiv:2410.08866}, 2024.
		
		\bibitem{MR2339952}
		X.~Ma and G.~Marinescu.
		\newblock {\em Holomorphic {M}orse inequalities and {B}ergman kernels}, volume
		254 of {\em Progress in Mathematics}.
		\newblock Birkh\"{a}user Verlag, Basel, 2007.
		
		\bibitem{MR1292020}
		H.~Moscovici and F.~Wu.
		\newblock Index theory without symbols.
		\newblock In {\em {$C^\ast$}-algebras: 1943--1993 ({S}an {A}ntonio, {TX},
			1993)}, volume 167 of {\em Contemp. Math.}, pages 304--351. Amer. Math. Soc.,
		Providence, RI, 1994.
		
		\bibitem{MR1254310}
		H.~Moscovici and F.-B. Wu.
		\newblock Localization of topological {P}ontryagin classes via finite
		propagation speed.
		\newblock {\em Geom. Funct. Anal.}, 4(1):52--92, 1994.
		
		\bibitem{MR159138}
		J.~Moser.
		\newblock On {H}arnack's theorem for elliptic differential equations.
		\newblock {\em Comm. Pure Appl. Math.}, 14:577--591, 1961.
		
		\bibitem{MR996446}
		J.~Roe.
		\newblock Partitioning noncompact manifolds and the dual {T}oeplitz problem.
		\newblock In {\em Operator algebras and applications, {V}ol.\ 1}, volume 135 of
		{\em London Math. Soc. Lecture Note Ser.}, pages 187--228. Cambridge Univ.
		Press, Cambridge, 1988.
		
		\bibitem{MR1147350}
		J.~Roe.
		\newblock Coarse cohomology and index theory on complete {R}iemannian
		manifolds.
		\newblock {\em Mem. Amer. Math. Soc.}, 104(497):x+90, 1993.
		
		\bibitem{MR720934}
		J.~Rosenberg.
		\newblock {$C\sp{\ast} $}-algebras, positive scalar curvature, and the
		{N}ovikov conjecture.
		\newblock {\em Inst. Hautes \'Etudes Sci. Publ. Math.}, (58):197--212, 1983.
		
		\bibitem{MR1157815}
		W.~Rudin.
		\newblock {\em Functional analysis}.
		\newblock International Series in Pure and Applied Mathematics. McGraw-Hill,
		Inc., New York, second edition, 1991.
		
		\bibitem{MR541332}
		R.~Schoen and S.-T. Yau.
		\newblock Existence of incompressible minimal surfaces and the topology of
		three-dimensional manifolds with nonnegative scalar curvature.
		\newblock {\em Ann. of Math. (2)}, 110(1):127--142, 1979.
		
		\bibitem{MR526976}
		R.~Schoen and S.-T. Yau.
		\newblock On the proof of the positive mass conjecture in general relativity.
		\newblock {\em Comm. Math. Phys.}, 65(1):45--76, 1979.
		
		\bibitem{MR4554426}
		W.~P. Thurston.
		\newblock {\em The geometry and topology of three-manifolds. {V}ol. {IV}}.
		\newblock American Mathematical Society, Providence, RI, [2022] \copyright
		2022.
		\newblock Edited and with a preface by Steven P. Kerckhoff and a chapter by J.
		W. Milnor.
		
		\bibitem{MR4411373}
		R.~Willett and G.~Yu.
		\newblock {\em Higher index theory}, volume 189 of {\em Cambridge Studies in
			Advanced Mathematics}.
		\newblock Cambridge University Press, Cambridge, 2020.
		
		\bibitem{MR626707}
		E.~Witten.
		\newblock A new proof of the positive energy theorem.
		\newblock {\em Comm. Math. Phys.}, 80(3):381--402, 1981.
		
		\bibitem{MR1626745}
		G.~Yu.
		\newblock The {N}ovikov conjecture for groups with finite asymptotic dimension.
		\newblock {\em Ann. of Math. (2)}, 147(2):325--355, 1998.
		
		\bibitem{MR3590530}
		G.~Yu.
		\newblock The {N}ovikov conjecture for algebraic {K}-theory of the group
		algebra over the ring of {S}chatten class operators.
		\newblock {\em Adv. Math.}, 307:727--753, 2017.
		
		\bibitem{MR3664818}
		W.~Zhang.
		\newblock Positive scalar curvature on foliations.
		\newblock {\em Ann. of Math. (2)}, 185(3):1035--1068, 2017.
		
	\end{thebibliography}

	\def\cprime{$'$} \def\cprime{$'$}

\end{document}